\documentclass[12pt]{article}

\usepackage{amsmath}
\usepackage{amsmath, amsthm, amsfonts, amssymb, color}
\usepackage{mathrsfs}
\usepackage{latexsym,bm}
\usepackage{graphicx}

\newtheorem{theorem}{Theorem}[section]
\newtheorem{proposition}[theorem]{Proposition}

\newtheorem{lemma}[theorem]{Lemma}
\newtheorem{corollary}[theorem]{Corollary}

\newtheorem{remark}[theorem]{Remark}
\newtheorem{example}[theorem]{Example}
\newtheorem{defi}[theorem]{Definition}

\def\eps{\varepsilon}

\newcommand{\me}{\mathcal{M}(E)}

\newcommand{\p}{\mathrm{P}}
\newcommand{\pp}{\mathbb{P}}

\newcommand{\q}{\mathrm{Q}}

\newcommand{\B}{\mathcal{B}}
\newcommand{\R}{\mathbb{R}}

\newcommand{\esup}{\mathrm{essup}}

\definecolor{wco}{rgb}{0.5,0.2,0.3}

\numberwithin{equation}{section} 

\begin{document}

\allowdisplaybreaks

\title{\bf Spine decomposition and $L\log L$ criterion for  superprocesses with non-local branching mechanisms
}
\author{ \bf  Yan-Xia Ren\footnote{LMAM School of Mathematical Sciences \& Center for
Statistical Science, Peking
University, Beijing, 100871, P.R. China.
Email: yxren@math.pku.edu.cn}
\hspace{1mm}\hspace{1mm}\hspace{1mm}
Renming Song\footnote{Department of Mathematics, University of Illinois, Urbana, IL 61801, U.S.A.
Email: rsong@illinois.edu} \hspace{1mm}\hspace{1mm} and \hspace{1mm}\hspace{1mm}
Ting Yang\footnote{School of Mathematics and Statistics, Beijing Institute of Technology, Beijing, 100081, P.R.China.
Email: yangt@bit.edu.cn}
}
\date{}
\maketitle

\begin{abstract}
In this paper, we provide a pathwise spine decomposition for superprocesses with both local and non-local branching mechanisms under a martingale change of measure.
This result complements earlier results established for superprocesses with purely local
 branching mechanisms and for multitype superprocesses.
As an application of this decomposition, we obtain necessary/sufficient conditions for the limit of the fundamental martingale to be non-degenerate. In particular, we obtain extinction properties of
superprocesses with non-local branching mechanisms
as well as a Kesten-Stigum $L\log L$ theorem for the fundamental martingale.
\end{abstract}

\medskip


\medskip


\section{Introduction}
The so-called \textit{spine decomposition} for superprocesses was introduced in terms of a semigroup decomposition by Evans
\cite{E}. To be more specific,  Evans \cite{E} described the semigroup of a superprocess with branching mechanism $\psi(\lambda)=\lambda^{2}$ under a martingale change of measure in terms of the semigroup of an immortal particle (called the \textit{spine}) and the semigroup of the original superprocess. Since then there has been a lot of interest in finding the spine decomposition for other types
of superprocesses due to a variety of applications. For example, Engl\"{a}nder and Kyprianou \cite{EK} used a similar semigroup decomposition to establish the $L^{1}$-convergence of martingales for superdiffusions with quadratic branching mechanisms. Later, Kyprianou et al. \cite{KLMR,KM} obtained a pathwise spine decomposition for a one-dimensional super-Brownian motion with spatially-independent local branching mechanism, in which independent copies of the original superprocess immigrate along the path of the immortal particle, and they used this decomposition to establish the $L^{p}$-boundedness ($p\in (1,2]$) of martingales.  A similar pathwise decomposition was obtained by Liu et al. \cite{LRS} for a class of superdiffusions in bounded domains with spatially-dependent local branching mechanisms, and it was used to establish a Kesten-Stigum $L\log L$ theorem, which gives the sufficient and necessary condition for the martingale limit to be non-degenerate.
In the set-up of branching Markov processes, such as branching diffusions and branching random walks,
an analogous decomposition has been introduced and used as a tool
to analyze branching Markov processes.
See, for example, \cite{HH} for a brief history of the spine
approach for branching Markov processes.
Until very recently such a spine decomposition for superprocesses was only available for superprocesses with local branching mechanisms.
In the recent paper \cite{KP}, Kyprianou and Palau established a spine decomposition
for a multitype continuous-state branching process (MCSBP) and used it to study the extinction properties. Concurrently to their work, a similar decomposition has been obtained by Chen et al. \cite{CRS} for a multitype superdiffusion.
However, in both papers,
only a very special kind of non-local branching mechanisms
are considered.
The first goal of
this paper is to close the gap by establishing a pathwise spine decomposition for superprocesses with
both local and general  non-local branching mechanisms.

In this paper, the Schr\"{o}dinger operator associated with the mean semigroup of the superprocess is characterised by its bilinear form.
Then some technical assumptions (Assumptions 1-2 below) are imposed to ensure the existence of a principal eigenvalue $\lambda_{1}$
and a positive ground state $h$, and hence to ensure the existence of a fundamental martingale (Theorem \ref{them1} below). These assumptions may look strong, but they hold for a large class of processes, and we illustrate this for several interesting examples, including MCSBP, in Section \ref{s:examples}.
Our result (Theorem \ref{them2} below)
shows that, for a superprocess with both local and non-local branching, under a martingale change of measure, the spine runs as a copy of a conservative process, which can be constructed by concatenating copies of
a subprocess of the $h$-transform of the original spatial motion via a transfer
kernel determined by the non-local branching mechanism,
and the general nature of the branching mechanism induces three different kinds of immigration: the \textit{continuous}, \textit{discontinuous} and \textit{revival-caused} immigration. The concatenating procedure and revival-caused immigration are consequences of non-local branching, and they do not occur when the branching mechanism is purely local.

In connection with the limit theory, it is natural to ask whether or not the limit of the fundament martingale is non-degenerate.
Using the spine decomposition, we establish sufficient and necessary conditions
for the martingale limit to be non-degenerate, respectively, in Theorem \ref{them3} and Theorem \ref{them4}.
A similar idea was used in \cite{EK,LRS} for (purely local branching) superdiffusions, and in \cite{KP} for MCSBP.
However, in this paper, we extend this idea much further by considering superprocesses where the
spatial motion may be discontinuous and the branching mechanism is allowed to
be generally non-local.
 Suppose that $\{Z_{n}:n\ge 1\}$ is a Galton-Watson branching process with
offspring distribution $\{p_n: n\ge 0\}$.
 Let $L$ stand for a random variable with this offspring distribution. Let
 $m:=\sum_{n=0}^{+\infty}np_{n}$ be the
 mean of the offspring distribution.
 Then $Z_{n}/m^{n}$ is a non-negative martingale.
Kesten and Stigum \cite{KeS} proved that
 when $1<m<+\infty$, the limit of $Z_{n}/m^{n}$ is non-degenerate if and only if $\mathrm{E}\left(L\log^{+}L\right)<+\infty$.
This result is usually referred to the Kesten-Stigum $L\log L$ theorem.
Our Corollary
\ref{cor2} shows that, in the case of $\lambda_{1}<0$, the martingale limit is non-degenerate if and only if
an $L\log L$-type condition holds.
This result extends an earlier result obtained in \cite{LRS} for superdiffusions and can be viewed as a natural analogue of the Kesten-Stigum $L\log L$ theorem for superprocesses.
Our Corollary \ref{cor8.3} says that, under suitable assumptions, the non-local branching superprocess exhibits weak local extinction if and only if $\lambda_{1}\ge 0$.
This result can be regarded as a general non-local branching
counterpart of \cite[Theorem 3]{EK},
where the same result is proved for a special class of superdiffusions in a domain $D\subset\R^d$
(the branching mechanism considered in \cite{EK} is $\psi(x,\lambda)=-\beta(x)\lambda+\alpha(x)\lambda^{2}$ with $\alpha, \beta$ being H\"{o}lder continuous functions in $D$ with order $\eta\in(0,1]$, $\alpha>0$  in $D$ and $\beta$ being bounded from above.

In this paper we assume the spatial motion to be a
symmetric Hunt process on a locally compact separable metric space.
This assumption is not really necessary.
An extension is possible. One direction is to assume the spatial motion to be a transient Borel right process on a Luzin space, whose Dirichlet form satisfies Silverstein's sector condition.
Definitions of smooth measures and Kato class can then be extended, while still preserving the properties used in this paper. We refer to \cite{Chen,CS} for Kato class measures defined in this way.
Nevertheless, we keep to the less general type of spatial motions to avoid unnecessary
technicalities.

The rest of this paper is organized as follows.
In Section \ref{Sec2} we review some
basic definitions and properties of non-local branching superprocesses,
including the definition of Kuznetsov measures which will be used later.
In Section \ref{Sec3}, we present our main working assumptions and the fundamental martingale.
Section \ref{Sec4} provides the spine decomposition and its proof.
The proof of Proposition \ref{prop3} is postponed to the Appendix.
In Sections \ref{Sec5} and \ref{Sec6} we use the spine decomposition to find sufficient and necessary conditions
for the limit of the fundamental martingale to be non-degenerate respectively. In particular, we obtain extinction properties of the non-local branching superprocess as well as a Kesten-Stigum $L\log L$ theorem for the martingale.
In the last section, we give some concrete examples to illustrate our results.

Notation and basic setting: Throughout this paper we use ``$:=$" as a definition.
We always assume that $E$ is a
locally compact separable metric space
with Borel $\sigma$-algebra $\mathcal{B}(E)$ and $m$ is a $\sigma$-finite measure on $(E,\mathcal{B}(E))$ with full support.
 Let $E_{\partial}:=E\cup \{\partial\}$ be the one-point compactification of $E$.
Any function $f$ on $E$
will be automatically extended to $E_{\partial}$ by
setting $f(\partial)=0$.  For a function $f$ on $E$, $\|f\|_{\infty}:=\sup_{x\in E}|f(x)|$ and $\esup_{x\in E}f:=\inf_{N:m(N)=0}\sup_{x\in E\setminus N}|f(x)|$.
Numerical functions $f$ and $g$ on $E$ are said to be $m$-equivalent ($f=g\ [m]$ in notation) if $m\left(\{x\in E:f(x)\not=g(x)\}\right)=0$. If $f(x,t)$ is a function on $E\times [0,+\infty)$, we say $f$ is \textit{locally bounded} if $\sup_{t\in [0,T]}\sup_{x\in E}|f(x,t)|<+\infty$ for every $T\in (0,+\infty)$. We denote by $f^{t}(\cdot)$ the function $x\mapsto f(x,t)$.
Let $\me$ denote the space of finite Borel measures on $E$ topologized by the weak convergence.
Let $\me^{0}:=\me\setminus\{0\}$ where $0$ denotes the null measure on $E$. When $\mu$ is a measure on $\mathcal{B}(E)$ and $f$, $g$ are measurable functions, let $\langle f,\mu\rangle:=\int_{E}f(x)\mu({\rm d}x)$ and $(f,g):=\int_{E}f(x)g(x)m({\rm d}x)$ whenever the right hand sides make sense. Sometimes we also write $\mu(f)$ for $\langle f,\mu\rangle$. We use $\mathcal{B}_{b}(E)$ (respectively,
$\mathcal{B}^{+}(E)$) to denote
the space of bounded (respectively,
non-negative) measurable functions on $(E,\mathcal{B}(E))$.
For $a,b\in \mathbb{R}$,
$a\wedge b:=\min\{a,b\}$,
$a\vee b:=\max\{a,b\}$,
and $\log^{+}a:=\log (a\vee 1)$.

\section{Preliminaries}\label{Sec2}

\subsection{Superprocess with non-local branching mechanisms}\label{sec1}
Let $\xi=(\Omega,\mathcal{H},\mathcal{H}_{t},\theta_{t},\xi_{t},\Pi_{x},\zeta)$ be an
$m$-symmetric Hunt process
on $E$. Here $\{\mathcal{H}_{t}:\ t\ge 0\}$ is the minimal admissible filtration, $\{\theta_{t}:\ t\ge 0\}$ the time-shift operator of $\xi$ satisfying $\xi_{t}\circ\theta_{s}=\xi_{t+s}$ for $s,t\ge 0$, and $\zeta:=\inf\{t>0:\
\xi_{t}=\partial\}$ the lifetime of $\xi$.  Let $\{\mathfrak{S}_{t}:\ t\ge 0\}$ be the transition semigroup of $\xi$, i.e., for any non-negative measurable function $f$,
 $$\mathfrak{S}_{t}f(x):=\Pi_{x}\left[f(\xi_{t})\right].$$
For $\alpha>0$ and $f\in\mathcal{B}^{+}(E)$,
let $G_{\alpha}f(x):=\int_{0}^{+\infty}e^{-\alpha t}\mathfrak{S}_{t}f(x){\rm d}t$.
It is known by \cite[Lemma 1.1.14]{CF}
that $\{\mathfrak{S}_{t}:t\ge 0\}$ can be uniquely extended to a strongly continuous contraction semigroup on $L^{2}(E,m)$, which we also denote by $\{\mathfrak{S}_{t}:t\ge 0\}$.
By the theory of Dirichlet forms, there exists a
regular symmetric Dirichlet form
$(\mathcal{E},\mathcal{F})$ on $L^{2}(E,m)$ associated with $\xi$:
 $$
 \mathcal{F}= \Big \{u\in L^{2}(E,m):\ \sup_{t>0}\frac{1}{t}\int_{E}\left(u(x)-\mathfrak{S}_{t}u(x)\right)u(x)m({\rm d}x)<+\infty \Big\},
 $$
 $$\mathcal{E}(u,v)=\lim_{t\to 0}\frac{1}{t}\int_{E}\left(u(x)-\mathfrak{S}_{t}u(x)\right)v(x)m({\rm d}x),\quad
 \forall u,v\in\mathcal{F}.
 $$
 Moreover, for all $f\in \mathcal{B}_{b}(E)\cap L^{2}(E,m)$ and $\alpha>0$,
 \begin{equation}
 G_{\alpha}f\in\mathcal{F}\mbox{ satisfies that }\mathcal{E}_{\alpha}(G_{\alpha}f,v)=(f,v)\quad\forall v\in \mathcal{F},\label{d1}
 \end{equation}
 where $\mathcal{E}_{\alpha}(u,v):=\mathcal{E}(u,v)+\alpha(u,v)$.
We assume that $\xi$ admits a transition density $p(t,x,y)$ with respect to the measure $m$, which is symmetric
in $(x,y)$ for each $t>0$.
Under this absolute continuity assumption, ``quasi everywhere" statements can be strengthened to ``everywhere" ones.  Moreover, we can define notions without
exceptional sets, for example, {\it positive
continuous additive functionals (PCAF in abbreviation) in the strict sense} (cf. \cite[Section 5.1]{FOT}).
In this paper, we will only deal with notions in the strict sense and
omit ``in the strict sense''.

It is well known (see \cite[Theorem A.3.21]{FOT}, for instance) there exist a kernel $N(x, dy)$
on $(E, \mathcal{B}(E))$ with $N(x, \{x\})=0$ for all $x\in E$ and a PCAF $H$ of $\xi$ with
$\int_E\Pi_x(H_t)\mu(dx)<\infty$ for all $t\ge 0$ and probability measure $\mu$ on
 $(E, \mathcal{B}(E))$ such that for any $x\in E$, any $t\ge0$,
and any non-negative Borel function $f$ on $E\times E$ vanishing on the diagonal $\{(y, y): y\in E\}$,
\begin{equation}
\Pi_{x}\left[\sum_{s\le t}f(\xi_{s-},\xi_{s})\right]=\Pi_{x}\left[\int_{0}^{t}
\int_{E}
f(\xi_{s},y)N(\xi_{s},{\rm d}y){\rm d}H_s\right].\label{levysys}
\end{equation}
The pair $(N, H)$ is called a L\'evy system of $\xi$.

In this paper, we consider a superprocess $X:=\{X_{t}:\ t\ge 0\}$ with spatial motion $\xi$ and a non-local branching mechanism $\psi$ given by
\begin{equation}
\psi(x,f)=\phi^{L}(x,f(x))+\phi^{NL}(x,f)\quad\mbox{ for }x\in E,\ f\in\mathcal{B}^{+}_{b}(E).\label{1.1}
\end{equation}
The first term $\phi^{L}$ in \eqref{1.1}
is called  the local branching mechanism and takes the form
\begin{equation}\label{eq:local part}
\phi^{L}(x,\lambda)=a(x)\lambda+b(x)\lambda^{2}+\int_{(0,+\infty)}\left(e^{-\lambda \theta}-1+\lambda \theta\right)\Pi^{L}(x,{\rm d}\theta),
\quad x\in E,\ \lambda\ge 0,
\end{equation}
where $a\in\mathcal{B}_{b}(E)$, $b\in\mathcal{B}^{+}_{b}(E)$
and $(\theta\wedge \theta^{2})\Pi^{L}(x,{\rm d}\theta)$ is a bounded kernel from $E$ to $(0,+\infty)$. The second term $\phi^{NL}$ in \eqref{1.1}
is called  the non-local branching mechanism and takes the form
\begin{equation}\label{eq:non-local part}
\phi^{NL}(x,f)=-c(x)\pi(x,f)-\int_{(0,+\infty)}\left(1-e^{-\theta\pi(x,f)}\right)
\Pi^{NL}(x,{\rm d}\theta), \quad x\in E,
\end{equation}
where $c(x)$ is a non-negative
bounded measurable function on $E$,  $\pi(x,{\rm d}y)$ is a probability kernel on $E$
with $\pi(x,\{x\})\not\equiv 1$, $\pi(x,f)$ stands for $\int f(y)\pi(x,{\rm d}y)$
and $\theta\Pi^{NL}(x,{\rm d}\theta)$ is a bounded kernel from $E$ to $(0,+\infty)$.
To be specific, $X$ is an $\me$-valued Markov process such that for every $f\in\mathcal{B}^{+}_{b}(E)$ and every $\mu\in\me$,
\begin{equation}
\p_{\mu}\left(e^{-\langle f,X_{t}\rangle}\right)=e^{-\langle u_{f}(\cdot,t),\mu\rangle}\quad\mbox{ for }t\ge 0,\label{1.4}
\end{equation}
where $u_{f}(x,t):=-\log\p_{\delta_{x}}\left(e^{-\langle f,X_{t}\rangle}\right)$ is the unique non-negative locally bounded solution to the integral equation
\begin{align}
u_{f}(x,t)
&=\mathfrak{S}_{t}f(x)-\Pi_{x}\left[\int_{0}^{t}\psi(\xi_{s},u^{t-s}_{f}){\rm d}s\right]\nonumber\\
&=\mathfrak{S}_{t}f(x)-\Pi_{x}\left[\int_{0}^{t}
\phi^{L}
(\xi_{s},u_{f}(t-s,\xi_{s})){\rm d}s\right]
-\Pi_{x}\left[\int_{0}^{t}
\phi^{NL}
(\xi_{s},u^{t-s}_{f}){\rm d}s\right].\label{1.3}
\end{align}
We refer to the process described
above as a $(\mathfrak{S}_{t},\phi^{L},\phi^{NL})$-superprocess.
Such a process is defined in \cite{Li} via its log-Laplace functional. Another usual way of constructing the $(\mathfrak{S}_{t},\phi^{L},\phi^{NL})$-superprocess is as the high
intensity limit of
a sequence of branching particle systems, where whenever a particle dies, it chooses from two different branching types: the local branching type
(when the particle dies at $x$, it is replaced by a random number of offspring situated at $x$), and the non-local branching type (when the particle dies, it gives birth to a random number of particles in $E$, and the offspring then start to move from their locations of birth). We refer to \cite{Li,DGL} for such a construction.

We define for $x\in E$,
\begin{equation}\label{def:gamma-rs}
\gamma(x,{\rm d}y):=\left(c(x)+\int_{(0,+\infty)}\theta\Pi^{NL}(x,{\rm d}\theta)\right)\pi(x,{\rm d}y),\quad\gamma(x):=\gamma(x,1).
\end{equation}
Clearly, $\gamma(x)$ is a non-negative bounded function on $E$ and $\gamma(x,{\rm d}y)$ is a bounded kernel on $E$.
Define
$A:=\{x\in E:\ \gamma(x)>0\}.$
Note that $\phi^{NL}(x,\cdot)=0$ for all $x\in E\setminus A$.
If $A=\emptyset$ (i.e., $\phi^{NL}\equiv 0$), we call $\psi$ a (purely) local branching mechanism. Without loss of generality, we always assume that $A\not=\emptyset$. The arguments and results of this paper also work for (purely) local branching mechanisms.

It follows from \cite[Theorem 5.12]{Li} that the $(\mathfrak{S}_{t},\phi^{L},\phi^{NL})$-superprocess has a right realization in $\me$. Let $\mathcal{W}^{+}_{0}$ denote the space of right continuous paths from $[0,+\infty)$ to $\me$ having zero as a trap. We may and do assume that $X$ is the coordinate process in $\mathcal{W}^{+}_{0}$ and that $(\mathcal{F}_{\infty},(\mathcal{F}_{t})_{t\ge 0})$ is the natural
filtration
on $\mathcal{W}^{+}_{0}$ generated by the coordinate process. The following proposition follows from \cite[Proposition 2.27 and Proposition 2.29]{Li}.
\begin{proposition}\label{prop0}
For all $\mu\in\me$ and $f\in\mathcal{B}_{b}(E)$,
\begin{equation*}
\p_{\mu}\left(\langle f,X_{t}\rangle\right)=
\langle \mathfrak{P}_{t}f,\mu\rangle,
\end{equation*}
where $\mathfrak{P}_{t}f(x)$ is the unique locally
bounded solution to the following integral equation:
\begin{align}
\mathfrak{P}_{t}f(x)&=\mathfrak{S}_{t}f(x)
-\Pi_{x}\left[\int_{0}^{t}a(\xi_{s})
\mathfrak{P}_{t-s}f(\xi_{s}){\rm d}s\right]
+\Pi_{x}\left[\int_{0}^{t}\gamma(\xi_{s},\mathfrak{P}_{t-s}f){\rm d}s\right].
\label{li2}
\end{align}
Moreover, for all $\mu\in\me$, $g\in\mathcal{B}^{+}_{b}(E)$ and $f\in\mathcal{B}_{b}(E)$,
\begin{equation}
\p_{\mu}\left(\langle f,X_{t}\rangle e^{-\langle g,X_{t}\rangle }\right)=e^{-\langle V_{t}g,\mu\rangle}\langle V^{f}_{t}g,\mu\rangle,\nonumber
\end{equation}
where $V_{t}g(x):=u_{g}(x,t)$ is the unique non-negative locally bounded solution to \eqref{1.3} with initial value $g$,
and $V^{f}_{t}g(x)$ is the unique locally bounded solution to the following integral equation
\begin{equation}
V^{f}_{t}g(x)=\mathfrak{S}_{t}f(x)-\Pi_{x}\left[\int_{0}^{t}\Psi(\xi_{s},V_{t-s}g,V^{f}_{t-s}g){\rm d}s\right],\label{li1}
\end{equation}
where
\begin{align}
\Psi(x,f,g)&:=g(x)\left(a(x)+2b(x)f(x)+\int_{(0,+\infty)}\theta\left(1-e^{-f(x)\theta}\right)\Pi^{L}(x,{\rm d}\theta)\right)\nonumber\\
&\quad -\pi(x,g)\left(
c(x)+\int_{(0,+\infty)}\theta e^{-\theta\pi(x,f)}\Pi^{NL}(x,{\rm d}\theta)\right).\nonumber
\end{align}
\end{proposition}

\subsection{Kuznetsov measures}\label{sec3}

Let $\{Q_{t}(\mu,\cdot):=\p_{\mu}\left(X_{t}\in\cdot\right):t\ge 0,\ \mu\in\me\}$ be the transition
kernel of the $(\mathfrak{S}_{t},\phi^{L},\phi^{NL})$-superprocess $X$. Then by \eqref{1.4}, we have
\begin{equation}
\int_{\me}e^{-\langle f,\nu\rangle}Q_{t}(\mu,{\rm d}\nu)=\exp
\left(
-\langle V_{t}f,\mu\rangle
\right)
\quad\mbox{ for }\mu\in\me\mbox{ and }t\ge 0.\nonumber
\end{equation}
It implies that $Q_{t}(\mu_{1}+\mu_{2},\cdot)=Q_{t}(\mu_{1},\cdot)*Q_{t}(\mu_{2},\cdot)$ for all $\mu_{1},\mu_{2}\in\me$, and hence $Q_{t}(\mu,\cdot)$ is an infinitely divisible probability measure on $\me$. By the semigroup property of $Q_{t}$, $V_{t}$ satisfies that
\begin{equation}
V_{s}V_{t}=V_{t+s}\quad\mbox{ for all }s,t\ge 0.\nonumber
\end{equation}
Moreover, by the infinite divisibility of $Q_{t}$, each operator $V_{t}$ has the representation
\begin{equation}
V_{t}f(x)=\lambda_{t}(x,f)+\int_{\me^{0}} \left(1-e^{-\langle f,\nu\rangle}\right)L_{t}(x,{\rm d}\nu), \quad t>0, f\in\mathcal{B}^{+}_{b}(E),
\label{3.1}
\end{equation}
where $\lambda_{t}(x,{\rm d}y)$ is a bounded kernel on $E$ and $(1\wedge \nu(1))L_{t}(x,{\rm d}\nu)$ is a bounded kernel from $E$ to $\me^{0}$.
Let $Q^{0}_{t}$ be the restriction of $Q_{t}$ to $\me^{0}$. Let
\begin{equation*}
E_{0}:=\{x\in E:\lambda_{t}(x,E)=0\mbox{ for all }t >0\}.
\end{equation*}
If $x\in E_{0}$, then we get from \eqref{3.1} that
\begin{equation}\nonumber
V_{t}f(x)=\int_{\me^{0}} \left(1-e^{-\langle f,\nu\rangle}\right)L_{t}(x,{\rm d}\nu)\quad\mbox{ for } t>0,\ f\in\mathcal{B}^{+}_{b}(E).\label{4}
\end{equation}
It then follows from \cite[Proposition 2.8 and Theorem A.40]{Li} that for every $x\in E_{0}$, the family of measures $\{L_{t}(x,\cdot):t>0\}$ on $\me^{0}$ constitutes an entrance law for the restricted semigroup $\{Q^{0}_{t}:t\ge 0\}$, and hence there corresponds a unique $\sigma$-finite measure $\mathbb{N}_{x}$ on
$(\mathcal{W}^{+}_{0},\mathcal{F}_{\infty})$
such that $\mathbb{N}_{x}(\{0\})=0$,
and that for any $0<t_{1}<t_{2}<\cdots <t_{n}<+\infty$,
  \begin{align}
 \mathbb{N}_{x}\left(X_{t_{1}}\in {\rm d}\nu_{1},\ X_{t_{2}}\in {\rm d}\nu_{2},\cdots, X_{t_{n}}\in {\rm d}\nu_{n}\right)
 =L_{t_{1}}(x,{\rm d}\nu_{1})Q^{0}_{t_{2}-t_{1}}(\nu_{1},{\rm d}\nu_{2})\cdots Q^{0}_{t_{n}-t_{n-1}}(\nu_{n-1},{\rm d}\nu_{n}).\nonumber
  \end{align}
It immediately follows that for all $t>0$ and $f\in\mathcal{B}^{+}_{b}(E)$,
\begin{equation}
\mathbb{N}_{x}\left(1-e^{-\langle f,X_{t}\rangle}\right)=\int_{\me^{0}} \left(1-e^{-\langle f,\nu\rangle}\right)L_{t}(x,{\rm d}\nu)=V_{t}f(x).\label{N}
\end{equation}
This measure $\mathbb{N}_{x}$ is called the \textit{Kuznetsov measure} corresponding to the entrance law $\{L_{t}(x,\cdot):t>0\}$ or the \textit{excursion law} for the $(\mathfrak{S}_{t},\phi^{L},\phi^{NL})$-superprocess.
When $\mu\in \me$ is supported by $E_0$ and $N({\rm d}w)$ is a Poisson random measure on
$\mathcal{W}^{+}_{0}$ with intensity measure
$$
\int_{E_0}\mu({\rm d}x)\mathbb{N}_{x}(\cdot),
$$
the process defined by
$$
\widetilde{X}_0=\mu; \quad \widetilde{X}_t:=\int_{\mathcal{W}^{+}_{0}}w_tN({\rm d}w),
\quad t>0,
$$
is a realization of the superprocess $(X, \p_\mu)$.
We refer to \cite[section 8.4]{Li} for more details on the Kuznetsov measures.
In the sequel, we assume that

\medskip

\noindent\textbf{Assumption 0.} $E_{+}:=\{x\in E:\ b(x)>0\}\subset E_{0}$.

\medskip

Under this assumption, the Kuznetsov measure $\mathbb{N}_{x}$ exists for every $x\in E_{+}$ when $E_{+}$ is nonempty. It is established in \cite{CRY} that Assumption 0 is automatically true for
superdiffusions with a (purely) local branching mechanism.
In the general case, \cite[Theorem 8.6]{Li} gives the following sufficient condition for Assumption 0:
If there is a spatially independent local branching mechanism $\phi(\lambda)$ taking the form
\begin{equation}
\phi(\lambda)=\alpha\lambda+\beta\lambda^{2}+\int_{(0,+\infty)}\left(e^{-\lambda \theta}-1+\lambda \theta\right)n({\rm d}\theta)\quad\mbox{ for }\lambda\ge 0,\nonumber
\end{equation}
where $\alpha\in\mathbb{R}$, $\beta\in\mathbb{R}^{+}$ and $(\theta\wedge \theta^{2})n({\rm d}\theta)$ is a bounded kernel on $(0,+\infty)$,
such that $
\phi'(\lambda)
\to+\infty$ as $\lambda\to+\infty$, and that the branching mechanism $\psi$ of $X$ is bounded below by $\phi$ in the sense that
\begin{equation}
\psi(x,f)\ge \phi(f(x))\quad\mbox{ for all }x\in E\mbox{ and } f\in\mathcal{B}^{+}_{b}(E),\nonumber
\end{equation}
then $E_{0}=E$.

\section{Fundamental martingale and weak local extinction}\label{Sec3}

In this section we will establish a fundamental martingale
of the form $e^{\lambda_1 t}\langle h, X_t\rangle$ for the superprocess $X$
 in terms of the principal eigenvalue $\lambda_{1}$ and
the corresponding positive eigenfunction $h$ of the
Schr\"{o}dinger operator associated with the mean semigroup.
For a MCSBP (resp. a multitype superdiffusion), if one considers the
$E$-valued spatial motion
on an
enriched state space $E\times I$,
where $I$ is the finite or countable set of types, then the mutation in types is the jumps in the $I$-coordinate, and the associated mean semigroup is generated by a matrix (resp. a coupled elliptic system). So, the spectral theory of matrices (resp. the potential theory for elliptic systems) can be applied. See, for example, \cite[Examples 3.7 and 3.8]{PY}.
For a general non-local branching superprocess considered in
Subsection  \ref{sec1},
the associated Schr\"{o}dinger operator takes the form
$\mathcal{L}-a+\gamma$,
where $\mathcal{L}$ is the generator of underlying spatial motion,
and $\gamma$ is an integral operator given by
$\gamma (f)(x)=\gamma(x,f)$.
Since the integral operator $\gamma$ can be quite general,
the method mentioned above is not applicable.
Instead we characterize the Schr\"{o}dinger operator in terms of the associated bilinear form, and impose some technical assumptions to ensure the existence of $\lambda_{1}$ and $h$.

\begin{defi}
 We call a non-negative measure $\mu$  on $E$ a smooth measure of $\xi$ if there
 is a PCAF $A^{\mu}_{t}$ of $\xi$ such that
 $$\int_{E}f(x)\mu({\rm d}x)=\lim_{t\to 0}\frac{1}{t}\Pi_{m}\left[ \int_{0}^{t}f(
  \xi_{s}
 ){\rm d}A^{\mu}_{s}\right]
 \quad \mbox{ for all }f\in \mathcal{B}^{+}(E).$$
Here $\Pi_{m}(\cdot):=\int_{E}\Pi_{x}(\cdot)m({\rm d}x)$.
This measure $\mu$ is called the Revuz measure of $A^{\mu}_{t}$. Moreover,
we say that a smooth measure $\mu$
belongs to the Kato class $\mathbf{K}(\xi)$,
if \begin{equation*}
\lim_{t\downarrow 0}\sup_{x\in
E}\int_{0}^{t}\int_{E}p(s,x,y)
\mu({\rm d}y){\rm d}s
=0.
\end{equation*}
A function $q$ is said to be in the class $\mathbf{K}(\xi)$ if $q(x)m({\rm d}x)$ is in $\mathbf{K}(\xi)$.
\end{defi}
Clearly all bounded measurable functions are included in $\mathbf{K}(\xi)$.
It is known
(see, e.g., \cite[Proposition 2.1.(i)]{AM} and \cite[Theorem 3.1]{SV})
that if $\nu\in\mathbf{K}(\xi)$, then for every $\eps>0$ there is some constant $A_{\eps}>0$ such that
\begin{equation}
\int_{E}u(x)^{2}\nu({\rm d}x)\le \eps\mathcal{E}(u,u)+A_{\eps}\int_{E}u(x)^{2}m({\rm d}x)\quad\forall u\in\mathcal{F}.\label{i1}
\end{equation}

\medskip

\noindent\textbf{Assumption 1.}
$\int_{E}\gamma(x,\cdot)m({\rm d}x)\in \mathbf{K}(\xi)$,
where $\gamma(x,\cdot)$ is the kernel defined in \eqref{def:gamma-rs}.

\medskip

Under Assumption 1, it follows from \eqref{i1}, the boundedness of $\gamma(x)$ and the inequality
$$|u(x)u(y)|\le \frac{1}{2}(u(x)^{2}+u(y)^{2})$$
that, for every $\eps>0$, there is a constant $K_{\eps}>0$ such that
\begin{equation}\nonumber
\int_{E}\int_{E}u(x)u(y)\gamma(x,{\rm d}y)m({\rm d}x)\le \eps \mathcal{E}(u,u)+K_{\eps}\int_{E}u(x)^{2}m({\rm d}x)\quad\forall u\in\mathcal{F}.
\end{equation}
It follows that the bilinear form $(\mathcal{Q},\mathcal{F})$ defined by
\begin{equation}\nonumber
\mathcal{Q}(u,v):=\mathcal{E}(u,v)+\int_{E}a(x)u(x)v(x)m({\rm d}x)-\int_{E}\int_{E}u(y)v(x)\gamma(x,{\rm d}y)m({\rm d}x), \quad u,v\in\mathcal{F},
\end{equation}
 is closed and that there are positive constants $K$ and $\beta_{0}$ such that $\mathcal{Q}_{\beta_{0}}(u,u):=\mathcal{Q}(u,u)+\beta_{0}(u,u)\ge 0$ for all $u\in\mathcal{F}$, and
$$|\mathcal{Q}(u,v)|\le K \mathcal{Q}_{\beta_{0}}(u,u)^{1/2}\mathcal{Q}_{\beta_{0}}(v,v)^{1/2}\quad\forall u,v\in\mathcal{F}.$$
It then follows from \cite{Kunita} that for such a closed form $(\mathcal{Q},\mathcal{F})$ on $L^{2}(E,m)$, there are unique strongly continuous semigroups $\{T_{t}:t\ge 0\}$ and
$\{\widehat{T}_{t}:t\ge 0\}$
on $L^{2}(E,m)$  such that $\|T_{t}\|_{L^{2}(E,m)}\le e^{\beta_{0}t}$, $\|\widehat{T}_{t}\|_{L^{2}(E,m)}\le e^{\beta_{0}t}$, and
\begin{equation}
(T_{t}f,g)=(f,\widehat{T}_{t}g)\quad\forall f,g\in L^{2}(E,m).\label{duality}
\end{equation}
Let $\{U_{\alpha}\}_{\alpha>\beta_{0}}$ and $\{\widehat{U}_{\alpha}\}_{\alpha>\beta_{0}}$ be given by $U_{\alpha}f:=\int_{0}^{+\infty}e^{-\alpha t}T_{t}f {\rm d}t$ and $\widehat{U}_{\alpha}f:=\int_{0}^{+\infty}e^{-\alpha t}\widehat{T}_{t}f {\rm d}t$ respectively. Then $\{U_{\alpha}\}_{\alpha>\beta_{0}}$ and $\{\widehat{U}_{\alpha}\}_{\alpha>\beta_{0}}$ are strongly continuous pseudo-resolvents in the sense that they satisfy the resolvent equations
$$U_{\alpha}-U_{\beta}+(\alpha-\beta)U_{\alpha}U_{\beta}=0,\quad \widehat{U}_{\alpha}-\widehat{U}_{\beta}+(\alpha-\beta)\widehat{U}_{\alpha}\widehat{U}_{\beta}=0$$
for all $\alpha,\beta>\beta_{0}$, and
\begin{equation}\label{QalphaU}
\mathcal{Q}_{\alpha}(U_{\alpha}f,g)=\mathcal{Q}_{\alpha}(g,\widehat{U}_{\alpha}f)=(f,g)\quad\forall f\in L^{2}(E,m),\ g\in\mathcal{F}.
\end{equation}
Recall from Proposition \ref{prop0} that $\mathfrak{P}_{t}$ is the mean
semigroup of the $(\mathfrak{S}_{t},\Phi^{L},\Phi^{NL})$-superprocess, which satisfies the equation \eqref{li2}.
Since
$\gamma(x,{\rm d}y)$ is a bounded kernel on $E$, by \eqref{li2}, we have for every $f\in\mathcal{B}_{b}(E)$,
$$
\|\mathfrak{P}_{t}f\|_{\infty}\le \|f\|_{\infty}+(\|a\|_{\infty}+\|\gamma(\cdot,1)\|_{\infty})\int_{0}^{t}\|\mathfrak{P}_{t-s}f\|_{\infty}{\rm d}s.
$$
By Gronwall's lemma, $\|\mathfrak{P}_{t}f\|_{\infty}\le e^{c_{1}t}\|f\|_{\infty}$
for some constant $c_{1}>0$.
For $f\in\mathcal{B}_{b}(E)$ and $\alpha>c_{1}$, define
$R_{\alpha}f(x):=\int_{0}^{+\infty}e^{-\alpha t}\mathfrak{P}_{t}f(x){\rm d}t$.
By taking Laplace transform on both sides of \eqref{li2}, we get
\begin{equation}
R_{\alpha}f(x)=G_{\alpha}f(x)-G_{\alpha}\left(\alpha R_{\alpha}f\right)(x)+G_{\alpha}\left(\gamma(\cdot,R_{\alpha}f)\right)(x),\label{r1}
\end{equation}
where $G_{\alpha}$ is the $\alpha$-resolvent of
$(\mathfrak{S}_{t})_{t\ge 0}$.
A particular case is when
$a(x),\gamma(x)\in L^{2}(E,m)$. In this case, for all $f\in \mathcal{B}_{b}(E)\cap L^{2}(E,m)$
and $\alpha$ sufficiently large, both $a(x) R_{\alpha}f(x)$ and $\gamma(x,R_{\alpha}f)$ are in $\mathcal{B}_{b}(E)\cap L^{2}(E,m)$. Then it follows from \eqref{d1} that $G_{\alpha}f$, $G_{\alpha}\left(\alpha R_{\alpha}f\right)$, $G_{\alpha}\left(\gamma(\cdot,R_{\alpha}f)\right)\in\mathcal{F}$, and then
by \eqref{r1}, \eqref{d1} and \eqref{QalphaU},
$$
\mathcal{Q}_{\alpha}(R_{\alpha}f,v)=(f,v)
=\mathcal{Q}_{\alpha}(U_{\alpha}f,v)
\quad\mbox{ for all }v\in\mathcal{F},
$$
which implies that $R_{\alpha}f$ is $m$-equivalent to $U_{\alpha}f$ for
$\alpha$ sufficiently large.
This indicates that there is some strong relation
between $\mathfrak{P}_{t}$ and $T_{t}$.
In fact we will show in Proposition \ref{prop4} below that
$\mathfrak{P}_{t}f=T_{t}f\ [m]$ for every $t>0$ and
every $f\in \mathcal{B}_{b}(E)\cap L^{2}(E,m)$.
This means that $\mathfrak{P}_{t}$ can be regarded as a
bounded linear operator on the space of bounded measurable functions in $L^{2}(E,m)$, which is dense in $L^{2}(E,m)$. Hence $T_{t}$ can be regarded as the unique bounded linear operator on $L^{2}(E,m)$
which is an extension of $\mathfrak{P}_{t}$.

\medskip

\noindent\textbf{Assumption 2.}
There
exist a constant $\lambda_{1}\in (-\infty,+\infty)$ and
positive functions $h,\widehat{h}\in\mathcal{F}$ with $h$ bounded continuous, $\|h\|_{L^{2}(E,m)}=1$ and $(h,\widehat{h})=1$ such that
\begin{equation}
\mathcal{Q}(h,v)=\lambda_{1}(h,v),\quad
\mathcal{Q}(v,\widehat{h})=
\lambda_{1}(v,\widehat{h})\quad\forall v\in\mathcal{F}.\label{e1}
\end{equation}

\medskip

In the Theorem \ref{them1} below,
we will prove that $e^{\lambda_1 t}\langle h, X_t\rangle$ is a non-negative martingale. To prove this we first prove that $h$ is invariant for some semigroup, see \eqref{3.13} below.
Since $h\in \mathcal{F}$ is continuous,
it follows from \cite[Theorem 4.2.6]{CF} that for
every $x\in E$, $\Pi_{x}$-a.s.
$$
h(\xi_{t})-h(\xi_{0})=M^{h}_{t}+N^{h}_{t},\quad    t\ge 0,
$$
where $M^{h}$ is a martingale additive functional of $\xi$ having
finite energy and $N^{h}_{t}$ is a continuous additive functional of
$\xi$ having zero energy.
The formula above is usually called Fukushima's decomposition.
 It follows from \eqref{e1} and \cite[Theorem 5.4.2]{FOT} that $N^{h}_{t}$ is of
bounded variation, and
$$N^{h}_{t}=-\lambda_{1}\int_{0}^{t}h(\xi_{s}){\rm d}s+\int_{0}^{t}a(\xi_{s})h(\xi_{s}){\rm d}s-\int_{0}^{t}\gamma(\xi_{s},h){\rm d}s,\quad\forall t\ge 0.$$
Following the idea of \cite[Section 2]{CFTYZ},  we define
a local martingale on the random time interval
$[0,\zeta_{p})$ by
\begin{equation} \label{(1)}
M_{t}:=\int_{0}^{t}\frac{1}{h(\xi_{s-})}{\rm d}M^{h}_{s},
\quad  t\in [0,\zeta_{p}),
\end{equation}
where $\zeta_{p}$ is the predictable part of the lifetime $\zeta$
of $\xi$, that is,
$$
\zeta_{p}=\begin{cases} \zeta & \mbox{ if } \zeta<\infty \mbox{ and } \xi_{\zeta-}=\xi_{\zeta},\cr
\infty, & \mbox{ otherwise},
\end{cases}
$$
see \cite[Theorem 44.5]{Sharpe}.
 Then the solution $H_{t}$ of the stochastic differential
equation
\begin{equation}
H_{t}=1+\int_{0}^{t}H_{s-}{\rm d}M_{s},  \quad
 t\in [0,\zeta_{p}),\label{(2)}
\end{equation}
is a positive local martingale on $[0,\zeta_{p})$ and hence a
supermartingale. Consequently, the formula
$$
{\rm d}\Pi^{h}_{x}=H_{t}\,{\rm d}\Pi_{x}\quad
\mbox{on }\mathcal{H}_{t}\cap\{t<\zeta\}\quad\mbox{ for }x\in E
$$
uniquely determines a family of subprobability measures
$\{\Pi^{h}_{x}:x\in E\}$ on $(\Omega,\mathcal{H})$. Hence we have
$$
\Pi^{h}_{x}\left[f(\xi_{t})\right]=\Pi_{x}\left[H_{t}f(\xi_{t});t<\zeta\right],
\quad t\ge 0, f\in\mathcal{B}^{+}(E).
$$
Note that by \eqref{(1)},
\eqref{(2)} and Dol\'{e}an-Dade's formula,
\begin{equation}\label{(3)}
H_{t}=\exp\left( M_{t}-\frac{1}{2}
\langle M^{c}\rangle_{t}\right) \prod_{0<s\le
t}\frac{h(\xi_{s})}{h(\xi_{s-})}\exp\left( 1-\frac{h(\xi_{s})}{h(\xi_{s-})}\right)
\quad \forall t\in [0, \zeta_p),
\end{equation}
where $M^{c}$ is the continuous martingale part of $M$. Applying
Ito's formula to $\log h(\xi_{t})$, we obtain that for every $x\in E$,
$\Pi_{x}$-a.s. on  $[0,\zeta)$,
\begin{align}
\log h(\xi_{t})-\log h(\xi_{0})
&=M_{t}-\frac{1}{2}\langle
M^{c}\rangle_{t}+\sum_{s\le t}
\left(\log\frac{h(\xi_{s})}{h(\xi_{s-})}-\frac{h(\xi_{s})-h(\xi_{s-})}{h(\xi_{s-})}\right)\nonumber\\
&\quad-\lambda_{1}t
+\int_{0}^{t}a(\xi_{s}){\rm d}s-\int_{0}^{t}\frac{\gamma(\xi_{s},h)}{h(\xi_{s})}\,{\rm d}s.\label{(4)}
\end{align}
Put
\begin{equation}\label{def:q-rs}
q(x):=\frac{\gamma(x,h)}{h(x)}\quad\mbox{ for }x\in E.
\end{equation}
By \eqref{(3)} and \eqref{(4)}, we get
$$
H_{t}=\exp\left( \lambda_{1}t-\int_{0}^{t}a(\xi_{s}){\rm d}s
+\int_{0}^{t}q(\xi_{s}){\rm d}s\right)
\frac{h(\xi_{t})}{h(\xi_{0})}.
$$
To emphasize, the process $\xi$ under $\{\Pi^{h}_{x},x\in E\}$ will be
denoted as $\xi^{h}$.
For a measurable function $g$, we set
 $$e_{g}(t):=\exp\left(-\int_{0}^{t}g(\xi_{s}){\rm d}s\right)\quad\forall
 t\ge 0,
 $$
 whenever it is well defined.
Then we have for all $f\in\mathcal{B}^{+}(E)$ and $t\ge 0$,
\begin{equation}
\mathfrak{S}^h_tf(x):=\Pi^{h}_{x}\left[f(\xi^{h}_{t})\right]
=\frac{e^{\lambda_{1}t}}{h(x)}\Pi_{x}\left[e_{a-q}(t)h(\xi_{t})f(\xi_{t})\right].\label{pihx}
\end{equation}
It follows from
\cite[Theorem 2.6]{CFTYZ} that
the transformed process $\xi^{h}$
is
a conservative and recurrent (in the sense of \cite{FOT})
$\widetilde m$-symmetric right Markov process on $E$ with $\widetilde{m}({\rm d}y):=h(y)^{2}m({\rm d}y)$.
Thus
\begin{equation}
\mathfrak{S}^{h}_{t}1=1\ [\widetilde{m}]\quad\mbox{ for all }t>0.\label{3.13'}
\end{equation}
Note that for all $t>0$ and $x\in E$, the measure $\mathfrak{S}^{h}_{t}(x,\cdot):=\Pi^{h}_{x}\left(\xi^{h}_{t}\in\cdot\right)$ is absolutely continuous with respect to $\widetilde{m}$, since $\mathfrak{S}^{h}_{t}(x,\cdot)$ is absolutely continuous with respect to
the measure $\mathfrak{S}_{t}(x,\cdot):=\Pi_{x}(\xi_{t}\in\cdot)$ by \eqref{pihx} and the latter is absolutely continuous with respect to $m$. Moreover, by the right continuity of the sample paths of $\xi^{h}$, one can easily verify that both $1$ and $\mathfrak{S}^{h}_{t}1(x)$ are excessive functions for $\{\mathfrak{S}^{h}_{t}:t>0\}$. Thus by \cite[Theorem A.2.17]{CF}, \eqref{3.13'} implies that
\begin{equation}
1=\mathfrak{S}^{h}_{t}1(x)
=\frac{e^{\lambda_{1}t}}{h(x)}\Pi_{x}\left[e_{a-q}(t)h(\xi_{t})\right]\quad\mbox{ for all }x\in E.\label{3.13}
\end{equation}

\begin{theorem}\label{them1}
Suppose Assumptions 1-2 hold.
Then for every $\mu\in\me$,
$W^{h}_{t}(X):=e^{\lambda_{1}t}\langle h,X_{t}\rangle$ is a non-negative $\p_{\mu}$-martingale with respect to the filtration $\{\mathcal{F}_{t}:t\ge 0\}$.
\end{theorem}
\proof
By the Markov property of $X$, it suffices to prove that for all $x\in E$ and $t\ge 0$,
\begin{equation}
\mathfrak{P}_{t}h(x)=\p_{\delta_{x}}\left(\langle h,X_{t}\rangle\right)
=e^{-\lambda_{1}t}h(x).\label{prop1.1}
\end{equation}
Let $A(s,t):=-\int_{s}^{t}a(\xi_{r}){\rm d}r+\int_{s}^{t}q(\xi_{r}){\rm d}r$ and $u(t,x):=\Pi_{x}
\left[
e^{A(0,t)}h(\xi_{t})
\right]
$. Clearly by \eqref{3.13}, $u(t,x)=e^{-\lambda_{1}t}h(x)$. Note that
\begin{equation}
e^{A(0,t)}-1=-(e^{A(t,t)}-e^{A(0,t)})=\int_{0}^{t}\left(-a(\xi_{s})+q(\xi_{s})\right)e^{A(s,t)}{\rm d}s.\label{prop1.2}
\end{equation}
By \eqref{prop1.2}, Fubini's theorem and the Markov property of $\xi$, we have
\begin{align}
&{ }u(t,x)\nonumber\\
&=\mathfrak{S}_{t}h(x)+\Pi_{x}\left[(e^{A(0,t)}-1)h(\xi_{t})\right]\nonumber\\
&=\mathfrak{S}_{t}h(x)-\Pi_{x}\left[\int_{0}^{t}a(\xi_{s})e^{A(s,t)}h(\xi_{t}){\rm d}s\right]+\Pi_{x}\left[\int_{0}^{t}q(\xi_{s})e^{A(s,t)}h(\xi_{t}){\rm d}s\right]\nonumber\\
&=\mathfrak{S}_{t}h(x)-\Pi_{x}\left[\int_{0}^{t}a(\xi_{s})u(t-s,\xi_{s}){\rm d}s\right]+\Pi_{x}\left[\int_{0}^{t}\frac{\gamma(\xi_{s},h)}{h(\xi_{s})}
u(t-s,\xi_{s}){\rm d}s\right]\nonumber\\
&=\mathfrak{S}_{t}h(x)-\Pi_{x}\left[\int_{0}^{t}a(\xi_{s})u(t-s,\xi_{s}){\rm d}s\right]+\Pi_{x}\left[\int_{0}^{t}\gamma(\xi_{s},u^{t-s}){\rm d}s\right].\nonumber
\end{align}
In the last equality above we use the fact that $u(t-s,x)=e^{-\lambda_{1}(t-s)}h(x)$.
Thus $u(t,x)$ is a locally bounded solution to \eqref{li2} with initial value $h$. By the uniqueness of the solution,
we get $u(t,x)=
\mathfrak{P}_{t}h(x)
=\p_{\delta_{x}}\left(\langle h,X_{t}\rangle\right)$.
\qed

\medskip

For $\mu\in\me$, we say that  the process $X$ exhibits \textit{weak local extinction} (resp. \textit{local extinction}) under $\p_{\mu}$ if for every nonempty relatively compact open subset $B$ of $E$, $\p_{\mu}\left(\lim_{t\to+\infty}X_{t}(B)=0\right)=1$ (resp. $\p_{\mu}\left(X_{t}(B)=0\mbox{ for sufficitently large $t$}\right)=1$).
It is proved in \cite{EP99} (see also \cite{EK}) that local extinction and weak local extinction coincide for superdiffusions in a domain $D\subset\R^d$ with local branching mechanism $\psi(x,\lambda)=-\beta(x)\lambda+\alpha(x)\lambda^{2}$, where  $\alpha$ and $\beta$ are H\"{o}lder continuous functions on $D$ with order $\eta\in(0,1]$, $\alpha>0$ in $D$ and $\beta$ is bounded from above.
However the two notions are different in general. In this paper we are only concerned with weak local extinction.

\begin{corollary}\label{cor1}
Suppose Assumptions 1-2 hold.
For all $\mu\in\me$ and nonempty
relatively compact open subset $B$ of $E$,
$$\p_{\mu}\left(\limsup_{t\to+\infty}e^{\lambda_{1}t}X_{t}(B)<+\infty\right)=1.$$
In particular, if $\lambda_{1}>0$, then $X$  exhibits weak local extinction under $\p_{\mu}$.
\end{corollary}
\proof This corollary follows immediately from Theorem \ref{them1} and the fact that
$$e^{\lambda_{1}t}X_{t}(B)\le e^{\lambda_{1}t}\langle \frac{h}{\inf_{x\in B}h(x)}1_{B},X_{t}\rangle \le \frac{1}{\inf_{x\in B}h(x)}W^{h}_{t}(X).$$\qed

\begin{remark}\rm
Corollary \ref{cor1} implies that the local mass of $X_{t}$ grows subexponentially and the growth rate can not exceed $-\lambda_{1}$.
However, when one considers the total mass process
$\langle 1,X_{t}\rangle$, the growth rate may actually exceed $-\lambda_{1}$. We refer to \cite{EK} and \cite{ERS} for more concrete examples.
\end{remark}

\section{Spine decomposition}\label{Sec4}

\subsection{Concatenation process}

We assume Assumptions 1-2 hold.
It is well-known (see, e.g., \cite[p. 286]{Sharpe}) that for every $x\in E$,
there is a unique (up to equivalence in law)
right process
$((\widehat{\xi}_{t})_{t\ge 0};\widehat{\Pi}^{h}_{x})$
on $E$ with lifetime $\widehat{\zeta}$ and cemetery point $\partial$,
such that
$$
\widehat{\Pi}^{h}_{x}\left(\widehat{\xi}_{t}\in B\right)=\Pi^{h}_{x}\left[e_{q}(t);\xi^{h}_{t}\in B\right]\quad\forall B\in\mathcal{B}(E),
$$
where $q$ is the nonnegative function defined in \eqref{def:q-rs}.
$\widehat{\xi}$ is called the $e_{q}(t)$-subprocess of $\xi^{h}$, which can be obtained by killing $\xi^h$ with rate $q$.
In fact, a version of the $e_{q}(t)$-subprocess can be obtained by the following method of curtailment of the lifetime.
Let $Z$ be an exponential random variable, of parameter 1, independent of
$\xi^{h}$.
Put
$$\widehat{\zeta}(\omega):=\inf\{t\ge 0:\ \int_{0}^{t}q\left(\xi^{h}_{s}(\omega)\right){\rm d}s\ge Z(\omega)\}
\ (=+\infty, \mbox{ if such $t$ does not exist})
,$$
and
\begin{equation} \nonumber
\widehat{\xi}_{t}(\omega):=\begin{cases}
         \xi^{h}_{t}(\omega)
         \quad &\hbox{if } t<\widehat{\zeta}(\omega) , \smallskip \\
         \partial
         \quad &\hbox{if } t\ge \widehat{\zeta}(\omega).
        \end{cases}
\end{equation}
 Then the process $((\widehat{\xi}_{t})_{t\ge 0},\Pi^{h}_{x})$ is equal in law to the $e_{q}(t)$-subprocess of $\xi^{h}$.
Now we define a probability on $E$ by
\begin{equation}\label{def:pihrs}
\pi^{h}(x,{\rm d}y):=\frac{h(y)\gamma(x,{\rm d}y)}{\gamma(x,h)}=\frac{h(y)\pi(x,{\rm d}y)}{\pi(x,h)}\quad\mbox{ for }x\in E.
\end{equation}
Let
$\widetilde{\xi}:=(\widetilde{\Omega},\widetilde{\mathcal{G}},\widetilde{\mathcal{G}}_{t},\widetilde{\theta}_{t},\widetilde{\xi_{t}},\widetilde{\Pi}_{x},\widetilde{\zeta})$
be the right process constructed from $\widehat{\xi}$ and
the instantaneous distribution
$\kappa(\omega,{\rm d}y):=\pi^{h}(\widehat{\xi}_{\widehat{\zeta}(\omega)-}(\omega),{\rm d}y)$
by using the so-called ``piecing out" procedure (cf. Ikeda et al.
\cite{INW}).
We will follow the terminology of \cite[Section II.14]{Sharpe} and call
$\widetilde{\xi}$ a \textit{concatenation} process defined from an infinite sequence of copies of $\widehat{\xi}$ and the transfer kernel $\kappa(\omega,{\rm d}y)$.
One can also refer to \cite[Section A.6]{Li} for a summary of concatenation processes. The intuitive idea of this concatenation is described as follows. The process $\widetilde{\xi}$ evolves as
the process $\xi^{h}$
until time $\widehat{\zeta}$, it is then revived by means of the kernel $\kappa(\omega,{\rm d}y)$ and evolves again as
$\xi^{h}$
and so on, until a countably infinite number of revivals have occurred.
Clearly in the case of purely local branching mechanism (i.e. $\gamma(x)\equiv 0$ on $E$), we have $\widehat{\zeta}=+\infty$ almost surely and hence $\widetilde{\xi}$ runs as a copy of $\xi^{h}$.

Let $\widetilde{\mathfrak{S}}_{t}$ be the transition semigroup of $\widetilde{\xi}$.
It satisfies the following renewal equation.
\begin{equation}
\widetilde{\mathfrak{S}}_{t}f(x)=\Pi_{x}^{h}\left[e_{q}(t)f(\xi^{h}_{t})\right]+\Pi_{x}^{h}\left[\int_{0}^{t}q(\xi^{h}_{s})e_{q}(s)
\pi^{h}(\xi^{h}_{s},\widetilde{\mathfrak{S}}_{t-s}f){\rm d}s\right], \quad f\in\mathcal{B}^{+}_{b}(E).
\label{s1}
\end{equation}
By \cite[Proposition 2.9]{Li}, the above equation can be rewritten as
\begin{align*}
\widetilde{\mathfrak{S}}_{t}f(x)&=\Pi_{x}^{h}\left[f(\xi^{h}_{t})\right]-\Pi_{x}^{h}\left[\int_{0}^{t}q(\xi^{h}_{s})\widetilde{\mathfrak{S}}_{t-s}f(\xi^{h}_{s}){\rm d}s\right]
+\Pi_{x}^{h}\left[\int_{0}^{t}q(\xi^{h}_{s})\pi^{h}(\xi^{h}_{s},\widetilde{\mathfrak{S}}_{t-s}f){\rm d}s\right].
\end{align*}

\begin{proposition}\label{prop2}
For all $f\in\mathcal{B}^{+}_{b}(E)$, $t\ge 0$ and $x\in E$,
\begin{equation}
\widetilde{\mathfrak{S}}_{t}f(x)=\frac{e^{\lambda_{1}t}}{h(x)}\mathfrak{P}_{t}(fh)(x).
\label{prop2.0}
\end{equation}
In particular, $\widetilde{\mathfrak{S}}_{t}1(x)\equiv 1$,
and hence $\widetilde{\xi}$ has
infinite lifetime.
Moreover, for each $t>0$ and $x\in E$, $\widetilde{\xi}$ has a transition density $\widetilde{p}(t,x,y)$ with respect to the probability measure $\rho({\rm d}y):=h(y)\widehat{h}(y)m({\rm d}y)$.
\end{proposition}
\proof
By \eqref{s1},  \eqref{pihx}, \eqref{def:gamma-rs}, \eqref{def:q-rs} and \eqref{def:pihrs}, we have
\begin{align}
\widetilde{\mathfrak{S}}_{t}f(x)
&=\frac{e^{\lambda_{1}t}}{h(x)}\Pi_{x}\left[e_{a}(t)h(\xi_{t})f(\xi_{t})\right]\nonumber\\
&\quad +\frac{e^{\lambda_{1}t}}{h(x)}\Pi_{x}\left[\int_{0}^{t}e_{a}(s)q(\xi_{s})h(\xi_{s})e^{-\lambda_{1}(t-s)}\pi^{h}(\xi_{s},
\widetilde{\mathfrak{S}}_{t-s}f){\rm d}s\right]\nonumber\\
&=\frac{e^{\lambda_{1}t}}{h(x)}\Pi_{x}\left[e_{a}(t)h(\xi_{t})f(\xi_{t})\right]\nonumber\\
&\quad +\frac{e^{\lambda_{1}t}}{h(x)}\Pi_{x}\left[\int_{0}^{t}e_{a}(s)\gamma(\xi_{s},e^{-\lambda_{1}(t-s)}h\widetilde{\mathfrak{S}}_{t-s}f){\rm d}s\right].\label{prop2.1}
\end{align}
Let $u(t,x):=e^{-\lambda_{1}t}h(x)\widetilde{\mathfrak{S}}_{t}f(x)$. Clearly $u(t,x)$ is a locally bounded function on $[0,+\infty)\times E$. Moreover, it follows from \eqref{prop2.1} and \cite[Proposition 2.9]{Li} that
\begin{align}
u(t,x)&=\Pi_{x}\left[e_{a}(t)h(\xi_{t})f(\xi_{t})\right]
+\Pi_{x}\left[\int_{0}^{t}e_{a}(s)\gamma(\xi_{s},u^{t-s}){\rm d}s\right]\nonumber\\
&=\Pi_{x}\left[h(\xi_{t})f(\xi_{t})\right]-\Pi_{x}\left[\int_{0}^{t}a(\xi_{s})u(t-s,\xi_{s}){\rm d}s\right]
+\Pi_{x}\left[\int_{0}^{t}\gamma(\xi_{s},u^{t-s}){\rm d}s\right].
\nonumber
\end{align}
This implies that $u(t,x)$ is a locally bounded solution to \eqref{li2} with initial value $fh$. Hence we get $e^{-\lambda_{1}t}h(x)\widetilde{\mathfrak{S}}_{t}f(x)
=u(t,x)=\mathfrak{P}_{t}(fh)(x)$ by the uniqueness of the solution. It then follows from
\eqref{prop1.1}
that $\widetilde{\mathfrak{S}}_{t}1(x)\equiv 1$ on $E$.

To prove the second part of this proposition, it suffices to prove that for all $t>0$ and $x\in E$, $\widetilde{\mathfrak{S}}_{t}1_{B}(x)=0$ for all $B\in \mathcal{B}(E)$ with $\rho(B)=0$ (or equivalently, $m(B)=0$).
Note that $\Pi_{x}\left[h1_{B}(\xi_{t})\right]=\int_{B}h(y)p(t,x,y)m({\rm d}y)=0$. It follows from the above argument that $e^{-\lambda_{1}t}h(x)\widetilde{\mathfrak{S}}_{t}1_{B}(x)=
\mathfrak{P}_{t}(h1_{B})(x)$ is the unique locally bounded solution to
\eqref{li2} with initial value $0$. Thus $\widetilde{\mathfrak{S}}_{t}1_{B}(x)\equiv 0$.
\qed

\medskip

\begin{remark}\rm
The formula \eqref{prop2.0} can be written as
\begin{equation}\label{many-to-one}
\frac{\p_{\delta_{x}}\left[\langle f h,X_{t}\rangle\right]}{\p_{\delta_{x}}\left[\langle h,X_{t}\rangle\right]}=\widetilde{\Pi}_{x}\left[f(\widetilde{\xi}_{t})\right]
\quad\mbox{ for }f\in\mathcal{B}^{+}_{b}(E)\mbox{ and } t\ge 0,
\end{equation}
which enables us to calculate the first moment of the superprocess in
terms of an auxiliary process.
An analogous formula for a special class of non-local branching Markov processes, which is called a ``many-to-one" formula, is established in \cite{BDMT}, but with a totally different method.
By \eqref{prop1.1}, we may rewrite \eqref{many-to-one} as
\begin{equation*}
\p_{\delta_{x}}\left[\langle f h,X_{t}\rangle\right]=e^{\lambda_{1}t}h(x)
\widetilde{\Pi}_{x}\left[f(\widetilde{\xi}_{t})\right]
\quad\mbox{ for }f\in\mathcal{B}^{+}_{b}(E)\mbox{ and } t\ge 0.
\end{equation*}
\end{remark}

\medskip

Let
$\tau_{1}$ be the first revival time of $\widetilde{\xi}$. For $n\ge 2$, define $\tau_{n}$ recursively by $\tau_{n}:=
\tau_{n-1}
+\tau_{1}\circ \widetilde{\theta}_{\tau_{n-1}}$.
Since $\widetilde{\xi}$ has infinite lifetime,
$\widetilde{\Pi}_{x}\left(\lim_{n\to+\infty}\tau_{n}=+\infty\right)=1$ for all $x\in E$.

\begin{proposition}\label{prop3}
For all $f(s,x,y),g(s,x,y)\in\mathcal{B}^{+}([0,+\infty)\times E\times E)$, $t>0$ and $x\in E$, we have
\begin{equation}
\widetilde{\Pi}_{x}\left[\sum_{\tau_{i}\le t}f(\tau_{i},\widetilde{\xi}_{\tau_{i}-},\widetilde{\xi}_{\tau_{i}})\right]
=\widetilde{\Pi}_{x}\left[\int_{0}^{t}{\rm d}s\int_{E}\pi^{h}(\widetilde{\xi}_{s},{\rm d}y)q(\widetilde{\xi}_{s})f(s,\widetilde{\xi}_{s},y)\right]\label{l1}
\end{equation}
and
\begin{multline}
\widetilde{\Pi}_{x}\left[\left(\sum_{\tau_{i}\le t}f(\tau_{i},\widetilde{\xi}_{\tau_{i}-},\widetilde{\xi}_{\tau_{i}})\right)\left(\sum_{\tau_{j}\le t}g(\tau_{j},\widetilde{\xi}_{\tau_{j}-},\widetilde{\xi}_{\tau_{j}})\right)\right]\\
=\widetilde{\Pi}_{x}\left[\sum_{\tau_{i}\le t}f g(\tau_{i},\widetilde{\xi}_{\tau_{i}-},\widetilde{\xi}_{\tau_{i}})\right]
+\widetilde{\Pi}_{x}\big[\int_{0}^{t}{\rm d}s\int_{E}\pi^{h}(\widetilde{\xi_{s}},{\rm d}y)q(\widetilde{\xi}_{s})f(s,\widetilde{\xi_{s}},y)\\
\shoveright{\cdot\widetilde{\Pi}_{y}\left(\int_{0}^{t-s}{\rm d}r\int_{E}\pi^{h}(\widetilde{\xi}_{r},{\rm d}z)q(\widetilde{\xi}_{r})g(s+r,\widetilde{\xi}_{r},z)\right)
\big]}\\
\shoveleft{\quad\quad\quad\quad\quad\quad+\widetilde{\Pi}_{x}\big[\int_{0}^{t}{\rm d}s\int_{E}\pi^{h}(\widetilde{\xi_{s}},{\rm d}y)q(\widetilde{\xi}_{s})g(s,\widetilde{\xi_{s}},y)}\\
{\shoveright{\cdot\widetilde{\Pi}_{y}\left(\int_{0}^{t-s}{\rm d}r\int_{E}\pi^{h}(\widetilde{\xi}_{r},{\rm d}z)q(\widetilde{\xi}_{r})f(s+r,\widetilde{\xi}_{r},z)\right)
\big].}}\label{l2}
\end{multline}
\end{proposition}

The proof of this proposition will be given in the Appendix below.

\subsection{Spine decomposition}

In this section we work under Assumptions 0-2.
Recall from Theorem \ref{them1} that the process $W^{h}_{t}(X)$
is a non-negative $\p_{\mu}$-martingale for every $\mu\in\me$.
We can define a new probability measure $\q_{\mu}$ for every $\mu\in\me^0$ by the following formula:
\begin{equation*}
\left.{\rm d}\q_{\mu}\right|_{\mathcal{F}_{t}}:=\frac{1}{\langle h,\mu\rangle }\left.W^{h}_{t}(X){\rm d}\p_{\mu}\right|_{\mathcal{F}_{t}}\quad\mbox{ for all }t\ge 0.
\end{equation*}
It then follows from Proposition \ref{prop0} that for any $f\in\mathcal{B}^{+}_{b}(E)$ and $t\ge 0$,
\begin{equation*}
\q_{\mu}\left(e^{-\langle f,X_{t}\rangle}\right)=\frac{e^{\lambda_{1}t}}{\langle h,\mu\rangle}
\p_{\mu}\left(\langle h,X_{t}\rangle e^{-\langle f,X_{t}\rangle}\right)=\frac{e^{\lambda_{1}t}}{\langle h,\mu\rangle}
e^{-\langle V_{t}f,\mu\rangle}\langle V^{h}_{t}f,\mu\rangle,
\end{equation*}
where $V^{h}_{t}f(x)$ is the unique locally bounded solution to \eqref{li1} with initial value $h$.
In this subsection we will establish the spine decomposition of $X$ under $\q_{\mu}$.

\begin{defi}\label{sd1}
For all $\mu\in\me$ and $x\in E$, there is a probability space with probability measure $\pp_{\mu,x}$ that carries the following processes.
\begin{description}
\item{\rm(i)}
$((\widetilde{\xi}_{t})_{t\ge 0};\pp_{\mu,x})$ is equal
in law to $\widetilde{\xi}$,
a copy of the concatenation process starting from $x$;

\item{\rm(ii)}
$(n;\pp_{\mu,x})$ is a random measure such that,
given $\widetilde{\xi}$ starting from $x$,
$n$ is a Poisson random measure which issues $\me$-valued processes $X^{n,t}:=(X^{n,t}_{s})_{s\ge 0}$ at space-time points $(\widetilde{\xi}_{t},t)$ with rate
$$
{\rm d}\mathbb{N}_{\widetilde{\xi}_{t}}\times 2b(\widetilde{\xi}_{t}){\rm d}t.
$$
Here for every $y\in E_{+}=\{z\in E:\ b(z)>0\}$, $\mathbb{N}_{y}$
denotes the Kuznetsov measure on $\mathcal{W}^{+}_{0}$ corresponding to the
$(\mathfrak{S}_{t},\phi^{L},\phi^{NL})$-superprocess, while for
$y\in E\setminus E_{+}$,
$\mathbb{N}_{y}$ denotes the null measure
on $\mathcal{W}^{+}_{0}$.
Note that, given $\widetilde{\xi}$,
immigration happens only at space-time points
$(\widetilde{\xi}_{t},t)$ with $b(\widetilde{\xi}_{t})>0$. Let $D^{n}$ denote the almost surely countable set of immigration times, and $D^{n}_{t}:=D^{n}\cap [0,t]$. Given $\widetilde{\xi}$, the processes $\{X^{n,t}:t\in D^{n}\}$ are mutually independent.

\item{\rm(iii)}
$(m;\pp_{\mu,x})$ is a random measure such that,
given $\widetilde{\xi}$ starting from $x$,
$m$ is a Poisson random measure which issues $\me$-valued processes $X^{m,t}:=(X^{m,t}_{s})_{s\ge 0}$ at space-time points $(\widetilde{\xi}_{t},t)$
with initial mass $\theta$ at rate
$$
\theta\Pi^{L}(\widetilde{\xi}_{t},{\rm d}\theta)\times {\rm d}\p_{\theta\delta_{\widetilde{\xi}_{t}}}\times {\rm d}t.
$$
Here $\p_{\theta\delta_{x}}$ denotes the law of the
$(\mathfrak{S}_{t},\phi^{L},\phi^{NL})$-superprocess starting from $\theta\delta_{x}$.
 Let $D^{m}$ denote the almost surely countable set of immigration times, and $D^{m}_{t}:=D^{m}\cap [0,t]$.
 Given $\widetilde{\xi}$, the processes $\{X^{m,t}:t\in D^{m}\}$ are mutually independent, also independent of $n$ and $\{X^{n,t}:t\in D^{n}\}$.

 \item{\rm (iv)}
$\{((X^{r,i}_{s})_{s\ge 0}; \pp_{\mu,x}),\ i\ge 1\}$
is a family of $\me$-valued processes such that,
given $\widetilde{\xi}$ starting from $x$ (including its revival times
$\{\tau_{i}:i\ge 1\}$),
$X^{r,i}:=(X^{r,i}_{s})_{s \ge 0}$ is equal
in law  to $((X_{s})_{s\ge 0},\p_{\pi_{i}})$
where $\p_{\pi_{i}}$ denotes the law of the $(\mathfrak{S}_{t},\phi^{L},\phi^{NL})$-superprocess starting from $\pi_{i}(\cdot):=\Theta_{i}\pi(\widetilde{\xi}_{\tau_{i}-},\cdot)$ and $\Theta_{i}$ is a
$[0,+\infty)$-valued random variable with distribution $\eta(\widetilde{\xi}_{\tau_{i}-},{\rm d}\theta)$ given by
\begin{align}
     \eta(x,{\rm d}\theta)
     :=\left(\frac{c(x)}{\gamma(x)}1_{A}(x)+1_{E\setminus A}(x)\right)\delta_{0}({\rm d}\theta)
     +\frac{1}{\gamma (x)}1_{A}(x)1_{(0,+\infty)}(\theta)\theta\Pi^{NL}(x,{\rm d}\theta).\label{eta}
\end{align}
    Moreover, given $\widetilde{\xi}$ starting from $x$ (including $\{\tau_{i}:i\ge 1\}$), $\{\Theta_{i}:i\ge 1\}$ are mutually independent, $\{X^{r,i}:i\ge 1\}$ are mutually independent, also independent of $n$, $m$, $\{X^{n,t}:t\in D^{n}\}$ and $\{X^{m,t}:t\in D^{m}\}$.

 \item{\rm(v)}
 $((X_{t})_{t\ge 0};\pp_{\mu,x})$ is equal
in law to $((X_{t})_{t\ge 0};\p_{\mu})$,
 a copy of the $(\mathfrak{S}_{t},\phi^{L},\phi^{NL})$-superprocess starting from $\mu$. Moreover, $((X_{t})_{t\ge 0};\pp_{\mu,x})$ is independent of $\widetilde{\xi}$, $n$, $m$ and all the immigration processes.
\end{description}
We denote by
$$I^{c}_{t}:=\sum_{s\in D^{n}_{t}}
X^{n,s}_{t-s}
,\quad I^{d}_{t}:=\sum_{s\in D^{m}_{t}}
X^{m,s}_{t-s}
\quad\mbox{and  }\quad I^{r}_{t}:=\sum_{\tau_{i}\le t}X^{r,i}_{t-\tau_{i}}$$
the continuous immigration, the discontinuous immigration and the revival-caused immigration, respectively. We define $\Gamma_{t}$ by
$$\Gamma_{t}:=X_{t}+I^{c}_{t}+I^{d}_{t}+I^{r}_{t},\qquad\forall t\ge 0.$$
The process $\widetilde{\xi}$ is called the \textit{spine} process, and the process $I_{t}:=I^{c}_{t}+I^{d}_{t}+I^{r}_{t}$ is called the \textit{immigration} process.
\end{defi}

For any $\mu\in \me$ and
any measure $\nu$ on $(E,\mathcal{B}(E))$ with $0<\langle h,\nu\rangle <+\infty$,
we randomize the law $\pp_{\mu,x}$ by replacing the deterministic choice of $x$ with an $E$-valued random variable having distribution $h(x)\nu({\rm d}x)/\langle h,\nu\rangle$. We denote the resulting law by $\pp_{\mu,\nu}$. That is to say,
\begin{equation*}
\pp_{\mu,\nu}(\cdot):=\frac{1}{\langle h,\nu\rangle}\int_{E}\pp_{\mu,x}(\cdot)h(x)\nu({\rm d}x).
\end{equation*}
Clearly $\pp_{\mu,\delta_{x}}=\pp_{\mu,x}$.
Since the laws of $X$ and $(\widetilde{\xi},I)$
under $\pp_{\mu,\nu}$ do not depend on $\nu$ and $\mu$ respectively, we sometimes write $\pp_{\mu,\cdot}$ or $\pp_{\cdot,\nu}$.
For simplicity we also write $\pp_{\mu}$ for $\pp_{\mu,\mu}$. Here we take the convention that $\pp_{0}(\Gamma_{t}=0\ \forall t\ge 0)=1$.

For $s\ge 0$, define
\begin{equation}\label{def:Lambda-rs}
\Lambda^{m}_{s}:=\langle 1,X^{m,s}_{0}\rangle, \qquad \mbox{if }s\in D^{m}\quad
\mbox{ and }\quad\Lambda^{m}_{s}:=0\quad
\mbox{otherwise}.
\end{equation}
Then, given $\widetilde{\xi}$, $\{\Lambda^{m}_{s}, s\ge 0\}$ is a Poisson point process with characteristic measure
$\theta\Pi^{L}(\widetilde{\xi}_{s},{\rm d}\theta)$.
Let $\mathcal{G}$ be the $\sigma$-field generated by $\widetilde{\xi}$
(including $\{\tau_{i}:i\ge 1\}$),
 $\{\Theta_i: i\ge 1\}$,
$\{D^{m}_{t}: t\ge 0\}$, $\{D^{n}_{t}: t\ge 0\}$,
and $\{\Lambda^{m}_{s}, s\ge 0\}$.
\begin{proposition}\label{p:prop4.4}
For $\mu\in\me^{0}$,
$f\in\B^{+}_{b}(E)$ and $t\ge 0$, we have $\pp_{\mu}\mbox{-a.s.}$
\begin{align}
\pp_{\mu}\left[\langle f,\Gamma_{t}\rangle | \mathcal{G}\right]
&=\langle \mathfrak{P}_{t}f,\mu\rangle+\sum_{s\in D^{n}_{t}}\mathfrak{P}_{t-s}f(\widetilde{\xi}_{s})\notag\\
&\quad+\sum_{s\in D^{m}_{t}}\Lambda^{m}_{s}\mathfrak{P}_{t-s}f(\widetilde{\xi}_{s})
+\sum_{\tau_{i}\le t}\Theta_{i}\pi(\widetilde{\xi}_{\tau_{i}-},\mathfrak{P}_{t-\tau_{i}}f).
\label{sd2}
\end{align}
\end{proposition}

\proof

By \eqref{N}, we have for every $x\in E_{+}=\{x\in E: b(x)>0\}$, $f\in\B^{+}_{b}(E)$ and $t>0$,
$$
\mathbb{N}_{x}\left(\langle f,X_{t}\rangle\right)=\p_{\delta_{x}}\left(\langle f,
X_{t}\rangle\right)=\mathfrak{P}_{t}f(x).
$$
Let $D^{r}_{t}:=\{\tau_{i}:\tau_{i}\le t\}$.
Then by the definition of $\Gamma_{t}$, under $\pp_{\mu}$,
\begin{align}
\pp_{\mu}\left[\langle f,\Gamma_{t}\rangle|\mathcal{G}\right]
&=\pp_{\mu}\left(\langle f,X_{t}\rangle\right)+\sum_{s\in D^{n}_{t}}\pp_{\mu}\left[\langle f,X^{n,s}_{t-s}\rangle|\mathcal{G}\right]\nonumber\\
&\quad +\sum_{s\in D^{m}_{t}}\pp_{\mu}\left[\langle f,X^{m,s}_{t-s}\rangle|\mathcal{G}\right]+\sum_{s=\tau_{i}\in D^{r}_{t}}\pp_{\mu}\left[\langle f,X^{r,i}_{t-s}\rangle|\mathcal{G}\right]
\nonumber\\
&=\p_{\mu}\left(\langle f,X_{t}\rangle\right)+\sum_{s\in D^{n}_{t}}\mathbb{N}_{\widetilde{\xi}_{s}}\left(\langle f,X_{t-s}\rangle\right)\nonumber\\
&\quad +\sum_{s\in D^{m}_{t}}\p_{\Lambda^{m}_{s}\delta_{\widetilde{\xi}_{s}}}\left(\langle f,X_{t-s}\rangle\right)+\sum_{s=\tau_{i}\in D^{r}_{t}}\p_{\pi_{i}}
\left(\langle f,X_{t-s}\rangle\right)\nonumber\\
&=\langle \mathfrak{P}_{t}f,\mu\rangle+\sum_{s\in D^{n}_{t}}\mathfrak{P}_{t-s}f(\widetilde{\xi}_{s})\nonumber\\
&\quad +\sum_{s\in D^{m}_{t}}\Lambda^{m}_{s}\mathfrak{P}_{t-s}f(\widetilde{\xi}_{s})+\sum_{s=\tau_{i}\in D^{r}_{t}}\Theta_{i}\pi(\widetilde{\xi}_{s-},\mathfrak{P}_{t-s}f).\nonumber
\end{align}\qed

The following is our main result on the spine decomposition for the $(\mathfrak{S}_{t},\phi^{L},\phi^{NL})$-superprocess. Its proof will be given in the next subsection.

\begin{theorem}\label{them2}
Suppose that Assumptions 0-2 hold.
For every $\mu\in\me^0$, the process $((\Gamma_{t})_{t\ge 0};\pp_{\mu})$
is Markovian and has the same law as $((X_{t})_{t\ge 0};\q_{\mu})$.
\end{theorem}

\begin{remark}\rm
In the case of a purely local branching mechanism, the revival-caused immigration does not occur. To be more specific, in that case the spine runs as a copy of the $h$-transformed process $\xi^{h}$ while only continuous and discontinuous immigration occur along the spine. The concatenating procedure and the revival-cased immigration are consequences of non-local branching. Similar phenomenon has been observed in \cite{KP} for multitype continuous-state branching processes and in \cite{CRS} for multitype superdiffusions.
\end{remark}

 \begin{remark}
 \rm
 The non-local branching mechanism $\psi$ given by \eqref{1.1}-\eqref{eq:non-local part} is not the most general form that can be assumed to establish a spine decomposition. In fact,
we can establish a spine decomposition for the class of branching mechanisms developed in \cite{DGL}:
 $$\psi(x,f)=\phi^{L}(x,f(x))+\phi^{NL}(x,f)\quad\mbox{ for }x\in E,\ f\in\mathcal{B}^{+}_{b}(E),$$
 where $\phi^{L}$ takes the same form of \eqref{eq:local part} and
\begin{align}
\phi^{NL}(x,f)
=-\int_{\mathcal{P}(E)}c(x,\pi)\pi(f)G(x,{\rm d}\pi)
-\int_{\mathcal{P}(E)}\int_{(0,+\infty)}\left(1-\mathrm{e}^{-\theta\pi(f)}\right)\Pi^{NL}(x,\pi,{\rm d}\theta)
G(x,{\rm d}\pi),\notag
\end{align}
where $\mathcal{P}(E)$ denotes the space of probability measures on $E$, $G(x,{\rm d}\pi)$ is a probability measure from $E$ to $\mathcal{P}(E)$, $c(x,\pi)$ is a nonnegative bounded measurable function on $E\times \mathcal{P}(E)$, and $\theta\Pi^{NL}(x,\pi,{\rm d}\theta)$ is a bounded kernel from $E\times \mathcal{P}(E)$ to
$(0,+\infty)$.
It is easy to see that $\phi^{NL}$ has the form given in \eqref{eq:non-local part}
when  $G(x,\mathrm{d}\pi)$ is a Dirac measure on $\mathcal{P}(E)$.

To establish the spine decomposition,
one should redefine
$\gamma(x,{\rm d}y)$ as
$$\gamma(x,{\rm d}y)=\int_{\mathcal{P}(E)}r(x,\pi)\pi({\rm d}y)G(x,{\rm d}\pi),$$
where $r(x,\pi):=c(x,\pi)+\int_{(0,+\infty)}\theta\Pi^{NL}(x,\pi,{\rm d}\theta)$.
As a result, the instantaneous distribution, defined by \eqref{def:pihrs}, of the concatenation process (the spine) $\widetilde{\xi}$, changes with $\gamma$ accordingly. Regarding the spine decomposition for the above branching mechanism,
there is a new feature of the revival-caused immigration, which is described as follows:
Given the spine $\widetilde{\xi}$ (including its revival times), at each revival time $\tau_{i}$, a probability measure $\pi_{i}$ is chosen
from $\mathcal{P}(E)$, independently, according to the distribution $G^{*}(\widetilde{\xi}_{\tau_{i}-},{\rm d}\pi)$, where
$$G^{*}(x,{\rm d}\pi):=\frac{\pi(h)r(x,\pi)G(x,{\rm d}\pi)}{\gamma(x,h)}.$$
An immigration $(X^{r,i}_{s})_{s\ge 0}$ then occurs at $\tau_{i}$, and it is equal in law to the process $((X_{s})_{s\ge 0},\p_{\Theta_{i}\pi_{i}})$, where $\Theta_{i}$
is an independent $[0,+\infty)$-valued random variable with distribution $\eta(\widetilde{\xi}_{\tau_{i}-},\pi_{i},{\rm d}\theta)$ given by
\begin{align*}
\eta(x,\pi,{\rm d}\theta)&:=\left(\frac{c(x,\pi)}{r(x,\pi)}1_{\{r(x,\pi)>0\}}+1_{\{r(x,\pi)=0\}}\right)\delta_{0}({\rm d}y)\\
&\quad +\frac{\theta 1_{\{\theta\in (0,+\infty)\}}\Pi^{NL}(x,\pi,{\rm d}\theta)}{r(x,\pi)} 1_{\{r(x,\pi)>0\}}.
\end{align*}
We omit the details of the proof here for brevity.
\end{remark}

\subsection{Proof of Theorem \ref{them2}}

In this subsection, we give the proof of Theorem \ref{them2}. In order to do this, we prove a few lemmas first.

\begin{lemma}\label{lem1}
For all $x\in E$, $t\ge 0$ and $f\in\mathcal{B}^{+}_{b}(E)$,
\begin{equation*}
\pp_{\cdot,\,x}\left[\exp\left(-\langle f,I^{c}_{t}+I^{d}_{t}\rangle\right)\,|\,\widetilde{\xi}_{s}:\ 0\le s\le t\right]=\exp\left(
-\int_{0}^{t}\Phi(\widetilde{\xi}_{s},V_{t-s}f(\widetilde{\xi}_{s})){\rm d}s\right),
\end{equation*}
where $\Phi(x,\lambda):=2b(x)\lambda+\int_{(0,+\infty)}\theta\left(1-e^{-\lambda \theta}\right)\Pi^{L}(x,{\rm d}\theta)$ for $x\in E$ and $\lambda\ge 0$.
\end{lemma}
\proof
This lemma follows from an argument which is almost identical to the one leading to (59)--(60) in \cite{KLMR}.
We omit the details here.\qed

\begin{lemma}\label{lem2}
Suppose $f,l\in\mathcal{B}^{+}_{b}(E)$ and $(x,s)\mapsto g_{s}(x)$ is a non-negative locally bounded measurable function on $E\times[0,+\infty)$. For all $x\in E$ and $t>0$, let
$$e^{-w(x,t)}:=\pp_{\cdot,\,x}\left[\exp\left(-\int_{0}^{t}g_{t-s}(\widetilde{\xi}_{s}){\rm d}s-\langle f,I^{r}_{t}\rangle-l(\widetilde{\xi}_{t})\right)\right].$$
Then $u(t,x):=e^{-\lambda_{1}t}h(x)e^{-w(x,t)}$ satisfies the following integral equation:
\begin{multline}
u(t,x)=\Pi_{x}\left[e^{-l(\xi_{t})}h(\xi_{t})\right]+\Pi_{x}\big[\int_{0}^{t}{\rm d}s\big(\Phi(\xi_{s},V_{t-s}f(\xi_{s}))u(t-s,\xi_{s})\\
-\Psi(\xi_{s},V_{t-s}f,u^{t-s})
-g_{t-s}(\xi_{s})u(t-s,\xi_{s})\big)\big],
\label{lemm2.2}
\end{multline}
where $\Psi$ and $\Phi$ are defined in Proposition \ref{prop0} and Lemma \ref{lem1} respectively.
\end{lemma}

\proof Following the idea of \cite{EO}, it suffices to prove the result in the case when $g$ does not depend on the time variable. Let $\tau_{1}$ denote the first revival time of $\widetilde{\xi}$.
We have the following fundamental equation:
\begin{align*}
e^{-w(x,t)}&=\Pi_{x}^{h}\left[e_{q+g}(t)e^{-l(\xi^{h}_{t})}\right]
+\Pi^{h}_{x}\big[\int_{0}^{t}{\rm d}s\,q(\xi^{h}_{s})e_{q+g}(s)\\ &\quad\cdot\pi^{h}(\xi^{h}_{s},e^{-w_{t-s}})\int_{[0,+\infty)}e^{-\theta\pi(\xi^{h}_{s},V_{t-s}f)}\eta(\xi^{h}_{s},{\rm d}\theta)\big].\nonumber
\end{align*}
The first term corresponds to the case when $\tau_{1}\ge t$, and the second term corresponds to the case when the first revival happens at time $s\in (0,t)$.
It then follows from Fubini's theorem and \eqref{pihx} that
\begin{align*}
e^{-\lambda_{1}t}h(x)e^{-w(x,t)}
&=\Pi_{x}\left[e_{a+g}(t)h(\xi_{t})e^{-l(\xi_{t})}\right]+
\Pi_{x}\big[\int_{0}^{t}{\rm d}s\,e_{a+g}(s)q(\xi_{s})h(\xi_{s})\\
&\quad\cdot \pi^{h}(\xi_{s},e^{-\lambda_{1}(t-s)}e^{-w_{t-s}})\int_{[0,+\infty)}e^{-\theta\pi(\xi_{s},V_{t-s}f)}\eta(\xi_{s},{\rm d}\theta)\big].\nonumber
\end{align*}
We continue the above calculation by \cite[Proposition 2.9]{Li} and \eqref{eta} to get
\begin{align}
&u(t,x)\nonumber\\
&=e^{-\lambda_{1}t}h(x)e^{-w(x,t)}\nonumber\\
&=\Pi_{x}\left[h(\xi_{t})e^{-l(\xi_{t})}\right]-\Pi_{x}\left[\int_{0}^{t}(a(\xi_{s})+g(\xi_{s}))e^{-\lambda_{1}(t-s)}h(\xi_{s})e^{-w(\xi_{s},t-s)}{\rm d}s\right]\nonumber\\
&\quad+\Pi_{x}\left[\int_{0}^{t}{\rm d}s\,q(\xi_{s})h(\xi_{s})e^{-\lambda_{1}(t-s)} \pi^{h}(\xi_{s},e^{-w_{t-s}})\int_{[0,+\infty)}e^{-\theta\pi(\xi_{s},V_{t-s}f)}\eta(\xi_{s},{\rm d}\theta)\right]\nonumber\\
&=\Pi_{x}\left[h(\xi_{t})e^{-l(\xi_{t})}\right]-\Pi_{x}\left[\int_{0}^{t}(a(\xi_{s})+g(\xi_{s}))e^{-\lambda_{1}(t-s)}h(\xi_{s})e^{-w(\xi_{s},t-s)}{\rm d}s\right]\nonumber\\
&\quad+\Pi_{x}\left[\int_{0}^{t}{\rm d}s\,\pi(\xi_{s},e^{-\lambda_{1}(t-s)}h e^{-w_{t-s}})
\left(c(\xi_{s})+\int_{(0,+\infty)}r e^{-\theta\pi(\xi_{s},V_{t-s}f)}\Pi^{NL}(\xi_{s},{\rm d}\theta)\right)\right].\nonumber
\end{align}
This directly leads to \eqref{lemm2.2}.
\qed

\begin{lemma}\label{lem3}
For all $f,g\in\mathcal{B}^{+}_{b}(E)$, $\mu\in\me$, $x\in E$ and  $t\ge 0$,
\begin{equation}
\pp_{\mu,x}\left[\exp\left(-\langle f,\Gamma_{t}\rangle-g(\widetilde{\xi}_{t})\right)\right]=\frac{e^{\lambda_{1}t}}{h(x)}e^{-\langle V_{t}f,\mu\rangle}V^{he^{-g}}_{t}f(x),\label{lem3.1}
\end{equation}
where $V^{he^{-g}}_{t}f(x)$ is the unique locally bounded solution to \eqref{li1} with initial value $he^{-g}$.
\end{lemma}
\proof
Recall from Definition \ref{sd1} that $(X;\pp_{\mu,x})$ is independent of $\widetilde{\xi}$ and all the immigration processes. Moreover, given $\widetilde{\xi}$ (including $\{\tau_{i}:i\ge 1\}$), $I^{r}$ is independent of $I^{c}$ and $I^{d}$. It then follows from Lemma \ref{lem1} that
\begin{align}
&\pp_{\mu,x}\left[\exp\left(-\langle f,\Gamma_{t}\rangle-g(\widetilde{\xi}_{t})\right)\right]\nonumber\\
&=\pp_{\mu,x}\left[\exp\left(-\langle f,X_{t}\rangle-\langle f,I^{c}_{t}+I^{d}_{t}\rangle-\langle f,I^{r}_{t}\rangle-g(\widetilde{\xi}_{t})\right)\right]\nonumber\\
&=\pp_{\mu,x}\left[e^{-\langle f,X_{t}\rangle}\right]
\pp_{\mu,x}\left\{e^{-g(\widetilde{\xi}_{t})}\pp_{\mu,x}\left[\exp\left(-\langle f,I^{c}_{t}+I^{d}_{t}\rangle\right)\,|\,\widetilde{\xi}_{s}:0\le s\le t\right]\right.\nonumber\\
&\quad\quad\left.\cdot\pp_{\mu,x}\left[\exp\left(-\langle f,I^{r}_{t}\rangle\right)|\{\widetilde{\xi}_{s}:0\le s\le t\}\cup\{\tau_{i}:\tau_{i}\le t\}\right]
\right\}
\nonumber\\
&=e^{-\langle V_{t}f,\mu\rangle}\pp_{\cdot,\,x}\left[\exp\left(-\int_{0}^{t}\Phi(\widetilde{\xi}_{s},V_{t-s}f(\widetilde{\xi}_{s})){\rm d}s-\langle f,I^{r}_{t}\rangle-g(\widetilde{\xi}_{t})\right)\right].\label{lem3.2}
\end{align}
Let $v(t,x):=e^{-\lambda_{1}t}h(x)\pp_{\cdot,\,x}\left[\exp\left(-\int_{0}^{t}\Phi(\widetilde{\xi}_{s},V_{t-s}f(\widetilde{\xi}_{s})){\rm d}s-\langle f,I^{r}_{t}\rangle-g(\widetilde{\xi}_{t})\right)\right]$. One can easily verify that $(x,s)\mapsto g_{s}(x):=\Phi(x,V_{s}f(x))$ is a locally bounded function. Thus by Lemma \ref{lem2}, $v(t,x)$ is a locally bounded solution to the equation \eqref{li1} with initial value $he^{-g}$. By the uniqueness of such a solution, we have $v(t,x)=V^{he^{-g}}_{t}f(x)$. This and \eqref{lem3.2} lead to \eqref{lem3.1}.\qed

\medskip

\noindent\textbf{Proof of Theorem \ref{them2}:}
First we claim that for every $\mu\in\me^0$, $((\Gamma_{t})_{t\ge 0};\pp_{\mu})$ has the same one dimensional distribution as $((X_{t})_{t\ge 0};\q_{\mu})$. This would follow if for every  $f\in\mathcal{B}^{+}_{b}(E)$ and every $t\ge 0$,
\begin{equation}
\pp_{\mu}\left(e^{-\langle f,\Gamma_{t}\rangle}\right)=\q_{\mu}\left(e^{-\langle f,X_{t}\rangle}\right).\label{them1.1}
\end{equation}
By the definition of $\q_{\mu}$ and Proposition \ref{prop0},
\begin{equation}
\q_{\mu}\left(e^{-\langle f,X_{t}\rangle}\right)=\frac{e^{\lambda_{1}t}}{\langle h,\mu\rangle}\p_{\mu}\left[\langle h,X_{t}\rangle e^{-\langle f,X_{t}\rangle}\right]=\frac{e^{\lambda_{1}t}}{\langle h,\mu\rangle}e^{-\langle V_{t}f,\mu\rangle}\langle V^{h}_{t}f,\mu\rangle,\label{them1.3}
\end{equation}
where $V^{h}_{t}f(x)$ is the unique locally bounded solution to \eqref{li1} with initial value $h$.
By Lemma \ref{lem3}, we have
$$\pp_{\mu,x}\left[\exp\left(-\langle f,\Gamma_{t}\rangle \right)\right]=\exp\left(\lambda_{1}t-\langle V_{t}f,\mu\rangle \right)h(x)^{-1}V^{h}_{t}f(x).$$
Thus
\begin{align}
\pp_{\mu}\left(e^{-\langle f,\Gamma_{t}\rangle}\right)
&=\frac{1}{\langle h,\mu\rangle}
\int_{E}\pp_{\mu,x}\left(e^{-\langle f,\Gamma_{t}\rangle}\right)h(x)\mu({\rm d}x)\nonumber\\
&=\frac{e^{\lambda_{1}t}}{\langle h,\mu\rangle}e^{-\langle V_{t}f,\mu\rangle}\langle V^{h}_{t}f,\mu\rangle.\label{them1.4}
\end{align}
Combining  \eqref{them1.3} and \eqref{them1.4}, we get \eqref{them1.1}. It follows that for every $\mu\in\me^0$,
\begin{align}
\pp_{\mu}(\Gamma_{t}=0)&=\q_{\mu}(X_{t}=0)\notag\\
&=\frac{1}{\langle h,\mu\rangle}\p_{\mu}\left(W^{h}_{t}(X);X_{t}=0\right)=0\quad\forall t>0.\label{them1.2}
\end{align}

It remains to prove the Markov property of $((\Gamma_{t})_{t\ge 0};\pp_{\mu})$. To do this, we apply \cite[Lemma 3.3]{EO} here. Recall that $E_{\partial}=E\cup\{\partial\}$ where $\partial$ is a cemetery point. We can extend the probability measure $\pp_{\mu,x}$ onto $\mu\times \{\partial\}$ by defining  that $\pp_{\mu,\partial}(\widetilde{\xi}_{t}=\partial,\ I_{t}=0\ \forall t\ge 0)=1$ for all $\mu\in\me$.
In the remainder of this proof, we call $\mathfrak{J}$ a Markov kernel if $\mathfrak{J}$ is a map from
the measurable space  $(S,\mathcal{S})$ to the measurable space $(S',\mathcal{S}')$
such that for every $y\in S$, $\mathfrak{J}(y,\cdot)$ is
a probability measure on $(S',\mathcal{S}')$,
and for every $B\in \mathcal{S}'$, $\mathfrak{J}(\cdot,B)\in b\mathcal{S}$ the space of bounded measurable functions on $S$. The kernel $\mathfrak{J}$ will also be viewed as an operator taking $f\in b\mathcal{S}'$ to $\mathfrak{J} f\in b\mathcal{S}$ where $\mathfrak{J} f(y):=\int_{S'}f(z)\mathfrak{J}(y,{\rm d}z)$.

Clearly $((Z_{t})_{t\ge 0}:=((\Gamma_{t},\widetilde{\xi_{t}}))_{t\ge 0};\pp_{\mu,x})$ is a Markov process on $\me\times E_{\partial}$. Denote by $S_{t}$ the transition semigroup of $Z_{t}$, by $\mathfrak{K}$ the Markov kernel from $\me\times E_{\partial}$ to $\me$ induced by the projection from $\me\times E_{\partial}$ onto $\me$, and by $\mathfrak{Q}$ the Markov kernel from $\me$ to $\me\times E_{\partial}$ given by
$$\mathfrak{Q}(\nu_{1},{\rm d}(\nu_{2}\times x)):=1_{\{\nu_{1}\not=0\}}\delta_{\nu_{1}}({\rm d}\nu_{2})\times 1_{E}(x)\frac{h(x)\nu_{1}({\rm d}x)}{\langle h,\nu_{1}\rangle}+
1_{\{\nu_{1}=0\}}\delta_{0}({\rm d}\nu_{2})\times \delta_{\partial}({\rm d}x).$$
Let $R_{t}:=\mathfrak{Q} S_{t} \mathfrak{K}$ for $t\ge 0$. One can easily verify that $\mathfrak{Q}\mathfrak{K}$ is the
identity kernel on $\me$
and $R_{t}(\nu_{1},{\rm d}\nu_{2})=\pp_{\nu_{1}}\left(\Gamma_{t}\in {\rm d}\nu_{2}\right)$ for all $\nu_{1}\in\me$.
By \cite[Lemma 3.3]{EO}, $((\Gamma_{t})_{t\ge 0};\pp_{\mu})$ is Markovian
as long as
$\mathfrak{Q} S_{t}=R_{t}\mathfrak{Q}$. This would follow if for all $f,g\in\mathcal{B}^{+}_{b}(E)$ and $\nu_{1}\in\me$,
\begin{align}
&\int_{\me}\int_{\me\times E_{\partial}}e^{-\langle f,\nu_{3}\rangle -g(y)}
\mathfrak{Q}(\nu_{2},{\rm d}(\nu_{3}\times y))R_{t}(\nu_{1},{\rm d}\nu_{2})\nonumber\\
&=\int_{\me\times E_{\partial}}\int_{\me\times E_{\partial}}e^{-\langle f,\nu_{3}\rangle -g(y)}
S_{t}(\nu_{2}\times x,{\rm d}(\nu_{3}\times y))\mathfrak{Q}(\nu_{1},{\rm d}(\nu_{2}\times x)).
\label{them1.5}
\end{align}
By the above definitions, we have
\begin{alignat*}{2}
\text{LHS of }\eqref{them1.5} &=\pp_{\nu_{1}}\left[e^{-\langle f,\Gamma_{t}\rangle}\frac{\langle h e^{-g},\Gamma_{t}\rangle}{\langle h,\Gamma_{t}\rangle}1_{\{\Gamma_{t}\not=0\}}\right]+\pp_{\nu_{1}}\left(\Gamma_{t}=0\right),\\
\text{RHS of }\eqref{them1.5}&=\pp_{\nu_{1}}\left[e^{-\langle f,\Gamma_{t}\rangle-g(\widetilde{\xi}_{t})}\right]1_{\{\nu_{1}\not=0\}}+1_{\{\nu_{1}=0\}}.
\end{alignat*}
In view of \eqref{them1.2}, to show \eqref{them1.5}, it suffices to show that for all $\mu\in\me^0$ and $f,g\in\mathcal{B}^{+}_{b}(E)$,
\begin{equation}
\pp_{\mu}\left[e^{-\langle f,\Gamma_{t}\rangle-g(\widetilde{\xi}_{t})}\right]=\pp_{\mu}\left[e^{-\langle f,\Gamma_{t}\rangle}\frac{\langle h e^{-g},\Gamma_{t}\rangle}{\langle h,\Gamma_{t}\rangle}1_{\{\Gamma_{t}\not=0\}}\right].\label{them1.6}
\end{equation}
It follows from Lemma \ref{lem3} that
\begin{align}
\pp_{\mu}\left[e^{-\langle f,\Gamma_{t}\rangle-g(\widetilde{\xi}_{t})}\right]
&=\frac{1}{\langle h,\mu\rangle}\int_{E}\pp_{\mu,x}\left[e^{-\langle f,\Gamma_{t}\rangle-g(\widetilde{\xi}_{t})}\right]h(x)\mu({\rm d}x)\nonumber\\
&=\frac{e^{\lambda_{1}t}}{\langle h,\mu\rangle}e^{-\langle V_{t}f,\mu\rangle}\langle V^{he^{-g}}_{t}f,\mu\rangle,\label{them1.7}
\end{align}
where $V^{he^{-g}}_{t}f(x)$ is the unique locally bounded solution to \eqref{li1} with initial value $he^{-g}$.
On the other hand, since $(\Gamma_{t},\pp_{\mu})$ and $(X_{t},\q_{\mu})$ are identically distributed for each $t\ge 0$, we have by the definition of $\q_{\mu}$ and Proposition \ref{prop0} that
\begin{align}
\pp_{\mu}\left[e^{-\langle f,\Gamma_{t}\rangle}\frac{\langle h e^{-g},\Gamma_{t}\rangle}{\langle h,\Gamma_{t}\rangle}1_{\{\Gamma_{t}\not=0\}}\right]
&=\q_{\mu}\left[e^{-\langle f,X_{t}\rangle}\frac{\langle h e^{-g},X_{t}\rangle}{\langle h,X_{t}\rangle}1_{\{X_{t}\not=0\}}\right]\nonumber\\
&=\frac{e^{\lambda_{1}t}}{\langle h,\mu\rangle}\p_{\mu}\left[e^{-\langle f,X_{t}\rangle}\langle h e^{-g},X_{t}\rangle\right]\nonumber\\
&=\frac{e^{\lambda_{1}t}}{\langle h,\mu\rangle}e^{-\langle V_{t}f,\mu\rangle}\langle V^{he^{-g}}_{t}f,\mu\rangle.\label{them1.8}
\end{align}
Combining \eqref{them1.7} and \eqref{them1.8},  we get \eqref{them1.6}. The proof is now complete.
\qed

\section{Sufficient conditions for a
non-degenerate martingale limit}\label{Sec5}

In this section, we will give sufficient conditions for the fundamental martingale to
have a non-degenerate limit. We start with an assumption.

\medskip

\noindent\textbf{Assumption 3.}
\begin{description}
\item{(i)} Either one of the following conditions holds.
\begin{description}
\item{(1)} $a(x),\gamma(x)\in L^{2}(E,m).$
\item{(2)} The L\'evy system $(N, H)$ of $\xi$ is of the form $(N,t)$, where $N$ is given by
$$
N(x,\mathrm{d}y)=N(x,y)m(\mathrm{d}y)
$$
with $N(x,y)$ being a symmetric Borel function on $E\times E$.
The probability kernel $\pi(x,{\rm d} y)$ has a density $\pi(x,y)$ with respect to the measure $m$ such that
$$\gamma(x)\pi(x,y)=F(x,y)N(x,y)\quad\forall x,y\in E$$
for some non-negative bounded Borel function $F(x,y)$ on $E\times E$ vanishing on the diagonal.
\end{description}

\item{(ii)} $(1_{A}\pi(\cdot,h),\widehat{h})<+\infty.$
 \item{(iii)} $x\mapsto \pi(x,h)/h$ is bounded from above on $A$.
\end{description}

\medskip

It is easy to see that Assumption 3.(iii) implies Assumption 3.(ii). In this section we will
use the first two items of this assumption. In the next section
we will use items (i) and (iii)
of this assumption.
The following theorem, giving an $L\log L$ type criterion for the martingale limit to be non-degenerate, is the main result of this section.

\begin{theorem}\label{them3}
Suppose Assumptions 0--2 and 3.(i)--(ii) hold.
Let $W^{h}_{\infty}(X)$ be the almost sure limit of the non-negative martingale $W^{h}_{t}(X)$.
Suppose that
\begin{multline}\label{llogl1}
\left( \int_{(0,+\infty)}rh(\cdot)\log^{+}(rh(\cdot))\Pi^{L}(\cdot,{\rm d}r),\widehat{h}\right)\\
+\left(\int_{(0,+\infty)}r\pi(\cdot,h)\log^{+}(r\pi(\cdot,h))\Pi^{NL}(\cdot,{\rm d}r),\widehat{h}\right)
<+\infty.
\end{multline}
We have
\begin{description}
\item{(i)}
if $\lambda_{1}<0$, then $W^{h}_{t}(X)$ converges to $W^h_\infty (X)$ as $t\to \infty$ in $L^{1}(\p_{\mu})$ for every $\mu\in\me$, and $W^{h}_{\infty}(X)$ is non-degenerate in the sense that $\p_{\mu}(W^{h}_{\infty}(X)>0)>0$ for
$\mu\in\me^{0}$;
\item{(ii)} if $\lambda_{1}>0$, then $W^{h}_{\infty}(X)=0$
$\p_{\mu}$-a.s. for every $\mu\in\me$.
 \end{description}
\end{theorem}

In the remainder of this section we will assume Assumptions 0-2 hold. Additional conditions used are stated explicitly. To prove Theorem \ref{them3}, we need a few preliminary results.

\begin{proposition}\label{prop4}
Suppose Assumption 3.(i) holds. For all $f\in \mathcal{B}_{b}(E)\cap L^{2}(E,m)$ and $s,t\in (0,+\infty)$,
\begin{equation}
\lim_{s\to t}\mathfrak{P}_{s}f=\mathfrak{P}_{t}f\quad\mbox{ in }L^{2}(E,m).
\label{prop4.0}
\end{equation}
Moreover,
\begin{equation}
\mathfrak{P}_{t}f=T_{t}f\ [m] \quad\mbox{ for all }t>0,
\label{prop4.1}
\end{equation}
where $(T_{t})_{t\ge 0}$ is the semigroup
associated with the bilinear form $(\mathcal{Q},\mathcal{F})$ by \eqref{QalphaU}.

\end{proposition}
\proof
First we suppose Assumption 3.(i2) holds. We note that $\gamma(x,{\rm d}y)=\gamma(x)\pi(x,{\rm d}y)=\gamma(x)\pi(x,y)m({\rm d}y). $ Then by Assumptions 1 and 3.(i2),
$$\int_{E}F(x,y)N(x,y)m({\rm d}x)=\int_{E}\gamma(x)\pi(x,y)m({\rm d}x)\in \mathbf{K}(\xi).$$
Thus $F$ is in the class $\mathbf{J}$ defined in \cite{CS}. Let $\hat{F}(x,y):=F(y,x)$. Then
$$\int_{E}\hat{F}(x,y)N(x,y)m({\rm d}x)=\int_{E}F(y,x)N(y,x)m({\rm d}x)=\int_{E}\gamma(y)\pi(y,x)m({\rm d}x)=\gamma(y).$$
Since $\gamma$ is a bounded function on $E$, the above equation implies that $\hat{F}$ is in the class $\mathbf{J}$.
One can also show easily that the functions $\log(1+F)$ and $\log(1+\hat{F})$ are in $\mathbf{J}$.
Define
$$A_{s,t}:=-\int_{s}^{t}a(\xi_{r})\mathrm{d}r+\sum_{s<r\le t}\log\left(1+F(\xi_{r-},\xi_{r})\right)\quad \forall 0\le s<t<+\infty.$$
It follows from \cite[Theorem 4.8]{CS} (see also \cite[p. 275]{CS2}) that the semigroup corresponding to the bilinear form $(\mathcal{Q},\mathcal{F})$ defined by \eqref{QalphaU} is
$$T_{t}f(x)=\Pi_{x}\left[\mathrm{e}^{A_{0,t}}f(\xi_{t})\right]\quad \forall t\ge 0, \ x\in E,\ f\in\mathcal{B}^{+}(E).$$
Furthermore, for any $1\le p\le +\infty$, $(T_{t})_{t\ge 0}$ is semigroup on $L^{p}(E,m)$, and for $1\le p<+\infty$, $(T_{t})_{t\ge 0}$ is strongly continuous semigroup on $L^{p}(E,m)$. Similar to \cite[(2.6)]{CS}, we have
$$\mathrm{e}^{A_{0,t}}-1=-\int_{0}^{t}\mathrm{e}^{A_{s,t}}a(\xi_{s})\mathrm{d}s+\sum_{s\le t}\mathrm{e}^{A_{s,t}}F(\xi_{s-},\xi_{s}).$$
Using this, the Markov property of $\xi$, \eqref{levysys} and Assumption 3.(i2),
one can show that for any $f\in\mathcal{B}_{b}(E)$ and $x\in E$,
\begin{align*}
T_{t}f(x)&=\Pi_{x}\left[f(\xi_{t})\right]-\Pi_{x}\left[\int_{0}^{t}a(\xi_{s})T_{t-s}f(\xi_{s})\mathrm{ds}\right]
+\Pi_{x}\left[\int_{0}^{t}\int_{E}T_{t-s}f(y)F(\xi_{s},y)N(\xi_{s},y)m(\mathrm{d}y)\mathrm{ds}\right]\\
&=\Pi_{x}\left[f(\xi_{t})\right]-\Pi_{x}\left[\int_{0}^{t}a(\xi_{s})T_{t-s}f(\xi_{s})\mathrm{ds}\right]
+\Pi_{x}\left[\int_{0}^{t}\int_{E}T_{t-s}f(y)\gamma(\xi_{s})\pi(\xi_{s},y)m(\mathrm{d}y)\mathrm{ds}\right]\\
&=\Pi_{x}\left[f(\xi_{t})\right]-\Pi_{x}\left[\int_{0}^{t}a(\xi_{s})T_{t-s}f(\xi_{s})\mathrm{ds}\right]
+\Pi_{x}\left[\int_{0}^{t}\int_{E}T_{t-s}f(y)\gamma(\xi_{s},\mathrm{d}y)\mathrm{ds}\right].
\end{align*}
This implies that $T_{t}f(x)$ satisfies the integral equation \eqref{li2}. By uniqueness $T_{t}f(x)=\mathfrak{P}_{t}f(x)$, and thus we conclude the result of this proposition.

Now we suppose Assumption 3.(i1) holds.
Fix $f\in \mathcal{B}_{b}(E)\cap L^{2}(E,m)$. We first prove \eqref{prop4.0}.
Without loss of generality, we assume $0<s<t<+\infty$.
Let $F_{r}(x):=-a(x)\mathfrak{P}_{r}f(x)+\gamma(x,\mathfrak{P}_{r}f)$.
We have shown in the argument below \eqref{QalphaU} that
$\|\mathfrak{P}_{r}f\|_{\infty}\le e^{c_{1}r}\|f\|_{\infty}$
for some constant $c_{1}>0$.
Thus by definition, $|F_{r}(x)|\le (|a(x)|+\gamma(x))
\|\mathfrak{P}_{r}f\|_{\infty}\le e^{c_{1}r}\|f\|_{\infty}(|a(x)|+\gamma(x))$.
Clearly by the boundedness of $a(x)$ and $\gamma(x)$, $(x,r)\mapsto F_{r}(x)$ is locally bounded on $E\times [0,+\infty)$ and by Assumption 3.(i), $x\mapsto F_{r}(x)\in \mathcal{B}_{b}(E)\cap L^{2}(E,m)$. By \eqref{li2}, we have
\begin{align}
\mathfrak{P}_{t}f(x)-\mathfrak{P}_{s}f(x)
&=
\mathfrak{S}_{t}f(x)-\mathfrak{S}_{s}f(x)+\Pi_{x}\left[\int_{s}^{t}F_{r}(\xi_{t-r}){\rm d}r\right]\nonumber\\
&\quad +\Pi_{x}\left[\int_{0}^{s}F_{r}(\xi_{t-r})-F_{r}(\xi_{s-r}){\rm d}r\right].\label{prop4.2}
\end{align}
Recall that $\{\mathfrak{S}_{t}:t\ge 0\}$ is a strongly continuous contraction semigroup on $L^{2}(E,m)$. Thus
\begin{align}
\|\mathfrak{S}_{t}f-\mathfrak{S}_{s}f\|_{L^{2}(E,m)}&=\|\mathfrak{S}_{s}\left(\mathfrak{S}_{t-s}f-f\right)\|_{L^{2}(E,m)}\notag\\
&\le\|\mathfrak{S}_{t-s}f-f\|_{L^{2}(E,m)}\to 0\text{ as }s\to t.\label{prop4.8}
\end{align}
Note that
\begin{equation}
\left|\Pi_{x}\left[\int_{s}^{t}F_{r}(\xi_{t-r}){\rm d}r\right]\right|\le \int_{s}^{t}\left|\Pi_{x}\left[F_{r}(\xi_{t-r})\right]\right|{\rm d}r=\int_{s}^{t}\left|\mathfrak{S}_{t-r}F_{r}(x)\right|{\rm d}r.\label{prop4.3}
\end{equation}
We have by
Minkowski's integral inequality
and the contractivity of $\mathfrak{S}_{t}$ that
\begin{align}
\|\int_{s}^{t}\left|\mathfrak{S}_{t-r}F_{r}\right|{\rm d}r\|_{L^{2}(E,m)}
&\le\int_{s}^{t}
\left\|\mathfrak{S}_{t-r}F_{r}\right\|_{L^{2}(E,m)}
{\rm d}r
  \le\int_{s}^{t}\|F_{r}\|_{L^{2}(E,m)}{\rm d}r\nonumber\\
&\le\|f\|_{\infty}(\|a\|_{L^{2}(E,m)}+\|\gamma\|_{L^{2}(E,m)})\int_{s}^{t}e^{c_{1}r}{\rm d}r\to 0\nonumber
\end{align}
as $s\to t$. This together with \eqref{prop4.3} implies that
\begin{equation}
\lim_{s\to t}\left\|\Pi_{x}\left[\int_{s}^{t}F_{r}(\xi_{t-r}){\rm d}r\right]\right\|_{L^{2}(E,m)}=0.\label{prop4.4}
\end{equation}
Note that by the Markov property of $\xi$,
\begin{align}
&\left|\Pi_{x}\left[\int_{0}^{s}F_{r}(\xi_{t-r})-F_{r}(\xi_{s-r}){\rm d}r\right]\right|\nonumber\\
&\le\int_{0}^{s}\Pi_{x}\left[\left|\Pi_{\xi_{s-r}}(F_{r}(\xi_{t-s}))-\Pi_{\xi_{s-r}}(F_{r}(\xi_{0}))\right|\right]{\rm d}r\nonumber\\
&=\int_{0}^{s}\mathfrak{S}_{s-r}\left(\left|\mathfrak{S}_{t-s}F_{r}-F_{r}\right|\right)(x){\rm d}r.\label{prop4.5}
\end{align}
It follows from the strong continuity and contractivity of the semigroup $\{\mathfrak{S}_{t}:t\ge 0\}$ that
\begin{alignat*}{2}
\lim_{s\to t}\|\mathfrak{S}_{t-s}F_{r}-F_{r}\|_{L^{2}(E,m)}&=0, \text{ and }\\
\|\mathfrak{S}_{t-s}F_{r}-F_{r}\|_{L^{2}(E,m)}\le 2\|F_{r}\|_{L^{2}(E,m)}&\le 2e^{c_{1}r}\|f\|_{\infty}(\|a\|_{L^{2}(E,m)}+\|\gamma\|_{L^{2}(E,m)}).
\end{alignat*}
Thus by
Minkowski's integral inequality
and the dominated convergence theorem, we have
\begin{align}
\left\|\int_{0}^{s}\mathfrak{S}_{s-r}\left(\left|\mathfrak{S}_{t-s}F_{r}-F_{r}\right|\right){\rm d}r\right\|_{L^{2}(E,m)}
&\le\int_{0}^{s}\left\|\mathfrak{S}_{s-r}\left(\left|\mathfrak{S}_{t-s}F_{r}-F_{r}\right|\right)\right\|_{L^{2}(E,m)}{\rm d}r\nonumber\\
&\le\int_{0}^{s}\left\|\mathfrak{S}_{t-s}F_{r}-F_{r}\right\|_{L^{2}(E,m)}{\rm d}r\nonumber\\
&\to 0\quad\mbox{ as }s\to t.\nonumber
\end{align}
 This together with \eqref{prop4.5} implies that
 \begin{equation}
 \lim_{s\to t}\left\|\Pi_{x}\left[\int_{0}^{s}F_{r}(\xi_{t-r})-F_{r}(\xi_{s-r}){\rm d}r\right]\right\|_{L^{2}(E,m)}=0.\label{prop4.6}
 \end{equation}
Combining \eqref{prop4.2}--\eqref{prop4.6}, we arrive at \eqref{prop4.0}.
To prove \eqref{prop4.1}, it suffices to prove that for every $t>0$ and every $g\in L^{2}(E,m)$,
\begin{equation}
\int_{E}\mathfrak{P}_{t}f(x)g(x)m({\rm d}x)=\int_{E}T_{t}f(x)g(x)m({\rm d}x).
\label{prop4.7}
\end{equation}
Note that by H\"{o}lder's inequality and \eqref{prop4.0}, for $s,t\in (0,+\infty)$,
\begin{align}\nonumber
&\left|\int_{E}\mathfrak{P}_{t}f(x)g(x)m({\rm d}x)-\int_{E}\mathfrak{P}_{s}f(x)g(x)m({\rm d}x)\right|\nonumber\\
&\le \|\mathfrak{P}_{t}f-\mathfrak{P}_{s}f\|_{L^{2}(E,m)}\|g\|_{L^{2}(E,m)}\to 0
\end{align}
as $s\to t$. This implies $t\mapsto \int_{E}\mathfrak{P}_{t}f(x)g(x)m({\rm d}x)$
is a continuous function on $(0,+\infty)$. Similarly, using the strong continuity of $\{T_{t}:t\ge 0\}$ on $L^{2}(E,m)$, one can prove that
$t\mapsto \int_{E}T_{t}f(x)g(x)m({\rm d}x)$ is also a continuous function on $(0,+\infty)$.
By taking the Laplace transform of $\int_{E}\mathfrak{P}_{t}f(x)g(x)m({\rm d}x)$
(respectively, $\int_{E}T_{t}f(x)g(x)m({\rm d}x)$), we get $\int_{E}R_{\alpha}f(x)g(x)m({\rm d}x)$ (respectively, $\int_{E}U_{\alpha}f(x)g(x)m({\rm d}x)$). It has been shown in the argument below \eqref{QalphaU} that under Assumption 3.(i), $R_{\alpha}f=U_{\alpha}f
\ [m]$ for $\alpha$ sufficiently large. So the Laplace transforms of both sides of \eqref{prop4.7} are identical for $\alpha$ sufficiently large. Hence \eqref{prop4.7} follows from Post's inversion theorem for Laplace transforms.
\qed

\begin{proposition}\label{prop:invariant measure}
Under the assumptions of Proposition \ref{prop4},
the measure $$\rho({\rm d}x)=h(x)\widehat{h}(x)m({\rm d}x)$$ is an invariant probability measure for the semigroup $\{\widetilde{\mathfrak{S}}_{t}:t\ge 0\}$, i.e., for all $t\ge 0$ and $f\in\mathcal{B}^{+}(E)$,
\begin{equation}
\int_{E}\widetilde{\mathfrak{S}}_{t}f(x)\rho({\rm d}x)=\int_{E}f(x)\rho({\rm d}x).\label{i2}
\end{equation}
\end{proposition}
\proof
By the monotone convergence theorem, we only need to prove \eqref{i2} for $f\in\mathcal{B}^{+}_{b}(E)$.
Clearly $fh\in \mathcal{B}^{+}_{b}(E)\cap L^{2}(E,m)$. It follows by \eqref{prop2.0},\eqref{prop4.1} and \eqref{duality} that
\begin{align}
\int_{E}\widetilde{\mathfrak{S}}_{t}f(x)\rho({\rm d}x)
&=\int_{E}e^{\lambda_{1}t}\mathfrak{P}_{t}(fh)(x)\widehat{h}(x)m({\rm d}x)\nonumber\\
&=\int_{E}e^{\lambda_{1}t}T_{t}(fh)(x)\widehat{h}(x)m({\rm d}x)\nonumber\\
&=\int_{E}e^{\lambda_{1}t}f(x)h(x)\widehat{T}_{t}\widehat{h}(x)m({\rm d}x)\nonumber\\
&=\int_{E}f(x)\rho({\rm d}x).\nonumber
\end{align}
\qed

\medskip

\begin{lemma}\label{lem2.2} The function
$g(x):=h(x)^{-1}\p_{\delta_{x}}\left[W^{h}_{\infty}(X)\right]$ satisfies that
\begin{equation}
\p_{\mu}\left[W^{h}_{\infty}(X)\right]=\langle gh,\mu\rangle
\quad\mbox{ for all }\mu\in \me.\label{7.8}
\end{equation}
 Moreover,
\begin{equation}
\widetilde{\mathfrak{S}}_{t}g(x)=g(x)\quad\mbox{ for all }t\ge 0 \mbox{ and }x\in E.\label{7.9}
\end{equation}
\end{lemma}

\proof
To prove the first claim, we note that for an arbitrary constant $\lambda>0$, by the bounded convergence theorem,
\begin{align}
\p_{\mu}\left[\exp\left(-\lambda W^{h}_{\infty}(X)\right)\right]
&=\lim_{t\to+\infty}\p_{\mu}\left[\exp\left(-\lambda W^{h}_{t}(X)\right)\right]\nonumber\\
&=\lim_{t\to+\infty}\exp\left(-\langle l_{\lambda}(t, \cdot),\mu\rangle\right)\nonumber\\
&=\exp\left(-\lim_{t\to+\infty}\langle l_{\lambda}(t, \cdot),\mu\rangle\right),\label{2.2.2}
\end{align}
where $l_{\lambda}(t,x):=-\log\p_{\delta_{x}}\left[\exp\left(-\lambda W^{h}_{t}(X)\right)\right]$. Let
$$
l_{\lambda}(x):=\lim_{t\to+\infty}l_{\lambda}(t,x)=-\log\p_{\delta_{x}}\left[\exp\left(-\lambda W^{h}_{\infty}(X)\right)\right].
$$
 We have by Jensen's inequality that
$$l_{\lambda}(t,x)\le
\lambda\p_{\delta_{x}}\left(W^{h}_{t}(X)\right)=
 \lambda e^{\lambda_{1}t}\mathfrak{P}_{t}h(x)
 =\lambda h(x)\quad\mbox{ for all }x\in E,\ t\ge 0.$$
Hence $l_{\lambda}(x)\le \lambda h(x)$ for all $x\in E$. This together with
\eqref{2.2.2} and the dominated convergence theorem yields that
\begin{equation}
\p_{\mu}\left[\exp\left(-\lambda W^{h}_{\infty}(X)\right)\right]=e^{-\langle l_{\lambda},\mu\rangle }.\label{5.25}
\end{equation}
Thus we get \eqref{7.8} by differentiating both sides of \eqref{5.25} with respect to $\lambda$
and then letting $\lambda\downarrow 0$. Note that $0\le g\le 1$ by Fatou's lemma.
By the Markov property of $X$ and \eqref{7.8}, we have for all $t\ge 0$ and $x\in E$,
\begin{align}
g(x)&=\frac{1}{h(x)}\p_{\delta_{x}}\left[\lim_{s\to+\infty}e^{\lambda_{1}(t+s)}\langle h,X_{t+s}\rangle\right]\nonumber\\
&=\frac{e^{\lambda_{1}t}}{h(x)}\p_{\delta_{x}}\left[\p_{X_{t}}\left(\lim_{s\to+\infty}W^{h}_{s}(X)\right)\right]=\frac{e^{\lambda_{1}t}}{h(x)}\p_{\delta_{x}}\left[\p_{X_{t}}\left(W^{h}_{\infty}(X)\right)\right]\nonumber\\
&=\frac{e^{\lambda_{1}t}}{h(x)}\p_{\delta_{x}}\left[\langle gh,X_{t}\rangle\right]=
\frac{e^{\lambda_{1}t}}{h(x)}\mathfrak{P}_{t}(gh)(x)
=\widetilde{\mathfrak{S}}_{t}g(x).\nonumber
\end{align}
Here we used \eqref{prop2.0} in the last equality.
\qed

\medskip

\begin{lemma}\label{lem4}
Suppose Assumption 3.(i) holds.
Let $D^{m}$ and $\Lambda^m_s$ be as in Definition \ref{sd1}.(iii)
and \eqref{def:Lambda-rs} respectively.
If condition \eqref{llogl1} holds, then for $m$-almost every $x\in E$,
\begin{equation}
\lim_{D^{m}\ni s\to+\infty}\frac{\log^{+}(\Lambda^{m}_{s}h(\widetilde{\xi}_{s}))}{s}
=\lim_{i\to+\infty}\frac{\log^{+}\Theta_{i}\pi(\widetilde{\xi}_{\tau_{i}-},h)}{\tau_{i}}=0\quad\pp_{\cdot,\,x}\mbox{-a.s.}\label{7.16}
\end{equation}
\end{lemma}

\proof To prove \eqref{7.16}, it suffices to prove that for any $\eps>0$ sufficiently small,
\begin{alignat}{2}
\pp_{\cdot,\,x}\left(\sum_{s\in D^{m}}1_{\{\Lambda^{m}_{s}h(\widetilde{\xi}_{s})>e^{\eps s}\}}=+\infty\right)
\text{ and }&\pp_{\cdot,\,x}\left(\sum_{i=1}^{+\infty}1_{\{\Theta_{i}\pi(\widetilde{\xi}_{\tau_{i}-},h)>e^{\eps \tau_{i}}\}}=+\infty\right)
=0.\label{7.17}
\end{alignat}
For any $B\in\mathcal{B}(E)$ with $0<m(B)<+\infty$, let $\mu_{B}({\rm d}x):=\widehat{h}(x)1_{B}(x)m({\rm d}x)$. Clearly
$\mu_{B}\in\me^{0}$.
Recall that given $\widetilde{\xi}$, $\{\Lambda^{m}_{s}:s\ge 0\}$ is a Poisson point process with characteristic measure
$\lambda\Pi^{L}(\widetilde{\xi}_{s},{\rm d}\lambda)$.
Thus by Fubini's theorem and the fact that $\rho({\rm d}x)=h(x)\widehat{h}(x)m({\rm d}x)$ is an invariant measure for $\widetilde{\mathfrak{S}}_{t}$, we have
\begin{align}
&\pp_{\mu_{B}}\left(\sum_{s\in D^{m}}1_{\{\Lambda^{m}_{s}h(\widetilde{\xi}_{s})>e^{\varepsilon s}\}}\right)\nonumber\\
&=\pp_{\mu_{B}}\left(\int_{0}^{+\infty}\int_{(0,+\infty)}\lambda1_{\{\lambda h(\widetilde{\xi}_{s})>e^{\varepsilon s}\}}\Pi^{L}(\widetilde{\xi}_{s},{\rm d}\lambda){\rm d}s\right)\nonumber\\
&=\frac{1}{\langle h,\mu_{B}\rangle}\int_{B}\pp_{\mu_{B},x}\left(\int_{0}^{+\infty}\int_{(0,+\infty)}\lambda 1_{\{\lambda h(\widetilde{\xi}_{s})>e^{\varepsilon s}\}}\Pi^{L}(\widetilde{\xi}_{s},{\rm d}\lambda){\rm d}s\right)\rho({\rm d}x)\nonumber\\
&\le\frac{1}{\langle h,\mu_{B}\rangle}\int_{0}^{+\infty}{\rm d}s\int_{E}\rho({\rm d}x)
\int_{(0,+\infty)}\lambda 1_{\{\lambda h(x)>e^{\varepsilon s}\}}\Pi^{L}(x,{\rm d}\lambda)\nonumber\\
&=\frac{1}{\langle h,\mu_{B}\rangle}\left(\int_{(0,+\infty)}\lambda\Pi^{L}(\cdot,{\rm d}\lambda )\int_{0}^{\log^{+}\lambda h(x)/\varepsilon}{\rm d}s,h\widehat{h}\right)\nonumber\\
&=\frac{1}{\varepsilon\langle h,\mu_{B}\rangle}\left(\int_{(0,+\infty)}\lambda h(\cdot)\log^{+}(\lambda h(\cdot))\Pi^{L}(\cdot,{\rm d}\lambda),\widehat{h}\right).\label{7.18}
\end{align}
The right hand side of \eqref{7.18} is finite by \eqref{llogl1}.
Thus we get
$$\pp_{\mu_{B}}\left(\sum_{s\in D^{m}}1_{\{\Lambda^{m}_{s}h(\widetilde{\xi}_{s})>e^{\varepsilon s}\}}<+\infty\right)=1.$$
Note that
\begin{align}\nonumber
&\pp_{\mu_{B}}\left(\sum_{s\in D^{m}}1_{\{\Lambda^{m}_{s}h(\widetilde{\xi}_{s})>e^{\varepsilon s}\}}<+\infty\right)
=\rho(B)^{-1}\int_{B}\pp_{\cdot,\,x}\left(\sum_{s\in D^{m}}1_{\{\Lambda^{m}_{s}h(\widetilde{\xi}_{s})>e^{\varepsilon s}\}}<+\infty\right)\rho({\rm d}x).
\end{align}
Thus $\pp_{\cdot,\,x}\left(\sum_{s\in D^{m}}1_{\{\Lambda^{m}_{s}h(\widetilde{\xi}_{s})>e^{\varepsilon s}\}}<+\infty\right)=1$ for $m$-almost every $x\in B$. Since $B$ is arbitrary, the first equality of \eqref{7.17} holds for $m$-almost every $x\in E$.

Recall from Definition \ref{sd1} that given $\widetilde{\xi}$ (including $\{\tau_{i}:i\ge 1\}$), $\Theta_{i}$ is distributed as $\eta(\widetilde{\xi}_{\tau_{i}-},{\rm d}\theta)$ given by \eqref{eta}. Thus by Fubini's theorem and \eqref{l1},
\begin{align}
&\pp_{\mu_{B}}\left(\sum_{i=1}^{+\infty}1_{\{\Theta_{i}\pi(\widetilde{\xi}_{\tau_{i}-},h)
> e^{\eps \tau_{i}} \}}\right)\nonumber\\
&=\pp_{\mu_{B}}\left[\sum_{i=1}^{+\infty}\int_{\theta\pi(\widetilde{\xi}_{\tau_{i}-},h)
> e^{\eps \tau_{i}} }\eta(\widetilde{\xi}_{\tau_{i}-},{\rm d}\theta)\right]\nonumber\\
&=\pp_{\mu_{B}}\left[\sum_{i=1}^{+\infty}\frac{1}{\gamma(\widetilde{\xi}_{\tau_{i}-})}1_{A}(\widetilde{\xi}_{\tau_{i}-})\int_{\theta\pi(\widetilde{\xi}_{\tau_{i}-},h)
> e^{\eps \tau_{i}} }\theta\Pi^{NL}(\widetilde{\xi}_{\tau_{i}-},{\rm d}\theta)\right]\nonumber\\
&=\pp_{\mu_{B}}\left[\int_{0}^{+\infty}q(\widetilde{\xi_{s}})\frac{1}{\gamma(\widetilde{\xi}_{s})}1_{A}(\widetilde{\xi}_{s}){\rm d}s
\int_{\theta\pi(\widetilde{\xi}_{s},h)
> e^{\eps s} }\theta\Pi^{NL}(\widetilde{\xi}_{s},{\rm d}\theta)\right]\nonumber\\
&=\frac{1}{\langle h,\mu_{B}\rangle}\int_{0}^{+\infty}{\rm d}s\int_{B}\pp_{\mu,x}\left[\frac{\pi(\widetilde{\xi}_{s},h)}{h(\widetilde{\xi}_{s})}
1_{A}(\widetilde{\xi}_{s})\int_{\theta\pi(\widetilde{\xi}_{s},h)
> e^{\eps s} }\theta\Pi^{NL}(\widetilde{\xi}_{s},{\rm d}\theta)\right]\rho({\rm d}x)\nonumber\\
&\le\frac{1}{\langle h,\mu_{B}\rangle}\int_{0}^{+\infty}{\rm d}s\int_{E}1_{A}(x)\pi(x,h)\widehat{h}(x)m({\rm d}x)\int_{\theta\pi(x,h)>e^{\eps s}}\theta\Pi^{NL}(x,{\rm d}\theta)\nonumber\\
&=\frac{1}{\langle h,\mu_{B}\rangle}\int_{E}1_{A}(x)\pi(x,h)\widehat{h}(x)m({\rm d}x)\int_{(0,+\infty)}\theta\Pi^{NL}(x,{\rm d}\theta)
\int_{0}^{\frac{\log^{+}(\theta\pi(x,h))}{\eps}}{\rm d}s\nonumber\\
&=\frac{1}{\eps\langle h,\mu_{B}\rangle}\left( \pi(\cdot,h)\int_{(0,+\infty)}\theta\log^{+}(\theta\pi(\cdot,h))\Pi^{NL}(\cdot,{\rm d}\theta),1_{A}\widehat{h}\right).\label{7.19}
\end{align}
The right hand side of \eqref{7.19} is finite by \eqref{llogl1}.
Thus we get
$$\pp_{\mu_{B}}\left(\sum_{i=1}^{+\infty}1_{\{\Theta_{i}\pi(\widetilde{\xi}_{\tau_{i}-},h)
> e^{\eps \tau_{i}} \}}<+\infty\right)=1.$$
Using
an argument similar to that at the end of the first paragraph of this proof,
one can prove that the second equality of \eqref{7.17} holds for $m$-almost every $x\in E$.\qed

\medskip

\noindent\textbf{Proof of Theorem \ref{them3}:}
Recall that, by assumption, Assumptions 0-2 and 3.(i)-(ii) hold.

\noindent (i) Suppose $\lambda_{1}<0$.
Without loss of generality, we assume $\mu\in\me^0$.
Since $W^{h}_{t}(X)$ is a non-negative martingale,
to show it is a closed martingale,
it suffices to prove
 \begin{equation}\label{4.7}
\p_{\mu}\left[W^{h}_{\infty}(X)\right]=\langle h,\mu\rangle.
\end{equation}
 First we claim that \eqref{4.7} is true for
 $
 \mu_{B}({\rm d}y):=
 1_{B}(y)
 \widehat{h}(y)
 m({\rm d}y)$ with $B\in\mathcal{B}(E)$ and $0<m(B)<+\infty$.
 It is
 straightforward
 to see from the change of measure methodology
(see, for example, \cite[Theorem 5.3.3]{Durrett})
 that the proof for this claim is complete as soon as we can show that
\begin{equation}
\q_{\mu_{B}}\left(\limsup_{t\to+\infty}W^{h}_{t}(X)<+\infty\right)=1.\label{4.8}
\end{equation}
Since $((X_{t})_{t\ge 0};\q_{\mu_{B}})$ is
equal in law
to $((\Gamma_{t})_{t\ge 0};\pp_{\mu_{B}})$, \eqref{4.8} is equivalent to that
\begin{equation}
\pp_{\mu_{B}}\left(\limsup_{t\to+\infty}W^{h}_{t}(\Gamma)<+\infty\right)=1.\label{4.9}
\end{equation}
In the remainder of this proof, we define a function $\log^{*}\theta:=\theta/e$ if $\theta\le e$ and $\log^{*}\theta:= \log\theta$ if $\theta>e$.
Under the assumptions of Theorem \ref{them3}, one can prove
by elementary computation that \eqref{llogl1} implies
\begin{multline}
\left( \int_{(0,+\infty)}rh(\cdot)\log^{*}(rh(\cdot))\Pi^{L}(\cdot,{\rm d}r),\widehat{h}\right)\\
+\left( \int_{(0,+\infty)}r\pi(\cdot,h)\log^{*}(r\pi(\cdot,h))\Pi^{NL}(\cdot,{\rm d}r),1_{A}\widehat{h}\right)<+\infty.\label{llogl1'}
\end{multline}
Recall that $\mathcal{G}$ is the $\sigma$-field generated by $\widetilde{\xi}$ (including $\{\tau_{i}:i\ge 1\}$),
$\{D^{m}_{t}: t\ge 0\}$, $\{D^{n}_{t}: t\ge 0\}$,
$\{\Theta_{i}:i\ge 1\}$ and $\{\Lambda^{m}_{s}:s\ge 0\}$.
By \eqref{sd2}, for any $t>0$,
\begin{align}
&\pp_{\mu_{B}}\left(W^{h}_{t}(\Gamma)\,|\,\mathcal{G}\right)\nonumber\\
&=e^{\lambda_{1}t}\left(\langle \mathfrak{P}_{t}h,\mu_{B}\rangle +\sum_{s\in D^{n}_{t}}\mathfrak{P}_{t-s}h(\widetilde{\xi}_{s})
+\sum_{s\in D^{m}_{t}}\Lambda^{m}_{s}\mathfrak{P}_{t-s}h(\widetilde{\xi}_{s})
+\sum_{\tau_{i}\le t}\Theta_{i}\pi(\widetilde{\xi}_{\tau_{i}-},\mathfrak{P}_{t-\tau_{i}}h)\right)
\nonumber\\
&=\langle h,\mu_{B}\rangle +\sum_{s\in D^{n}_{t}}e^{\lambda_{1}s}h(\widetilde{\xi}_{s})+\sum_{s\in D^{m}_{t}}e^{\lambda_{1}s}\Lambda^{m}_{s}h(\widetilde{\xi}_{s})+\sum_{\tau_{i}\le t}e^{\lambda_{1}\tau_{i}}\Theta_{i}\pi(\widetilde{\xi}_{\tau_{i}-},,h)\nonumber\\
&\le  \langle h,\mu_{B}\rangle +\sum_{s\in D^{n}}e^{\lambda_{1}s}h(\widetilde{\xi}_{s})+\sum_{s\in D^{m}}e^{\lambda_{1}s}\Lambda^{m}_{s}h(\widetilde{\xi}_{s})+\sum_{i=1}^{+\infty}e^{\lambda_{1}\tau_{i}}\Theta_{i}\pi(\widetilde{\xi}_{\tau_{i}-},,h).\label{5.8}
\end{align}
We begin with the second term on the right hand side of \eqref{5.8}. Let $\varepsilon \in (0,-\lambda_{1})$ be an arbitrary constant.
\begin{equation}
\sum_{s\in D^{n}}e^{\lambda_{1}s}h(\widetilde{\xi}_{s})=\sum_{s\in D^{n}}e^{\lambda_{1}s}h(\widetilde{\xi}_{s})1_{\{h(\widetilde{\xi}_{s})>e^{\varepsilon s}\}}+\sum_{s\in D^{n}}e^{\lambda_{1}s}h(\widetilde{\xi}_{s})1_{\{h(\widetilde{\xi}_{s})\le e^{\varepsilon s}\}}=:\mathrm{I}+\mathrm{II}.\nonumber
\end{equation}
Recall that given $\widetilde{\xi}$,
the random measure $\sum_{s\in D^{n}}\delta_{s}(\cdot)$ on $[0,+\infty)$ is a Poisson random measure with intensity $2 b(\widetilde{\xi}_{t}){\rm d}t$,
and that $\rho({\rm d}x)=h(x)\widehat{h}(x)m({\rm d}x)$ is an invariant probability measure for $\widetilde{\mathfrak{S}}_{t}$. We have by Fubini's theorem,
\begin{align}
\pp_{\mu_{B}}\left(\mathrm{II}\right)
&=\pp_{\mu_{B}}\left[
\pp_{\mu_{B}}\left(\left.\sum_{s\in D^{n}}e^{\lambda_{1}s}h(\widetilde{\xi}_{s})1_{\{h(\widetilde{\xi}_{s})\le e^{\varepsilon s}\}}
\right|
\widetilde{\xi}_{r}:r\ge 0\right)\right]\nonumber\\
&=\pp_{\mu_{B}}\left(\int_{0}^{+\infty}2 b(\widetilde{\xi}_{s})e^{\lambda_{1}s}h(\widetilde{\xi}_{s})1_{\{h(\widetilde{\xi}_{s})\le e^{\varepsilon s}\}}{\rm d}s\right)\nonumber\\
&=\frac{2}{\langle h,\mu_{B}\rangle}\int_{0}^{+\infty}e^{\lambda_{1}s}{\rm d}s\int_{B}\pp_{\mu,x}\left[b(\widetilde{\xi}_{s})h(\widetilde{\xi}_{s})
1\{h(\widetilde{\xi}_{s})\le e^{\eps s}\}\right]\rho({\rm d}x)\nonumber\\
&\le \frac{2}{\langle h,\mu_{B}\rangle}\int_{0}^{+\infty}e^{\lambda_{1}s}{\rm d}s\int_{E} b(x)h(x)1_{\{h(x)\le e^{\varepsilon s}\}}\rho({\rm d}x)\nonumber\\
&\le \frac{2\|b\|_{\infty}}{\langle h,\mu_{B}\rangle}\int_{0}^{+\infty}e^{(\lambda_{1}+\varepsilon)s}{\rm d}s\int_{E}\rho({\rm d}x)<+\infty.\nonumber
\end{align}
Thus we have $\pp_{\mu_{B}}\left(\mathrm{II}<+\infty\right)=1$. On the other hand,
\begin{align}
\pp_{\mu_{B}}\left(\sum_{s\in D^{n}}1_{\{h(\widetilde{\xi}_{s})>e^{\varepsilon s}\}}\right)
&=\pp_{\mu_{B}}\left(\int_{0}^{+\infty}2 b(\widetilde{\xi}_{s})1_{\{h(\widetilde{\xi}_{s})>e^{\varepsilon s}\}}{\rm d}s\right)\nonumber\\
&=\frac{2}{\langle h,\mu_{B}\rangle}\int_{0}^{+\infty}{\rm d}s\int_{B}\pp_{\mu,x}\left[b(\widetilde{\xi}_{s})1_{\{h(\widetilde{\xi}_{s})>e^{\eps s}\}}\right]\rho({\rm d}x)\nonumber\\
&\le \frac{2}{\langle h,\mu_{B}\rangle}\int_{0}^{+\infty}{\rm d}s\int_{E} b(x)1_{\{h(x)>e^{\varepsilon s}\}}
\rho({\rm d}x)\nonumber\\
&=\frac{2}{\langle h,\mu_{B}\rangle}\int_{E} b(x)\rho({\rm d}x)\int_{0}^{\frac{\log^{+}h(x)}{\varepsilon}}{\rm d}s\nonumber\\
&\le \frac{2}{\varepsilon \langle h,\mu_{B}\rangle}\|b\|_{\infty}\log^{+}\|h\|_{\infty}<+\infty.\nonumber
\end{align}
This implies that $\mathrm{I}$ is the sum of finitely many terms. Thus $\pp_{\mu_{B}}\left(\mathrm{I}<+\infty\right)=1$.
For the third term in \eqref{5.8}, we have
\begin{align}
&\sum_{s\in D^{m}}e^{\lambda_{1}s}\Lambda^{m}_{s}h(\widetilde{\xi}_{s})\nonumber\\
&=\sum_{s\in D^{m}}e^{\lambda_{1}s}\Lambda^{m}_{s}h(\widetilde{\xi}_{s})1_{\{\Lambda^{m}_{s}h(\widetilde{\xi}_{s})\le e^{\varepsilon s}\}}+\sum_{s\in D^{m}}e^{\lambda_{1}s}\Lambda^{m}_{s}h(\widetilde{\xi}_{s})1_{\{\Lambda^{m}_{s}h(\widetilde{\xi}_{s})> e^{\varepsilon s}\}}\nonumber\\
&=:\mathrm{III}+\mathrm{IV}.\nonumber
\end{align}
In view of Definition \ref{sd1}.(iii), for $\mathrm{III}$, we have
\begin{align}
&\pp_{\mu_{B}}\left(\mathrm{III}\right)\nonumber\\
&=\pp_{\mu_{B}}\left(\int_{0}^{+\infty}\int_{E}e^{\lambda_{1}s}r^{2}h(\widetilde{\xi}_{s})1_{\{rh(\widetilde{\xi}_{s})\le e^{\varepsilon s}\}}\Pi^{L}(\widetilde{\xi}_{s},{\rm d}r){\rm d}s\right)\nonumber\\
&=\frac{1}{\langle h,\mu_{B}\rangle}\int_{0}^{+\infty}e^{\lambda_{1}s}{\rm d}s\int_{B}\pp_{\mu,x}\left(\int_{E}r^{2}h(\widetilde{\xi}_{s})1_{\{rh(\widetilde{\xi}_{s})\le e^{\varepsilon s}\}}\Pi^{L}(\widetilde{\xi}_{s},{\rm d}r)\right)\rho({\rm d}x)\nonumber\\
&\le\frac{1}{\langle h,\mu_{B}\rangle}\int_{0}^{+\infty}e^{\lambda_{1}s}{\rm d}s\int_{E}\rho({\rm d}x)\int_{E}r^{2}h(x)1_{\{rh(x)\le e^{\varepsilon s}\}}\Pi^{L}(x,{\rm d}r)\nonumber\\
&=\frac{1}{\langle h,\mu_{B}\rangle}\int_{E}\rho({\rm d}x) \int_{(0,+\infty)}r^{2}h(\cdot)\Pi^{L}(\cdot,{\rm d}r)\int_{\log^{+}rh(\cdot)/\varepsilon}^{+\infty}e^{\lambda_{1}s}{\rm d}s\nonumber\\
&=\frac{-1}{\lambda_{1}\langle h,\mu_{B}\rangle}\,\big(\widehat{h},\,\int_{(0,+\infty)}r^{2}h^{2}(\cdot)\left(rh(\cdot)\vee 1\right)^{\lambda_{1}/\varepsilon}\Pi^{L}(\cdot,{\rm d}r)\big)\nonumber\\
&=\frac{-1}{\lambda_{1}\langle h,\mu_{B}\rangle}\,\left(\widehat{h},\,\int_{(0,+\infty)}\Pi^{L}(\cdot,{\rm d}r)rh(\cdot)\log^{*}(rh(\cdot))
\cdot\left(\frac{rh(\cdot)}{(rh(\cdot)\vee 1)^{-\lambda_{1}/\eps}\log^{*}(rh(\cdot))}\right)\right).
\label{4.13}
\end{align}
Note that the function $r\mapsto \frac{r}{(r\vee 1)^{-\lambda_{1}/\varepsilon}\log^{*}r}$ is bounded from above on $(0,+\infty)$. This together with \eqref{llogl1'} implies that the right hand side of \eqref{4.13} is finite.
It follows that
$\pp_{\mu_{B}}(\mathrm{III}<+\infty)=1$. It has been shown by \eqref{7.18} that $\pp_{\mu_{B}}\left(\sum_{s\in D^{m}}1_{\{
\Lambda^{m}_{s}
h(\widetilde{\xi}_{s})>e^{\varepsilon s}\}}<+\infty\right)=1$.
This implies that $\mathrm{IV}$ is the sum of finitely many terms. Thus we have $\pp_{\mu_{B}}(\mathrm{IV}<+\infty)=1$.
The fourth term on the right hand side of \eqref{5.8} can be dealt with similarly. In fact, we have
\begin{align}
\sum_{i=1}^{+\infty}e^{\lambda_{1}\tau_{i}}\Theta_{i}\pi(\widetilde{\xi}_{\tau_{i}-},h)
&=\sum_{i=1}^{+\infty}e^{\lambda_{1}\tau_{i}}\Theta_{i}\pi(\widetilde{\xi}_{\tau_{i}-},h)1_{\{\Theta_{i}\pi(\widetilde{\xi}_{\tau_{i}-},h)
\le e^{\eps \tau_{i}} \}}\nonumber\\
&\quad+\sum_{i=1}^{+\infty}e^{\lambda_{1}\tau_{i}}\Theta_{i}\pi(\widetilde{\xi}_{\tau_{i}-},h)1_{\{\Theta_{i}\pi(\widetilde{\xi}_{\tau_{i}-},h)
> e^{\eps \tau_{i}} \}}\nonumber\\
&=:\mathrm{V}+\mathrm{VI}.\nonumber
\end{align}
Recall that given $\widetilde{\xi}$ (including $\{\tau_{i}:i\ge 1\}$), $\Theta_{i}$ is distributed according to
$\eta(\widetilde{\xi}_{\tau_{i}-},{\rm d}\theta)$
given by \eqref{eta}. Thus by Fubini's theorem and \eqref{l1},
\begin{align}
&\pp_{\mu_{B}}(\mathrm{V})\nonumber\\
&=\pp_{\mu_{B}}\left[\sum_{i=1}^{+\infty}e^{\lambda_{1}\tau_{i}}\pi(\widetilde{\xi}_{\tau_{i}-},h)
\int_{[0,+\infty)}\theta 1_{\{\theta\pi(\widetilde{\xi}_{\tau_{i}-},h)\le e^{\eps \tau_{i}}\}}\eta(\widetilde{\xi}_{\tau_{i}-},{\rm d}\theta)\right]\nonumber\\
&=\pp_{\mu_{B}}\left[\sum_{i=1}^{+\infty}e^{\lambda_{1}\tau_{i}}\frac{\pi(\widetilde{\xi}_{\tau_{i}-},h)}{\gamma(\widetilde{\xi}
_{\tau_{i}-})}1_{A}(\widetilde{\xi}_{\tau_{i}-})
\int_{\{0<\theta\pi(\widetilde{\xi}_{\tau_{i}-},h)\le e^{\eps \tau_{i}}\}}\theta^{2}\Pi^{NL}(\widetilde{\xi}_{\tau_{i}-},{\rm d}\theta)\right]\nonumber\\
&=\frac{1}{\langle h,\mu_{B}\rangle}\int_{B}\rho({\rm d}x)\pp_{\mu,x}\big[\int_{0}^{+\infty}q(\widetilde{\xi}_{s})e^{\lambda_{1}s}\frac{\pi(\widetilde{\xi}_{s},h)}{\gamma(\widetilde{\xi}
_{s})}1_{A}(\widetilde{\xi}_{s}){\rm d}s\nonumber\\
&\quad\quad\quad\quad\quad\cdot\int_{\{0<\theta\pi(\widetilde{\xi}_{s},h)\le e^{\eps s}\}}\theta^{2}\Pi^{NL}(\widetilde{\xi}_{s},{\rm d}\theta)\big]\nonumber\\
&=\frac{1}{\langle h,\mu_{B}\rangle}\int_{0}^{+\infty}e^{\lambda_{1}s}{\rm d}s\int_{B}\rho({\rm d}x)
\pp_{\mu,x}\big[\frac{\pi(\widetilde{\xi}_{s},h)^{2}}{h(\widetilde{\xi_{s}})}1_{A}(\widetilde{\xi}_{s})\notag\\
&\quad\quad\quad\quad\quad\cdot\int_{\{0<\theta\pi(\widetilde{\xi}_{s},h)\le e^{\eps s}\}}\theta^{2}\Pi^{NL}(\widetilde{\xi}_{s},{\rm d}\theta)\big]\nonumber\\
&\le\frac{1}{\langle h,\mu_{B}\rangle}\int_{0}^{+\infty}e^{\lambda_{1}s}{\rm d}s\int_{E}1_{A}(x)\pi(x,h)^{2}\widehat{h}(x)m({\rm d}x)\int_{\{0<\theta\pi(x,h)\le e^{\eps s}\}}\theta^{2}\Pi^{NL}(x,{\rm d}\theta)\nonumber\\
&=\frac{1}{\langle h,\mu_{B}\rangle}\int_{E}1_{A}(x)\pi(x,h)^{2}\widehat{h}(x)m({\rm d}x)\int_{(0,+\infty)}\theta^{2}\Pi^{NL}(x,{\rm d}\theta)
\int_{\frac{\log^{+}(\theta\pi(x,h))}{\eps}}^{+\infty}e^{\lambda_{1}s}{\rm d}s\nonumber\\
&=\frac{-1}{\lambda_{1}\langle h,\mu_{B}\rangle}\big(\int_{(0,+\infty)}1_{A}(\cdot)\pi(\cdot,h)^{2}\theta^{2}\left(\pi(\cdot,h)\theta\vee 1\right)^{\lambda_{1}/\eps}\Pi^{NL}(\cdot,{\rm d}\theta),\widehat{h}\big)\nonumber\\
&=\frac{-1}{\lambda_{1}\langle h,\mu_{B}\rangle}\big(\int_{(0,+\infty)}\Pi^{NL}(\cdot,{\rm d}\theta)1_{A}(\cdot)\pi(\cdot,h)\theta\log^{*}(\pi(\cdot,h)\theta)
\left(\frac{\pi(\cdot,h)\theta}{
(\pi(\cdot,h)\theta\vee 1)^{-\lambda_{1}/\eps}\log^{*}(\pi(\cdot,h)\theta)}\right),\widehat{h}\big).\nonumber
\end{align}
Since $\theta\mapsto \frac{\theta}{(\theta\vee 1)^{-\lambda_{1}/\eps}\log^{*}\theta}$ is bounded from above on $(0,+\infty)$, we get $\pp_{\mu_{B}}(\mathrm{V})<+\infty$ by \eqref{llogl1'}, and hence $\pp_{\mu_{B}}(\mathrm{V}<+\infty)=1$.
We have shown in \eqref{7.19} that $$\pp_{\mu_{B}}\left(\sum_{i=1}^{+\infty}1_{\{\Theta_{i}\pi(\widetilde{\xi}_{\tau_{i}-},h)
> e^{\eps \tau_{i}} \}}<+\infty\right)=1.$$
Thus $\mathrm{VI}$ is the sum of finitely many terms and  $\pp_{\mu_{B}}\left(\mathrm{VI}<+\infty\right)=1$.
The above arguments show that the right hand side of \eqref{5.8}
is finite a.s., and hence $\limsup_{t\to +\infty}\pp_{\mu_{B}}\left(W^{h}_{t}(\Gamma)\,|\,\mathcal{G}\right)<+\infty$
$\pp_{\mu_{B}}$-a.s. By Fatou's lemma, $\pp_{\mu_{B}}\left(\liminf_{t\to+\infty}W^{h}_{t}(\Gamma)\,|\,\mathcal{G}\right)<+\infty$ $\pp_{\mu_{B}}$-a.s. Let $$A_{n}:=
\left\{
\pp_{\mu_{B}}\left(\liminf_{t\to+\infty}W^{h}_{t}(\Gamma)\,|\,\mathcal{G}\right)\le n
\right\}
\in \mathcal{G}\quad\mbox{ for }n\ge 1.$$ Then $\pp_{\mu_{B}}(\cup_{n=1}^{+\infty}A_{n})=1$.
 Since
 $$\int_{A_{n}}\liminf_{t\to+\infty}W^{h}_{t}(\Gamma)d\pp_{\mu_{B}}
 =\int_{A_{n}}\pp_{\mu_{B}}\left(\liminf_{t\to+\infty}W^{h}_{t}(\Gamma)\,|\,\mathcal{G}\right)d\pp_{\mu_{B}}\le n,$$
 we get $\liminf_{t\to+\infty}W^{h}_{t}(\Gamma)<+\infty$ $\pp_{\mu_{B}}$-a.s on $A_{n}$ for all $n\ge 1$. Thus
$$\pp_{\mu_{B}}\left(\liminf_{t\to+\infty}W^{h}_{t}(\Gamma)<+\infty\right)=1.$$
Note that by \cite[Proposition 2]{HR} $W^{h}_{t}(\Gamma)^{-1}$
 is a non-negative $\pp_{\mu_{B}}$-supermartingale,
  which implies that  $\lim_{t\to+\infty}W^{h}_{t}(\Gamma)^{-1}$ exists $\pp_{\mu_{B}}$-a.s.
 It follows that
 $$\pp_{\mu_{B}}\left(\limsup_{t\to+\infty}W^{h}_{t}(\Gamma)<+\infty\right)=1.$$
This proves \eqref{4.9} and consequently
$\p_{\mu_{B}}\left[W^{h}_{\infty}(X)\right]=\langle h,\mu_{B}\rangle.$
 Recall that $\p_{\mu_{B}}\left[W^{h}_{\infty}(X)\right]=\langle g h,\mu_{B}\rangle$ where
 $g(x)=h(x)^{-1}\p_{\delta_{x}}\left[W^{h}_{\infty}(X)\right]$.
 Thus we have
  \begin{equation}\label{4.15'}
  \langle g h,\mu_{B}\rangle=\langle h,\mu_{B}\rangle.\end{equation}
 Note that $0\le g(x)\le 1 $ for every $x\in E$.
 We get by
 \eqref{4.15'}
  that $g(x)=1$ $m$-a.e. on $B$. Since $B$ is arbitrary, $g(x)=1$ $m$-a.e. on $E$. It then follows from \eqref{7.9} that $g(x)=\widetilde{\mathfrak{S}}_{t}g(x)=\int_{E}\widetilde{p}(t,x,y)g(y)\rho({\rm d}y)=1$ for every $x\in E$. Therefore by \eqref{7.8}, $\p_{\mu}\left[W^{h}_{\infty}(X)\right]=\langle h,\mu\rangle$ holds for all $\mu\in\me$. This completes the proof for Theorem \ref{them3}.(i).

\medskip

\noindent(ii) Suppose $\lambda_{1}>0$.
Clearly $\p_{\mu}\left(W^{h}_{\infty}(X)=0\right)=1$ if and only if $\p_{\mu}\left[W^{h}_{\infty}(X)\right]=0$.
By \eqref{7.8}, this would follow if $g(x)=0$ for every $x\in E$. Recall that $g(x)=\widetilde{\mathfrak{S}}_{t}g(x)=\int_{E}\widetilde{p}(t,x,y)g(y)\rho({\rm d}y)$. It suffices to prove that $g(x)=0$ for $m$-almost every $x\in E$, or equivalently,
\begin{equation}
\p_{\delta_{x}}
\left[
W^{h}_{\infty}(X)
\right]
=0\quad\mbox{ for $m$-a.e. $x\in E$.}\label{7.28}
\end{equation}
By the change of measure methodology (see, for example,
\cite[Theorem 5.3.3]{Durrett}),
\eqref{7.28} would follow if
\begin{equation}
\pp_{\delta_{x}}(\limsup_{t\to+\infty}W^{h}_{t}(\Gamma)=+\infty)=1\quad\mbox{ for $m$-a.e. $x\in E$.}\label{7.29}
\end{equation}
By the definition of $\Gamma_{t}$, we have
$$W^{h}_{s}(\Gamma)=e^{\lambda_{1}s}\langle h,\Gamma_{s}\rangle\ge e^{\lambda_{1}s}
\Lambda^{m}_{s}
h(\widetilde{\xi}_{s})\quad\mbox{ for }s\in D^{m},$$
and
$$W^{h}_{\tau_{i}}(\Gamma)\ge e^{\lambda_{1}\tau_{i}}
\Theta_{i}
\pi(\widetilde{\xi}_{\tau_{i}-},h)\quad\mbox{ for }i\ge 1.$$
Thus under $\pp_{\delta_{x}}$,
\begin{multline}\label{7.30}
\limsup_{t\to+\infty}W^{h}_{t}(\Gamma)
\ge \limsup_{D^{m}\ni s\to+\infty}e^{\lambda_{1}s}\Lambda^{m}_{s}h(\widetilde{\xi}_{s})1_{\{\Lambda^{m}_{s}h(\widetilde{\xi}_{s})\ge 1\}}\\
\vee\limsup_{i\to+\infty}e^{\lambda_{1}\tau_{i}}\Theta_{i}\pi(\widetilde{\xi}_{\tau_{i}-},h)1_{\{\Theta_{i}\pi(\widetilde{\xi}_{\tau_{i}-},h)\ge 1\}}.
\end{multline}
Lemma \ref{lem4} implies that for $m$-a.e. $x\in E$, both $\Lambda^{m}_{s}h(\widetilde{\xi}_{s})1_{\{\Lambda^{m}_{s}h(\widetilde{\xi}_{s})\ge 1\}}$ and $\Theta_{i}\pi(\widetilde{\xi}_{\tau_{i}-},h)1_{\{\Theta_{i}\pi(\widetilde{\xi}_{\tau_{i}-},h)\ge 1\}}$ grow subexponentially. Thus when $\lambda_{1}>0$, the right hand side of \eqref{7.30} goes to infinity. Hence we get \eqref{7.29} for $m$-a.e. $x\in E$.
\qed

\section{Necessary conditions for a
non-degenerate martingale limit}\label{Sec6}

In this section we will give necessary conditions for the fundamental martingale to
have a non-degenerate limit.
Recall from Proposition \ref{prop2} that
$\widetilde{p}(t,x,y)$ is the transition density of the spine $\widetilde{\xi}$ with respect to the measure $\rho$.
We start with the following assumption.

\medskip

\noindent\textbf{Assumption 4.} $$\lim_{t\to+\infty}\sup_{x\in E}\esup_{y\in E}|\widetilde{p}(t,x,y)-1|=0.$$

\medskip

\begin{proposition}
  Suppose that Assumptions 0-4 hold.
Then $\rho$ is an ergodic measure for $(\widetilde{\mathfrak{S}}_{t})_{t\ge 0}$, that is, $\rho$ is an invariant probability measure for $(\widetilde{\mathfrak{S}}_{t})_{t\ge 0}$ and
for any invariant set $B$, either $\rho(B)=0$ or $\rho(B)=1$.
\end{proposition}

\proof Recall from Proposition \ref{prop:invariant measure} that $\rho$ is an invariant probability measure for
$(\widetilde{\mathfrak{S}}_{t})_{t \ge 0}$.
By \cite[Theorem 3.2.4]{Daprato}, it suffices to prove that for any $\varphi\in L^{2}(E,\rho)$,
\begin{equation}
\lim_{t\to+\infty}\frac{1}{t}\int_{0}^{t}\widetilde{\mathfrak{S}}_{s}\varphi {\rm d}s=\langle \varphi,\rho\rangle\quad\mbox{  in }L^{2}(E,\rho).\label{prop5.1}
\end{equation}
It follows from Assumption 4 that for any $\eps>0$, there is $t_{0}>0$ such that
\begin{equation}
\sup_{x\in E}\esup_{y\in E}|\widetilde{p}(s,x,y)-1|\le \eps\quad \mbox{for all }s\ge t_{0}.\label{prop5.2}
\end{equation}
For $x\in E$ and $t>t_{0}$,
\begin{align}
\frac{1}{t}\int_{0}^{t}\widetilde{\mathfrak{S}}_{s}\varphi {\rm d}s-\langle \varphi,\rho\rangle
&=\frac{1}{t}\int_{0}^{t_{0}}\widetilde{\mathfrak{S}}_{s}\varphi {\rm d}s-\frac{t_{0}}{t}\langle \varphi,\rho\rangle\nonumber\\
&\quad+\frac{1}{t}\int_{t_{0}}^{t}{\rm d}s\int_{E}\left(\widetilde{p}(s,x,y)-1\right)\varphi(y)\rho({\rm d}y).\label{prop5.3}
\end{align}
By \eqref{prop5.2} and Jensen's inequality, we have
\begin{align}
&\|\frac{1}{t}\int_{t_{0}}^{t}{\rm d}s\int_{E}\left(\widetilde{p}(s,x,y)-1\right)\varphi(y)\rho({\rm d}y)\|^{2}_{L^{2}(E,\rho)}\nonumber\\
&=\frac{1}{t^{2}}\int_{E}\rho({\rm d}x)\left(\int_{t_{0}}^{t}{\rm d}s\int_{E}\left(\widetilde{p}(s,x,y)-1\right)\varphi(y)\rho({\rm d}y)\right)^{2}\nonumber\\
&\le\frac{t-t_{0}}{t^{2}}\int_{E}\rho({\rm d}x)\int_{t_{0}}^{t}{\rm d}s\int_{E}\left(\widetilde{p}(s,x,y)-1\right)^{2}\varphi(y)^{2}\rho({\rm d}y)\nonumber\\
&\le\frac{(t-t_{0})^{2}}{t^{2}}\,\eps^{2}\|\varphi\|^{2}_{L^{2}(E,\rho)}.\label{prop5.4}
\end{align}
Moreover, by Jensen's inequality and \eqref{i2},
\begin{align}
\|\frac{1}{t}\int_{0}^{t_{0}}\widetilde{\mathfrak{S}}_{s}\varphi {\rm d}s\|^{2}_{L^{2}(E,\rho)}
&=\frac{1}{t^{2}}\int_{E}\rho({\rm d}x)\left(\int_{0}^{t_{0}}\widetilde{\mathfrak{S}}_{s}\varphi(x) {\rm d}s\right)^{2}\nonumber\\
&\le\frac{t_{0}}{t^{2}}\int_{E}\rho({\rm d}x)\int_{0}^{t_{0}}\widetilde{\mathfrak{S}}_{s}(\varphi^{2})(x){\rm d}s\nonumber\\
&=\frac{t^{2}_{0}}{t^{2}}\int_{E}\varphi(x)^{2}\rho({\rm d}x)=\frac{t^{2}_{0}}{t^{2}}\|\varphi\|^{2}_{L^{2}(E,\rho)}.\label{prop5.5}
\end{align}
By \eqref{prop5.3}--\eqref{prop5.5}, we have
\begin{align}
&\|\frac{1}{t}\int_{0}^{t}\widetilde{\mathfrak{S}}_{s}\varphi {\rm d}s-\langle \varphi,\rho\rangle\|_{L^{2}(E,\rho)}\nonumber\\
&\le\frac{t_{0}}{t}\|\varphi\|_{L^{2}(E,\rho)}
+\frac{t_{0}}{t}\left|\langle \varphi,\rho\rangle\right|
+\frac{t-t_{0}}{t}\,\eps\|\varphi\|_{L^{2}(E,\rho)}\nonumber\\
&\le\frac{2t_{0}}{t}\|\varphi\|_{L^{2}(E,\rho)}+\frac{t-t_{0}}{t}\,\eps\|\varphi\|_{L^{2}(E,\rho)}.\nonumber
\end{align}
Letting $t\to+\infty$ and then $\eps\to 0$, we get \eqref{prop5.1}.
\qed

\medskip

Define
\begin{alignat}{2}
E_1&:=
\{x\in E:\ \mathrm{supp}\,\Pi^{L}(x,\cdot)\supseteq [N,+\infty)\mbox{ for some }N\ge 0\},\label{E1} \\
E_2&:=\{x\in A:\ \mathrm{supp}\,\Pi^{NL}(x,\cdot)\supseteq [N,+\infty)\mbox{ for some }N\ge 0\}.\label{E2}
\end{alignat}
The main result of this section is the following theorem.

\begin{theorem}\label{them4}
Suppose that Assumptions 0-4 hold. Then $\p_{\mu}\left(W^{h}_{\infty}(X)=0\right)=1$ for all $\mu\in\me$ if either \eqref{llogl1} fails or the
 following conditions hold.
\center{$\lambda_{1}\ge 0$ and
$m(E_1\cup E_2)>0$.}
\end{theorem}

To prove Theorem \ref{them4}, we need the following lemma.
\begin{lemma}\label{lem5}
Suppose that Assumptions 0-4 hold.
\begin{description}
\item{(i)} If
$m(E_1\cup E_2)>0$,
then for $m$-almost every $x\in E$,
\begin{equation}
\limsup_{D^{m}\ni s\to+\infty}\Lambda^{m}_{s}h(\widetilde{\xi}_{s})\vee \limsup_{i\to+\infty}\Theta_{i}\pi(\widetilde{\xi}_{\tau_{i}-},h)=+\infty \quad\pp_{\cdot,\,x}\mbox{-a.s.};\label{8.7''}
\end{equation}
\item{(ii)}
if \eqref{llogl1} fails,
then for $m$-almost every $x\in E$,
\begin{equation}
\limsup_{D^{m}\ni s\to+\infty}\frac{\log^{+}\Lambda^{m}_{s}h(\widetilde{\xi}_{s})}{s}\vee \limsup_{i\to+\infty}\frac{\log^{+}\Theta_{i}\pi(\widetilde{\xi}_{\tau_{i}-},h)}{\tau_{i}}=+\infty\quad\pp_{\cdot,\,x}\mbox{-a.s.}\label{8.7'}
\end{equation}
\end{description}
\end{lemma}

\proof
It is easy to see that \eqref{8.7'} is equivalent to saying that for $m$-almost every $x\in E$ and all $\lambda<0$,
\begin{equation}
\limsup_{D^{m}\ni s\to+\infty} e^{\lambda s}\Lambda^{m}_{s}h(\widetilde{\xi}_{s})\vee \limsup_{i\to+\infty}e^{\lambda \tau_{i}}\Theta_{i}\pi(\widetilde{\xi}_{\tau_{i}-},h)=+\infty\quad\pp_{\cdot,\,x}\mbox{-a.s.}\nonumber
\end{equation}
We divide the conditions of this lemma into two cases, and prove the results separately.

\noindent Case I: Suppose either one of the following conditions holds:
\begin{description}
\item{(I.a)}
$m(E_1)>0$;
\item{(I.b)}
$(\displaystyle\int_{(0,+\infty)}r h(\cdot)\log^{+}(r h(\cdot))\Pi^{L}(\cdot,{\rm d}r),\widehat{h}) =+\infty.$
\end{description}
Let $\lambda<0$ be an arbitrary constant. To prove \eqref{8.7''} ( resp. \eqref{8.7'}) under condition (I.a) (resp. (I.b)), it suffices to prove that for
$m$-a.e. $x\in E$ and any $M\ge 1$,
\begin{alignat}{2}
\pp_{\cdot,\, x}\left(\sum_{s\in D^{m}}1_{\{\Lambda^{m}_{s}h(\widetilde{\xi}_{s})\ge M\}}=+\infty\right)&=1\nonumber\\
\text{(resp. }\pp_{\cdot,\, x}\left(\sum_{s\in D^{m}}1_{\{e^{\lambda s}\Lambda^{m}_{s}h(\widetilde{\xi}_{s})\ge M\}}=+\infty\right)&=1\text{ )}.
\label{8.10}
\end{alignat}
For $0\le s\le t<+\infty$, $\theta\le 0$ and $M\ge 1$, let
$$\mathrm{I}_{\theta}(s,t):=\int_{s}^{t}{\rm d}r\int_{(0,+\infty)}u1_{\{e^{\theta r}uh(\widetilde{\xi}_{r})\ge M\}}\Pi^{L}(\widetilde{\xi}_{r},{\rm d}u),$$
and $\mathrm{I}_{\theta}(t):=\mathrm{I}_{\theta}(0,t)$. Recall that, given $\widetilde{\xi}$, for any $T>0$, $\#\{s\in D^{m}_{T}:\ e^{\theta s}\Lambda^{m}_{s}h(\widetilde{\xi}_{s})\ge M\}$ is a Poisson random variable with parameter $\mathrm{I}_{\theta}(T)$. Hence \eqref{8.10} would follow if for $m$-a.e. $x\in E$,
\begin{equation}
\pp_{\cdot,\, x}\left(\mathrm{I}_{0}(\infty)=+\infty\right)=1\quad(\mbox{resp. }\pp_{\cdot,\, x}\left(\mathrm{I}_{\lambda}(\infty)=+\infty\right)=1)\label{8.11}
\end{equation}
under condition (I.a) (resp. (I.b)).
Let $\nu({\rm d}x):=\widehat{h}(x)m({\rm d}x)$.
Clearly $\pp_{\cdot,\nu}=\int_{E}\pp_{\cdot,\, x}\,\rho({\rm d}x)$.
Recall that $\rho$ is an invariant measure for $\widetilde{\mathfrak{S}}_{t}$. By Fubini's theorem,
\begin{align}
\pp_{\cdot,\nu}\left(\mathrm{I}_{\theta}(T)\right)
&=\int_{E}\pp_{\cdot,\, x}\left[\int_{0}^{T}{\rm d}r\int_{(0,+\infty)}u1_{\{e^{\theta r}uh(\widetilde{\xi}_{r})\ge M\}}\Pi^{L}(\widetilde{\xi}_{r},{\rm d}u)\right]\rho({\rm d}x)\nonumber\\
&=\int_{0}^{T}{\rm d}r\int_{E}\pp_{\cdot,\, x}\left[\int_{(0,+\infty)}u1_{\{e^{\theta r}uh(\widetilde{\xi}_{r})\ge M\}}\Pi^{L}(\widetilde{\xi}_{r},{\rm d}u)\right]\rho({\rm d}x)\nonumber\\
&=\int_{0}^{T}{\rm d}r\int_{E}\rho({\rm d}x)\int_{(0,+\infty)}u1_{\{e^{\theta r}uh(x)\ge M\}}\Pi^{L}(x,{\rm d}u).\label{I.2}
\end{align}
By the boundedness of $h$ and $x\mapsto
\int_{(0,+\infty)}(u\wedge u^{2}) \Pi^{L}(x,{\rm d}u)
$, we have
\begin{align}
&\pp_{\cdot,\nu}\left(\mathrm{I}_{\theta}(T)\right)\nonumber\\
&\le T\int_{E}\widehat{h}(x)m({\rm d}x)\int_{u\ge M/h(x)}h(x)u\Pi^{L}(x,{\rm d}u)\nonumber\\
&=T\int_{E}\widehat{h}(x)m({\rm d}x)\int_{u\ge M/h(x)}h(x)(1\vee \frac{1}{u})(u\wedge u^{2})\Pi^{L}(x,{\rm d}u)\nonumber\\
&\le T\left(1\vee \frac{\|h\|_{\infty}}{M}\right)\|\int_{(0,+\infty)}(u\wedge u^{2})\Pi^{L}(x,{\rm d}u)\|_{\infty}
\int_{E}h(x)\widehat{h}(x)m({\rm d}x)<+\infty.\nonumber
\end{align}
Thus $\pp_{\cdot,\nu}(\mathrm{I}_{\theta}(T)<+\infty)=1$.
On the other hand, by the Markov property of $\widetilde{\xi}$ and \eqref{I.2},
\begin{align}
&\pp_{\cdot,\nu}\left(\mathrm{I}_{\theta}(T)^{2}\right)\nonumber\\
&=\int_{E}\pp_{\cdot,\, x}\left(\mathrm{I}_{\theta}(T)^{2}\right)\rho({\rm d}x)\nonumber\\
&=2\int_{E}\rho({\rm d}x)\pp_{\cdot,\,x}\left[\int_{0}^{T}{\rm d}s\int_{(0,+\infty)}u1_{\{e^{\theta s}uh(\widetilde{\xi}_{s})\ge M\}}\Pi^{L}(\widetilde{\xi}_{s},{\rm d}u)
\right.\nonumber\\
&\qquad\left.\cdot\int_{s}^{T}{\rm d}r\int_{(0,+\infty)}v1_{\{e^{\theta r}vh(\widetilde{\xi}_{r})\ge M\}}\Pi^{L}(\widetilde{\xi}_{r},{\rm d}v)\right]\nonumber\\
&=2\int_{E}\rho({\rm d}x)\pp_{\cdot,\, x}\left[\int_{0}^{T}{\rm d}s\int_{(0,+\infty)}u1_{\{e^{\theta s}uh(\widetilde{\xi}_{s})\ge M\}}\Pi^{L}(\widetilde{\xi}_{s},{\rm d}u)\right.\nonumber\\
&\quad\quad\cdot
\left.\pp_{\cdot,\widetilde{\xi}_{s}}\left(\int_{0}^{T-s}{\rm d}r\int_{(0,+\infty)}
v1_{\{e^{\theta (r+s)}vh(\widetilde{\xi}_{r})\ge M\}}\Pi^{L}(\widetilde{\xi}_{r},{\rm d}v)\right)\right]\nonumber\\
&=2\int_{E}\rho({\rm d}x)\int_{0}^{T}{\rm d}s\int_{(0,+\infty)}u1_{\{e^{\theta s}uh(x)\ge M\}}\Pi^{L}(x,{\rm d}u)
\nonumber\\
&\qquad\cdot\pp_{\cdot,\, x}\left(\int_{0}^{T-s}{\rm d}r\int_{(0,+\infty)}
v1_{\{e^{\theta (r+s)}vh(\widetilde{\xi}_{r})\ge M\}}\Pi^{L}(\widetilde{\xi}_{r},{\rm d}v)\right)\nonumber\\
&\le 2\int_{E}\rho({\rm d}x)\int_{0}^{T}{\rm d}s\int_{(0,+\infty)}u1_{\{e^{\theta s}uh(x)\ge M\}}\Pi^{L}(x,{\rm d}u)
\nonumber\\
&\quad\quad\cdot\pp_{\cdot,\, x}\left(\int_{0}^{T}{\rm d}r\int_{(0,+\infty)}
v1_{\{e^{\theta r}vh(\widetilde{\xi}_{r})\ge M\}}\Pi^{L}(\widetilde{\xi}_{r},{\rm d}v)\right)\nonumber\\
&=2\int_{E}\rho({\rm d}x)\int_{0}^{T}{\rm d}s\int_{(0,+\infty)}u1_{\{e^{\theta s}uh(x)\ge M\}}\Pi^{L}(x,{\rm d}u)\pp_{\cdot,\,x}\left(\mathrm{I}_{\theta}(T)\right).\label{8.14}
\end{align}
Assumption 4 implies that there are constants $t_{1},\delta>0$ such that
\begin{equation}
\sup_{x\in E}\esup_{y\in E}\widetilde{p}(t,x,y)\le 1+\delta\quad\mbox{ for all }t\ge t_{1}.\label{II.0}
\end{equation}
Using Fubini's theorem, \eqref{II.0} and \eqref{I.2}, we have for $T>t_{1}$,
\begin{align}
\pp_{\cdot,\, x}\left[\mathrm{I}_{\theta}(t_{1},T)\right]
&=\int_{t_{1}}^{T}{\rm d}r\int_{E}\widetilde{p}(r,x,y)\rho({\rm d}y)\int_{(0,+\infty)}v1_{\{e^{\theta r}vh(y)\ge M\}}\Pi^{L}(y,{\rm d}v)\nonumber\\
&\le(1+\delta)\int_{t_{1}}^{T}{\rm d}r\int_{E}\rho({\rm d}y)\int_{(0,+\infty)}v1_{\{e^{\theta r}vh(y)\ge M\}}\Pi^{L}(y,{\rm d}v)\nonumber\\
&\le(1+\delta)\pp_{\cdot,\nu}\left(\mathrm{I}_{\theta}(T)\right).\label{8.17}
\end{align}
On the other hand, for $x\in E$,
\begin{align}
&\pp_{\cdot,\, x}\left(\mathrm{I}_{\theta}(t_{1})\right)\nonumber\\
&=\pp_{\cdot,\, x}\left[\int_{0}^{t_{1}}{\rm d}r\int_{(0,+\infty)}v1_{\{e^{\theta r}vh(\widetilde{\xi}_{r})\ge M\}}\Pi^{L}(\widetilde{\xi}_{r},{\rm d}v)\right]\nonumber\\
&=\int_{0}^{t_{1}}{\rm d}r\int_{E}\widetilde{p}(r,x,y)\rho({\rm d}y)\int_{(0,+\infty)}v1_{\{e^{\theta r}vh(y)\ge M\}}\Pi^{L}(y,{\rm d}v)\nonumber\\
&\le\int_{0}^{t_{1}}{\rm d}r\int_{E}\widetilde{p}(r,x,y)\rho({\rm d}y)\int_{v\ge M/h(y)}\left(1\vee \frac{1}{v}\right)\left(v\wedge v^{2}\right)\Pi^{L}(y,{\rm d}v)\nonumber
\end{align}
\begin{align}
&\le t_{1}\left(1\vee\frac{\|h\|_{\infty}}{M}\right)\|\int_{(0,+\infty)}(v\wedge v^{2})\Pi^{L}(\cdot,{\rm d}v)\|_{\infty}=:c_{1}<+\infty.\label{8.18}
\end{align}
It follows from \eqref{8.17} and \eqref{8.18} that for $T>t_{1}$,
$$\pp_{\cdot,\,x}\left(\mathrm{I}_{\theta}(T)\right)=\pp_{\cdot,\,x}\left(\mathrm{I}_{\theta}(t_{1})\right)+\pp_{\cdot,\,x}\left(\mathrm{I}_{\theta}(t_{1},T)\right)\le c_{1}+(1+\delta)\pp_{\cdot,\nu}\left(\mathrm{I}_{\theta}(T)\right).$$
This together with \eqref{I.2} and \eqref{8.14} implies that
\begin{equation}
\pp_{\cdot,\nu}\left(\mathrm{I}_{\theta}(T)^{2}\right)\le 2c_{1}\pp_{\cdot,\nu}\left(\mathrm{I}_{\theta}(T)\right)+2(1+\delta)\pp_{\cdot,\nu}\left(\mathrm{I}_{\theta}(T)\right)^{2}.\nonumber
\end{equation}
Hence by
the Cauchy-Schwarz inequality,
we have
\begin{equation}
\pp_{\cdot,\nu}\left(\mathrm{I}_{\theta}(T)\ge
\frac{1}{2}\pp_{\cdot,\nu}(\mathrm{I}_{\theta}(T))\right)\ge
\frac{\pp_{\cdot,\nu}(\mathrm{I}_{\theta}(T))^{2}}{4\pp_{\cdot,\nu}(\mathrm{I}_{\theta}(T)^{2})}\ge \frac{\pp_{\cdot,\nu}(I_{\theta}(T))}{8c_{1}+8(1+\delta)\pp_{\cdot,\nu}(I_{\theta}(T))}.\label{8.19}
\end{equation}
Note that $\pp_{\cdot,\nu}(\mathrm{I}_{0}(T))=T\int_{E}\rho({\rm d}x)\int_{u\ge h(x)/M}u\Pi^{L}(x,{\rm d}u)$.
Condition (I.a) implies that the integral on the right hand side is positive. Hence $\pp_{\cdot,\nu}(\mathrm{I}_{0}(T))\to +\infty$ as $T\to +\infty$.
On the other hand, note that by \eqref{I.2} and Fubini's theorem, for $\lambda<0$,
\begin{align}
&\pp_{\cdot,\nu}\left(\mathrm{I}_{\lambda}(T)\right)\nonumber\\
&\le
\int_{E}\widehat{h}(x)m({\rm d}x)\int_{(0,+\infty)}h(x)u\Pi^{L}(x,{\rm d}u)\int_{0}^{T}1_{\{s\le \frac{\log^{+}(h(x)u)-\log M}{-\lambda }\}}{\rm d}s\nonumber\\
&=\int_{E}\widehat{h}(x)m({\rm d}x)\int_{(0,+\infty)}h(x)u\left(T\wedge \frac{\log^{+}(h(x)u)-\log M}{-\lambda }\right)^{+}\Pi^{L}(x,{\rm d}u).\nonumber
\end{align}
Clearly condition (I.b) implies that $\lim_{T\to+\infty}\pp_{\cdot,\nu}(\mathrm{I}_{\lambda}(T))=+\infty$. Thus by letting $T\to +\infty$ in \eqref{8.19}, we get
$\pp_{\cdot,\nu}(\mathrm{I}_{0}(\infty)=+\infty)>0$ (resp. $\pp_{\cdot,\nu}(\mathrm{I}_{\lambda}(\infty)=+\infty)>0$) under condition (I.a) (resp. (I.b)).
Since $\{\mathrm{I}_{0}(\infty)=+\infty\}$ (resp. $\{\mathrm{I}_{\lambda}(\infty)=+\infty\}$) is an invariant event of
the canonical dynamic system associated with
$(\widetilde{\mathfrak{S}}_{t})_{t\ge 0}$
and ergodic measure $\rho$, it follows from
\cite[Theorem 1.2.4]{Daprato} that
$\pp_{\cdot,\nu}(\mathrm{I}_{0}(\infty)=+\infty)=1$ (resp. $\pp_{\cdot,\nu}(\mathrm{I}_{\lambda}(\infty)=+\infty)=1$) under condition (I.a) (resp. (I.b)). Hence we prove \eqref{8.11}.

\medskip

\noindent Case II. Suppose either one of the following conditions holds:
\begin{description}
\item{(II.a)}
$m(E_2)>0$;
\item{(II.b)} $\displaystyle(
\int_{(0,+\infty)}\pi(\cdot,h)r\log^{+}(\pi(\cdot,h)r)\Pi^{NL}(\cdot,{\rm d}r),\widehat{h}) =+\infty.$
\end{description}
Let $\lambda<0$ be an arbitrary constant. To prove \eqref{8.7''} ( resp. \eqref{8.7'}) under condition (II.a) (resp. (II.b)), it suffices to prove that for
$m$-a.e. $x\in E$ and any $M\ge 1$,
\begin{alignat}{2}
\pp_{\cdot,\,x}\left(\sum_{i=1}^{+\infty}1_{\{\Theta_{i}\pi(\widetilde{\xi}_{\tau_{i}-},h)\ge M\}}=+\infty\right)&=1\nonumber\\
\text{(resp. }\pp_{\cdot,\,x}\left(\sum_{i=1}^{+\infty}1_{\{e^{\lambda\tau_{i}}\Theta_{i}\pi(\widetilde{\xi}_{\tau_{i}-},h)\ge M\}}=+\infty\right)&=1\text{ ).}
\label{8.21}
\end{alignat}
The main idea of this proof is similar to that of Case I.
For all $T>0$, $\theta\le 0$ and $M\ge 1$, let $\mathrm{II}_{\theta}(T):=\sum_{\tau_{i}\le T}1_{\{e^{\theta \tau_{i}}\Theta_{i}\pi(\widetilde{\xi}_{\tau_{i}-},h)\ge M\}}$.
For any $s\ge 0$ and $x\in E$, define
\begin{alignat}{2}
f_{\theta}(s,x)&:=\int_{[0,+\infty)}1_{\{e^{\theta s}u\pi(x,h)\ge M\}}\eta(x,{\rm d}u)\nonumber\\
&=\frac{1}{\gamma(x)}1_{A}(x)\int_{(0,+\infty)}1_{\{e^{\theta s}u\pi(x,h)\ge M\}}u\Pi^{NL}(x,{\rm d}u),\nonumber\\
g_{\theta}(s,x)&:=q(x)f_{\theta}(s,x)\nonumber\\
&=\frac{\pi(x,h)}{h(x)}1_{A}(x)\int_{(0,+\infty)}1_{\{e^{\theta s}u\pi(x,h)\ge M\}}u\Pi^{NL}(x,{\rm d}u).\nonumber
\end{alignat}
Recall from Definition \ref{sd1} that, given $\widetilde{\xi}$ (including $\{\tau_{i}:i\ge 1\}$), $\Theta_{i}$ is distributed according to $\eta(\widetilde{\xi}_{\tau_{i}-},{\rm d}r)$. By \eqref{l1}, we have for $x\in E$,
\begin{align}
\pp_{\cdot,\,x}\left(\mathrm{II}_{\theta}(T)\right)
&=\pp_{\cdot,\,x}\left[\sum_{\tau_{i}\le T}\int_{[0,+\infty)}1_{\{e^{\theta \tau_{i}}r\pi(\widetilde{\xi}_{\tau_{i}-},h)\ge M\}}\eta(\widetilde{\xi}_{\tau_{i}-},{\rm d}r)\right]\nonumber\\
&=\pp_{\cdot,\,x}\left[\sum_{\tau_{i}\le T}f_{\theta}(\tau_{i},\widetilde{\xi}_{\tau_{i}-})\right]
=\pp_{\cdot,\,x}\left[\int_{0}^{T}q(\widetilde{\xi}_{s})f_{\theta}(s,\widetilde{\xi}_{s}){\rm d}s\right]\nonumber\\
&=\pp_{\cdot,\,x}\left[\int_{0}^{T}g_{\theta}(s,\widetilde{\xi}_{s}){\rm d}s\right].\label{II.6}
\end{align}
We still use $\nu$ to denote the measure $\widehat{h}(x)m({\rm d}x)$. Since $\rho$ is an invariant measure for $\widetilde{\mathfrak{S}}_{t}$, by Fubini's theorem,
\begin{align}
\pp_{\cdot,\nu}\left(\mathrm{II}_{\theta}(T)\right)
&=\int_{E}\pp_{\cdot,\,x}\left(\mathrm{II}_{\theta}(T)\right)\rho({\rm d}x)\nonumber\\
&=\int_{0}^{T}{\rm d}s\int_{E}\pp_{\cdot,\,x}\left(g_{\theta}(s,\widetilde{\xi}_{s})\right)\rho({\rm d}x)\nonumber\\
&=\int_{0}^{T}{\rm d}s\int_{E}g_{\theta}(s,x)\rho({\rm d}x)\label{8.22}\\
&=\int_{A}\widehat{h}(x)m({\rm d}x)\int_{(0,+\infty)}\pi(x,h)r\Pi^{NL}(x,{\rm d}r)\int_{0}^{T}1_{\{e^{\theta s}r\pi(x,h)\ge M\}}{\rm d}s.\label{II.1}
\end{align}
It then follows by Assumption 3.(ii) that
$$\mbox{RHS of }\eqref{II.1}\le T\|\int_{(0,+\infty)}y\Pi^{NL}(\cdot,{\rm d}y)\|_{\infty}\int_{A}\pi(x,h)\widehat{h}(x)m({\rm d}x)<+\infty.$$
Therefore $\pp_{\cdot,\nu}\left(\mathrm{II}_{\theta}(T)<+\infty\right)=1$.
Recall that given $\widetilde{\xi}$ (including $\{\tau_{i}:i\ge 1\}$),
$\{\Theta_{i}:i\ge 1\}$
are mutually independent, we have
\begin{align}
&\pp_{\cdot,\,x}\left(\mathrm{II}_{\theta}(T)^{2}\right)-\pp_{\cdot,\,x}\left(\mathrm{II}_{\theta}(T)\right)\nonumber\\
&=\pp_{\cdot,\,x}\left[\sum_{\tau_{i},\tau_{j}\le T,i\not=j}
1_{\{e^{\theta \tau_{i}}\Theta_{i}\pi(\widetilde{\xi}_{\tau_{i}-},h)\ge M\}}1_{\{e^{\theta \tau_{j}}\Theta_{j}\pi(\widetilde{\xi}_{\tau_{j}-},h)\ge M\}}
\right]\nonumber\\
&=\pp_{\cdot,\,x}\left[\sum_{\tau_{i},\tau_{j}\le T,i\not=j}\int_{[0,+\infty)}1_{\{e^{\theta \tau_{i}}y\pi(\widetilde{\xi}_{\tau_{i}-},h)\ge M\}}\eta(\widetilde{\xi}_{\tau_{i}-},{\rm d}y)
\int_{[0,+\infty)}1_{\{e^{\theta \tau_{j}}z\pi(\widetilde{\xi}_{\tau_{j}-},h)\ge M\}}\eta(\widetilde{\xi}_{\tau_{j}-},{\rm d}z)\right]\nonumber\\
&=\pp_{\cdot,\,x}\left[\sum_{\tau_{i},\tau_{j}\le T,i\not=j}f_{\theta}(\tau_{i},\widetilde{\xi}_{\tau_{i}-})f_{\theta}(\tau_{j},\widetilde{\xi}_{\tau_{j}-})\right].\nonumber
\end{align}
Thus by \eqref{l2},
\begin{align}
&\pp_{\cdot,\,x}\left(\mathrm{II}_{\theta}(T)^{2}\right)-\pp_{\cdot,\,x}\left(\mathrm{II}_{\theta}(T)\right)\nonumber\\
&=2\pp_{\cdot,\,x}\left[\int_{0}^{T}q(\widetilde{\xi}_{s})f_{\theta}(s,\widetilde{\xi}_{s}){\rm d}s\int_{E}
\widetilde{\Pi}_{y}\left(\int_{0}^{T-s}q(\widetilde{\xi}_{r})
f_{\theta}
(s+r,\widetilde{\xi}_{r}){\rm d}r\right)\pi^{h}(\widetilde{\xi}_{s},{\rm d}y)\right]\nonumber\\
&=2\pp_{\cdot,\,x}\left[\int_{0}^{T}g_{\theta}(s,\widetilde{\xi}_{s}){\rm d}s\int_{E}
\widetilde{\Pi}_{y}\left(\int_{0}^{T-s}g_{\theta}(s+r,\widetilde{\xi}_{r}){\rm d}r\right)\pi^{h}(\widetilde{\xi}_{s},{\rm d}y)\right].\label{II.2}
\end{align}
Note that for all $x\in E$ and $\theta\le 0$, $s\mapsto g_{\theta}(s,x)$ is non-increasing. Thus it follows from \eqref{II.2} that
\begin{align}
\pp_{\cdot,\,x}\left(\mathrm{II}_{\theta}(T)^{2}\right)&\le 2\pp_{\cdot,\,x}\left[\int_{0}^{T}g_{\theta}(s,\widetilde{\xi}_{s}){\rm d}s\int_{E}
\widetilde{\Pi}_{y}\left(\int_{0}^{T}g_{\theta}(r,\widetilde{\xi}_{r}){\rm d}r\right)\pi^{h}(\widetilde{\xi}_{s},{\rm d}y)\right]\nonumber\\
&\quad+\pp_{\cdot,\,x}\left(\mathrm{II}_{\theta}(T)\right).\label{II.3}
\end{align}
By Fubini's theorem, \eqref{II.0} and
\eqref{8.22},
 we have for $y\in E$ and $T>t_{1}$,
\begin{align}
\widetilde{\Pi}_{y}\left(\int_{t_{1}}^{T}g_{\theta}(s,\widetilde{\xi}_{s}){\rm d}s\right)
&=\int_{t_{1}}^{T}{\rm d}s\int_{E}\widetilde{p}(s,y,z)g_{\theta}(s,z)\rho({\rm d}z)\nonumber\\
&\le(1+\delta)\int_{t_{1}}^{T}{\rm d}s\int_{E}g_{\theta}(s,z)\rho({\rm d}z)\nonumber\\
&\le(1+\delta)\pp_{\cdot,\nu}\left(\mathrm{II}_{\theta}(T)\right).\label{II.4}
\end{align}
On the other hand, by Assumption 3.(iii),
\begin{align}
&\sup_{y\in E}\widetilde{\Pi}_{y}\left(\int_{0}^{t_{1}}g_{\theta}(s,\widetilde{\xi}_{s}){\rm d}s\right)\nonumber\\
&=\sup_{y\in E}\int_{0}^{t_{1}}{\rm d}s\int_{A}\widetilde{p}(s,y,z)\frac{\pi(z,h)}{h(z)}\rho({\rm d}z)\int_{(0,+\infty)}r
1_{\{e^{\theta s}r\pi(z,h)\ge M\}}\Pi^{NL}(z,{\rm d}r)\nonumber\\
&\le\|\frac{\pi(\cdot ,h)}{h}1_{A}\|_{\infty}\|\int_{(0,+\infty)}r\Pi^{NL}(\cdot,{\rm d}r)\|_{\infty}\sup_{y\in E}\int_{0}^{t_{1}}{\rm d}s\int_{A}\widetilde{p}(s,y,z)\rho({\rm d}z)\nonumber\\
&\le c_{2}t_{1}=:c_{3}<+\infty.\nonumber
\end{align}
This and \eqref{II.4} imply that
$$\widetilde{\Pi}_{y}\left(\int_{0}^{T}g_{\theta}(s,\widetilde{\xi}_{s}){\rm d}s\right)\le c_{3}+(1+\delta)\pp_{\cdot,\nu}(\mathrm{II}_{\theta}(T))\quad\mbox{ for all $y\in E$ and $T>t_{1}$.}$$
This together with \eqref{II.3} and \eqref{II.6} implies that
$$\pp_{\cdot,\,x}\left(\mathrm{II}_{\theta}(T)^{2}\right)\le (1+2c_{3})\pp_{\cdot,\,x}\left(\mathrm{II}_{\theta}(T)\right)+2(1+\delta)\pp_{\cdot,\nu}\left(\mathrm{II}_{\theta}(T)\right)\pp_{\cdot,\,x}\left(\mathrm{II}_{\theta}(T)\right).$$
Consequently,
\begin{align}
\pp_{\cdot,\nu}\left(\mathrm{II}_{\theta}(T)^{2}\right)&=\int_{E}\pp_{\cdot,\,x}\left(\mathrm{II}_{\theta}(T)^{2}\right)\rho({\rm d}x)\nonumber\\
&\le (1+2c_{3})\pp_{\cdot,\nu}\left(\mathrm{II}_{\theta}(T)\right)+2(1+\delta)\pp_{\cdot,\nu}\left(\mathrm{II}_{\theta}(T)\right)^{2}.\nonumber
\end{align}
Recall that $\pp_{\cdot,\nu}(\mathrm{II}_{0}(T))=T\int_{A}\pi(x,h)\widehat{h}(x)m({\rm d}x)\int_{r\ge M/\pi(x,h)}r\Pi^{NL}(x,{\rm d}r)$. Condition (II.a) implies that the integral on the right hand side is positive. Thus $\pp_{\cdot,\nu}(\mathrm{II}_{0}(T))\to+\infty$ as $T\to+\infty$.
On the other hand, note that by \eqref{II.1} and Fubini's theorem, for $\lambda<0$,
\begin{multline}
\pp_{\cdot,\nu}\left(\mathrm{II}_{\lambda}(T)\right)
=\int_{A}\widehat{h}(x)m({\rm d}x)\int_{(0,+\infty)}\Pi^{NL}(x,{\rm d}r)\pi(x,h)r\left(\frac{\log^{+}(\pi(x,h)r)-\log M}{-\lambda }\wedge T\right)^{+}.\nonumber
\end{multline}
Clearly condition (II.b) implies that $\lim_{T\to+\infty}\pp_{\cdot,\nu}(\mathrm{II}_{\lambda}(T))=+\infty$.  Similarly by using
the Cauchy-Schwarz inequality
and letting $T\to +\infty$, we get
$\pp_{\cdot,\nu}(\mathrm{II}_{0}(\infty)=+\infty)>0$ (resp. $\pp_{\cdot,\nu}(\mathrm{II}_{\lambda}(\infty)=+\infty)>0$) under condition (II.a) (resp. (II.b)).

For each $n\ge 1$, we denote by $\widetilde{\mathcal{G}}_{n}$ the $\sigma$-field generated by $\widetilde{\xi}$ up to time $\tau_{n}$ (including $\{\tau_1, \cdots, \tau_n\}$) and $\{\Theta_{i}: i\le n\}$. Obviously for each $i\ge 1$, both
$$
1_{\{e^{\theta \tau_{i}}\Theta_{i}\pi(\widetilde{\xi}_{\tau_{i}-},h)\ge M\}}
\quad \mbox{ and } \quad
\int_{[0,+\infty)}1_{\{e^{\theta \tau_{i}}r\pi(\widetilde{\xi}_{\tau_{i}-},h)\ge M\}}\eta(\widetilde{\xi}_{\tau_{i}-},{\rm d}r)
$$
are $\widetilde{\mathcal{G}}_{i}$-measurable. Moreover, for every $x\in E$, under $\pp_{\cdot,\,x}$,
\begin{align}
&\pp_{\cdot,\,x}\left(1_{\{e^{\theta \tau_{i+1}}\Theta_{i+1}\pi(\widetilde{\xi}_{\tau_{i+1}-},h)\ge M\}}\,|\,\widetilde{\mathcal{G}}_{i}\right)\notag\\
&=\pp_{\cdot,\,x}\left(\int_{[0,+\infty)}
1_{\{e^{\theta \tau_{i+1}}r\pi(\widetilde{\xi}_{\tau_{i+1}-},h)\ge M\}}\eta(\widetilde{\xi}_{\tau_{i+1}-},{\rm d}r)\,|\,\widetilde{\mathcal{G}}_{i}\right).\notag
\end{align}
Applying the second Borel-Cantelli lemma (see, for example, \cite[Corollary 5.3.2]{Durrett}) to both sides of the above equality, we get that
$$
\left\{\sum_{i=1}^{+\infty}1_{\{e^{\theta \tau_{i}}\Theta_{i}\pi(\widetilde{\xi}_{\tau_{i}-},h)\ge M\}}=+\infty\right\}
=\left\{\sum_{i=1}^{+\infty}\int_{[0,+\infty)}1_{\{e^{\theta \tau_{i}}r\pi(\widetilde{\xi}_{\tau_{i}-},h)\ge M\}}\eta(\widetilde{\xi}_{\tau_{i}-},{\rm d}r)=+\infty\right\}$$
under $\pp_{\cdot,\,x}$.
It is easy to see from the above representation that $\{\mathrm{II}_{0}(\infty)=+\infty\}$ (resp. $\{\mathrm{II}_{\lambda}(\infty)=+\infty\}$) is an invariant event of the canonical dynamic system associated with
$(\widetilde{\mathfrak{S}}_{t})_{t\ge 0}$
and ergodic measure $\rho$, so it follows from
\cite[Theorem 1.2.4]{Daprato} that
$\pp_{\cdot,\nu}(\mathrm{II}_{0}(\infty)=+\infty)=1$ (resp. $\pp_{\cdot,\nu}(\mathrm{II}_{\lambda}(\infty)=+\infty)=1$) under condition (II.a) (resp. (II.b)). Thus \eqref{8.21} is valid.
\qed

\medskip

\noindent\textbf{Proof of Theorem \ref{them4}.}
Applying the same argument as in the beginning of the proof of Theorem \ref{them3}.(ii) here, we only need to show that under
the assumptions of Theorem \ref{them4},
\begin{equation*}
\pp_{\delta_{x}}(\limsup_{t\to+\infty}W^{h}_{t}(\Gamma)=+\infty)=1\quad\mbox{ for $m$-a.e. $x\in E$.}
\end{equation*}
In view of \eqref{7.30}, this would follow if for $m$-a.e. $x\in E$,
\begin{equation}
\limsup_{D^{m}\ni s\to+\infty}e^{\lambda_{1}s}\Lambda^{m}_{s}h(\widetilde{\xi}_{s})\vee \limsup_{i\to+\infty}e^{\lambda_{1}\tau_{i}}\Theta_{i}\pi(\widetilde{\xi}_{\tau_{i}-},h)=+\infty\quad\pp_{\cdot,\,x}\mbox{-a.s.}\nonumber
\end{equation}
which, under the assumptions of this theorem, is automatically true by Lemma \ref{lem5}. Hence we complete the proof.\qed

\medskip

The following corollaries follow directly from Theorem \ref{them3} and Theorem \ref{them4}.
\begin{corollary}\label{cor8.3}
Suppose that Assumptions 0-4 hold and
that $m(E_1\cup E_2)>0$ with $E_1$ and $E_2$ defined in
\eqref{E1} and \eqref{E2} respectively.
For every $\mu\in\me^0$,
$W^{h}_{\infty}(X)$ is non-degenerate
if and only if
$\lambda_{1}<0$ and
condition \eqref{llogl1} holds.
Moreover, $X_{t}$ under $\p_{\mu}$ exhibits weak local extinction if $\lambda_{1}\ge 0$.
\end{corollary}

\begin{corollary}\label{cor2}
Suppose Assumptions 0-4 hold and $\lambda_{1}<0$.  For every $\mu\in\me^0$,
$W^{h}_{\infty}(X)$ is non-degenerate
if and only if
condition \eqref{llogl1} holds.
\end{corollary}

\begin{remark}\rm
Note that in the case of purely local branching mechanism, Assumption 4 can be written as
\begin{equation}
\lim_{t\to +\infty}\sup_{x\in E}\esup_{y\in E}\left|p^{h}(t,x,y)-1\right|=0,\nonumber
\end{equation}
where $p^{h}(t,x,y)$ denotes the transition density function of $\xi^{h}$ with respect to the measure $\rho$.
If $E$ is a bounded domain in $\mathbb{R}^{d}$, $m$ is the Lebesgue measure on $\mathbb{R}^{d}$ and $\xi$ is a symmetric diffusion on $E$, then $a(x)\in \mathcal{B}_{b}(E)\subset \mathbf{K}(\xi)\cap L^{2}(E,m)$. Hence for the class of superdiffusions with local branching mechanisms considered in \cite{LRS},
our Assumptions 0-4 hold naturally and
Corollary \ref{cor2}
generalizes \cite[Theorem 1.1]{LRS}.
\end{remark}

\section{Examples}\label{s:examples}

In this section, we will give examples satisfying
Assumptions 0-4.
We will not try to give
the most general examples possible.

\begin{example}\label{e:KP}{\rm
Suppose $E=\{1,2,\cdots,K\}\ (K\ge 2)$, $m$ is the counting measure on $E$ and
$\mathfrak{S}_{t}f(i)= f(i)$ for all $i\in E$, $t\ge 0$ and $f\in\mathcal{B}^{+}(E)$ (that is, there is no spatial motion).
Suppose
$$\phi^{L}(i,\lambda):=a(i)\lambda+b(i)\lambda^{2}+\int_{(0,+\infty)}\left(e^{-\lambda r}-1+\lambda r\right)\Pi^{L}(i,{\rm d}r),$$
$$\phi^{NL}(i,f):=-c(i)\pi(i,f)-\int_{(0,+\infty)}\left(1-e^{-r\pi(i,f)}\right)\Pi^{NL}(i,{\rm d}r),$$
where for each $i\in E$, $a(i)\in (-\infty,+\infty)$, $b(i),c(i)\ge 0$, $(r\wedge r^{2})\Pi^{L}(i,{\rm d}r)$ and $r\Pi^{NL}(i,{\rm d}r)$ are bounded kernels from $E$ to $(0,+\infty)$ with $\{i\in E: \int_{(0,+\infty)}r\Pi^{NL}(i,{\rm d}r)>0\}\not=\emptyset$, and $\pi(i,dj)$ is a probability kernel on $E$ with $\pi(i,\{i\})=0$ for every $i\in E$. As a special case of the model given in Section \ref{sec1}, we have a non-local branching superprocess $\{X_{t}:t\ge 0\}$ in $\me$ with transition probabilities given by
\begin{equation}\nonumber
\p_{\mu}\left[\exp\left(-\langle f,X_{t}\rangle\right)\right]=\exp\left(-\langle V_{t}f,\mu\rangle\right)\quad\mbox{ for }\mu\in\me,\ t\ge 0\mbox{ and }f\in \mathcal{B}^{+}_{b}(E),
\end{equation}
where $V_{t}f(i)$ is the unique non-negative locally bounded solution to the following integral equation:
\begin{equation}\nonumber
V_{t}f(i)=f(i)-\int_{0}^{t}\left(\phi^{L}(i,V_{s}f(i))+\phi^{NL}(i,V_{s}f)\right){\rm d}s\quad\mbox{ for } t\ge 0,\ i\in E.
\end{equation}
For every $i\in E$ and $\mu\in\me$, we define $\mu^{(i)}:=\mu(\{i\})$. The map $\mu\mapsto (\mu^{(1)},\cdots,\mu^{(K)})^{\mathrm{T}}$ is clearly a homeomorphism between $\me$ and the $K$-dimensional product space $[0,+\infty)^{K}$. Hence $\{(X^{(1)}_{t},\cdots,X^{(K)}_{t})^{\mathrm{T}}:t\ge 0\}$ is a Markov process in $[0,+\infty)^{K}$, which is called a $K$-type continuous-state branching process. (Clearly the $1$-type continuous-state branching process defined in a similar way coincides with the classical one-dimensional continuous-state branching process, see, for example, \cite[Chapter 3]{Li}.)
For simplicity, we assume $b(i)\equiv 0$.
For $i,j\in E$, let $p_{ij}:=\pi(i,\{j\})$ and $\gamma(i):=c(i)+\int_{(0,+\infty)}r\Pi^{NL}(i,{\rm d}r)$.
Define the $K\times K$ matrix $
\boldsymbol{M(t)}=(M(t)_{ij})_{ij}$ by $M(t)_{ij}:=\p_{\delta_{i}}\left[X^{(j)}_{t}\right]$ for $i,j\in E$.
Let $\mathfrak{P}_{t}$
denote the mean semigroup of $X$, that is
$$
\mathfrak{P}_{t}f(i)
:=\p_{\delta_{i}}\left[\langle f,X_{t}\rangle\right]=\sum_{j=1}^{K}M(t)_{ij}f(j)\quad\mbox{ for }i\in E,\ t\ge 0\mbox{ and }f\in\mathcal{B}^{+}(E).$$
By the Markov property and \eqref{li2}, $\boldsymbol{M(t)}$ satisfies that
\begin{equation}\nonumber
\boldsymbol{M(0)}=\boldsymbol{I},\quad \boldsymbol{M(t+s)}=\boldsymbol{M(t)}\boldsymbol{M(s)}\quad\mbox{ for }t,s\ge 0,
\end{equation}
\begin{equation}\nonumber
\mbox{ and }\quad M(t)_{ij}=\delta_{j}(i)-a(i)\int_{0}^{t}M(s)_{ij}{\rm d}s+\gamma(i)\sum_{k=1}^{K}p_{ik}\int_{0}^{t}M(s)_{kj}{\rm d}s
\end{equation}
for $i,j\in E$.
This implies that $\boldsymbol{M(t)}$ has a formal matrix generator $\boldsymbol{A}:=(A_{ij})_{ij}$ given by
$$\boldsymbol{M(t)}=e^{\boldsymbol{A}t}=\sum_{n=0}^{+\infty}t^{n}\frac{\boldsymbol{A}^{n}}{n!}, \quad \mbox{ and }A_{ij}=\gamma(i)p_{ij}-a(i)\delta_{i}(j)\mbox{ for }i,j\in E.$$
We assume $\boldsymbol{A}$ is an irreducible matrix. It then follows by \cite[Lemma A.1]{BP} that $M(t)_{ij}>0$ for all $t>0$ and $i,j\in E$. Let $\Lambda:=\sup_{\lambda\in\sigma(\boldsymbol{A})}\mathrm{Re}(\lambda)$ where $\sigma(\boldsymbol{A})$ denotes the set of eigenvalues of $\boldsymbol{A}$. The Perron-Frobenius theory (see, for example, \cite[Lemma A.3]{BP}) tells us that for every $t>0$, $e^{\Lambda t}$ is a simple eigenvalue of $\boldsymbol{M(t)}$, and there exist a unique positive right eigenvector $\boldsymbol{u}=(u_{1},\cdots,u_{K})^{\mathrm{T}}$ and a unique positive left eigenvector $\boldsymbol{v}=(v_{1},\cdots,v_{K})^{\mathrm{T}}$ such that
\begin{equation}\nonumber
\sum_{i=1}^{K}u_{i}=\sum_{i=1}^{K}u_{i}v_{i}=1,\quad \boldsymbol{M(t)}\boldsymbol{u}=e^{\Lambda t}\boldsymbol{u},\quad \boldsymbol{v}^{T}\boldsymbol{M(t)}=e^{\Lambda t}\boldsymbol{v}.
\end{equation}
Moreover it is known by \cite[Lemma A.3]{BP} that for each $i,j\in E$,
\begin{equation}
e^{-\Lambda t}M(t)_{ij}\to u_{i}v_{j}\quad\mbox{ as }t\to +\infty.\label{eg1}
\end{equation}
One can easily verify that Assumptions 0-3 hold with $\lambda_{1}=-\Lambda$, $h(i)=cu_{i}$ and $\widehat{h}(i)=c^{-1}v_{i}$, where $c:=\left(\sum_{j=1}^{K}u^{2}_{j}\right)^{-1/2}$ is a positive constant. Thus $W^{h}_{t}(X):=c e^{-\Lambda t}\sum_{i=1}^{K}u_{i}X^{(i)}_{t}$ is a non-negative martingale. Applying Theorem \ref{them2} here, we can deduce that under the martingale change of measure the spine process $\widetilde{\xi}$ is a continuous-time Markov process on $E$ with
intensity matrix
$\boldsymbol{Q}=(q_{ij})_{ij}$ given by
$$q_{ii}:=-\frac{\gamma(i)\sum_{j=1}^{K}p_{ij}u_{j}}{u_{i}}=-(\Lambda+a(i)),\quad q_{ij}:=\frac{\gamma(i)p_{ij}u_{j}}{u_{i}}\quad\mbox{ for }i\not=j.$$
Let $\rho(dj):=u_{j}v_{j}m(dj)=
\sum_{i=1}^{K}
u_{j}v_{j}\delta_{i}(dj)$. Let $\widetilde{\mathfrak{S}}_{t}$ denote the transition semigroup of the spine $\widetilde{\xi}$ and $\widetilde{p}(t,i,j)$ denote its transition density with respect to $\rho$. It follows by Proposition \ref{prop2} that for each $i,j\in E$,
\begin{align}\nonumber
\widetilde{p}(t,i,j)u_{j}v_{j}&=\int_{E}\widetilde{p}(t,i,k)\delta_{j}(k)\rho({\rm d}k)=\widetilde{\mathfrak{S}}_{t}\delta_{j}(i)\nonumber\\
&=\frac{e^{-\Lambda t}}{h(i)}
\mathfrak{P}_{t}(h\delta_{j})(i)
=\frac{e^{-\Lambda t}}{u_{i}}M(t)_{ij}u_{j}.
\end{align}
Thus $\widetilde{p}(t,i,j)=e^{-\Lambda t}(u_{i}v_{j})^{-1}M(t)_{ij}$. By \eqref{eg1}, we have for each $i,j\in E$.
\begin{equation}\nonumber
\widetilde{p}(t,i,j)\to 1\quad\mbox{ as }t\to +\infty.
\end{equation}
 Hence Assumption 4 also holds for this example. Applying Corollary \ref{cor8.3} here, we conclude that for every non-trivial $\mu\in\me$, the martingale limit $$W_{\infty}^{h}(X):=\lim_{t\to+\infty}W^{h}_{t}(X)=\lim_{t\to+\infty}ce^{-\Lambda t}\sum_{i=1}^{K}u_{i}X^{(i)}_{t}$$ is non-degenerate if and only if $\Lambda >0$ and
 \begin{multline}\notag
\sum_{i=1}^{K}u_{i}v_{i}\int_{(0,+\infty)}r\log^{+}(r u_{i})\Pi^{L}(i,{\rm d}r)\\
+\sum_{i=1}^{K}\sum_{j=1}^{K}p_{ij}u_{j}v_{i}\int_{(0,+\infty)}r\log^{+}(r\sum_{k=1}^{K}p_{ik}u_{k})\Pi^{NL}(i,{\rm d}r)<+\infty.
 \end{multline}
 Using elementary computation
 , one can reduce the above condition to
 \begin{equation}
 \int_{(0,+\infty)}r\log^{+}r\Pi^{L}(i,{\rm d}r)+\int_{(0,+\infty)}r\log^{+}r\Pi^{NL}(i,{\rm d}r)<+\infty\quad\forall i\in E.\label{eg.2}
 \end{equation}
 In particular, under condition \eqref{eg.2},
 $\p_{\mu}\left(\lim_{t\to+\infty}X^{(i)}_{t}=0\right)=1$ for every $i\in E$ and every non-trivial $\mu\in\me$ if and only if $\Lambda\le 0$.
This result is also proved in \cite[Theorem 6]{KP}.}
\end{example}

Now we give some other examples.

\begin{example}\label{e:hering}{\rm
Suppose that $E$ is a bounded $C^3$ domain in
$\R^{d}\ (d\ge 1)$
, $m$ is the Lebesgue
measure on $E$ and that $\xi=(\xi_t, \Pi_x)$
is the killed Brownian motion in $E$. Suppose that $\phi^L$ and $\phi^{NL}$ are as
given in Subsection \ref{sec1}.
We assume Assumption 0 holds.
We further assume that the probability kernel
$\pi(x, {\rm d}y)$ has a bounded density with respect to the Lebesgue measure $m$, i.\/e.,
$\pi(x, {\rm d}y)=\pi(x, y){\rm d}y$ with $\pi(x, y)$ being bounded on $E\times E$. Assumption
1  and Assumption 3.(i1)
are trivially satisfied.
Let $(\mathfrak{P}_{t})_{t\ge 0}$ be the semigroup on $\mathcal{B}_{b}(E)$ uniquely determined by the integral equation \eqref{li2}.
It follows from \cite[Theorem]{Hering}
that Assumption 2, Assumption 3.(ii) are satisfied, and that $(\mathfrak{P}_{t})_{t\ge 0}$ is {\it uniformly primitive} in the sense of \cite{Hering}.
Thus for all $t>0$, $f\in \mathcal{B}^{+}_{b}(E)$ and $x\in E$,
\begin{equation}\label{unipri}
\left|\mathfrak{P}_{t}f(x)-e^{-\lambda_{1}t}(f,\widehat{h})h(x)\right|\le c_{t}e^{-\lambda_{1}t}(f,\widehat{h})h(x),
\end{equation}
where $c_{t}\ge 0$ satisfying $c_{t}\downarrow 0$ as $t\uparrow +\infty$, $\lambda_{1}$ is the constant in Assumption 2, and $h$, $\widehat{h}$ are the functions in Assumption 2. Let $\widetilde{\mathfrak{S}}_{t}f(x):= e^{\lambda_{1}t}h(x)^{-1}\mathfrak{P}_{t}(fh)(x)$ for $f\in\mathcal{B}^{+}(E)$, $t\ge 0$ and $x\in E$. Let $\widetilde{p}(t,x,y)$ be the density of $\widetilde{\mathfrak{S}P}_{t}$ with respect to the measure $\rho({\rm d}y):=h(y)\widehat{h}(y){\rm d}y$ on $E$. By \eqref{unipri}, we have
for every $t>0$, $f\in \mathcal{B}^{+}_{b}(E)$ and $x\in E$,
\begin{equation}\nonumber
\left|\widetilde{\mathfrak{S}}_{t}f(x)-\langle f,\rho\rangle\right|=\left|\int_{E}\left(\widetilde{p}(t,x,y)-1\right)f(y)\rho({\rm d}y)\right|\le c_{t}\langle f,\rho\rangle.
\end{equation}
It follows from this that
$$\sup_{x\in E}\esup_{y\in E}\left|\widetilde{p}(t,x,y)-1\right|\le c_{t}\to 0\quad\mbox{ as }t\to +\infty.$$
Hence Assumption 4 is satisfied.
Assumption
3.(iii) will be satisfied if the function $\pi(x, y)$ satisfies
$$
\int_E\pi(x, y)h(y){\rm d}y
\le ch(x) \qquad
\forall x\in \{z\in E:\ \gamma(z)>0\}
$$
for some constant $c>0$, where $h$ is the function in Assumption 2
and $\gamma(z)$ is as given in Subsection \ref{sec1}.

}
\end{example}

\begin{example}\label{examp:CKS}{\rm
Suppose that $E$ is a bounded $C^{1, 1}$ open set in
$\R^{d}\ (d\ge 1)$, $m$ is the Lebesgue
measure on $E$, $\alpha\in (0, 2)$, $\beta\in [0, \alpha\wedge d)$ and that $\xi=(\xi_t, \Pi_x)$ is an $m$-symmetric Hunt process on $E$ satisfying the following conditions: (1) $\xi$ has a L\'evy system
$(N,t)$
where $N=N(x, {\rm d}y)$ is a kernel given by
$$
N(x, {\rm d}y)=\frac{C_1}{|x-y|^{d+\alpha}}\,{\rm d}y \qquad x, y\in E
$$
for some constant $C_1>0$.
(2)  $\xi$ admits a
jointly continuous transition density $p(t, x, y)$ with respect to the Lebesgue measure
and that there exists a constant $C_2>1$ such that
$$
C_2^{-1}q_{\alpha, \beta}(t, x, y)\le p(t, x, y)\le C_2q_{\alpha, \beta}(t, x, y)
\qquad
\forall
(t, x, y)\in (0, 1]\times E\times E,
$$
where
\begin{equation}\label{e:CKSassump}
q_{\alpha, \beta}(t, x, y)
=\left(1\wedge\frac{\delta_E(x)}{t^{1/\alpha}}\right)^\beta
\left(1\wedge\frac{\delta_E(y)}{t^{1/\alpha}}\right)^\beta
\left(t^{-d/\alpha}\wedge\frac{t}{|x-y|^{d+\alpha}}\right).
\end{equation}
Here $\delta_E(x)$ stands for
the Euclidean distance between $x$ and the boundary of $E$.
Suppose that $\phi^L$ and $\phi^{NL}$ are as
given in Subsection \ref{sec1}.
We assume Assumption 0 holds.
We further assume that the probability kernel
$\pi(x, {\rm d}y)$ has a density $\pi(x, y)$ with respect to the Lebesgue measure $m$
satisfying the condition
$$
\pi(x, y)\le C_3|x-y|^{\epsilon-d} \qquad \forall x, y\in E
$$
for some positive constants $C_3$ and $\epsilon$.
Define
$$F(x,y):=C^{-1}_{1}|x-y|^{d+\alpha}\gamma(x)\pi(x,y)\quad\forall x,y\in E.$$
One can show easily that Assumption 1 and Assumption 3.(i2) are satisfied.
Define
$$F^{*}(x, y):=\log\left(1+F(x,y)\right) \qquad\forall x, y\in E.$$
It is obvious that there exists
$C_{4}>0$ such that
$$
0\le F^{*}(x, y)\le
C_{4}
\left(|x-y|^{\epsilon+\alpha}
\wedge 1\right) \qquad
\forall
x, y\in E,
$$
and thus, by \cite[Proposition 4.2]{CKS},
$F^{*}$ belongs to the Kato class ${\bf J}_{\alpha, \beta}$ defined in
\cite{CKS}, i.e., $\lim_{t\downarrow0}N^{\alpha, \beta}_{F^*}(t)=0$, where
$$
N^{\alpha, \beta}_{F^*}(t):=\sup_{x\in E}\int^t_0\int_{E\times E}q_{\alpha, \beta}(s, x, z)\left(1+\frac{|z-y|\wedge t^{1/\alpha}}{|x-z|}\right)^\beta\frac{F^*(y, z)+F^*(z, y)}{|z-y|^{d+\alpha}}dydzds.
$$
The measure $\mu({\rm d}x):=-a(x){\rm d}x$ obviously
belongs to
the Kato class ${\bf K}_{\alpha, \beta}$ defined
in \cite{CKS}, i.e., $\lim_{t\downarrow0}N^{\alpha, \beta}_\mu(t)=0$, where
$$
N^{\alpha, \beta}_\mu(t)=\sup_{x\in E}\int^t_0\int_Eq_{\alpha, \beta}(s, x, y)|a|(y)dyds,
$$
since $a$ is a bounded function.
For $0\le t<+\infty$, let $A_{t}:=-\int_{0}^{t}a(\xi_{r}){\rm d}r+\sum_{0<r \le t}F^{*}(\xi_{r-},\xi_{r})$.
Let $(T_{t})_{t\ge 0}$ be the Feynman-Kac semigroup of $\xi$ given by
$$T_{t}f(x):=\Pi_{x}\left[\exp\left(A_{t}\right)f(\xi_{t})\right],\qquad t\ge 0,\ x\in E,
 f\in\mathcal{B}^{+}(E).$$
Now it follows from
\cite[Theorem 1.3]{CKS}
that the semigroup
$(T_{t})_{t\ge 0}$
has a jointly continuous density $q(t,x,y)$ with respect to the Lebesgue measure
and there exists a constant
$C_{5}>1$
such that
\begin{equation}\label{e:CKS}
C^{-1}_{5}
q_{\alpha, \beta}(t, x, y)\le q(t, x, y)\le C_{5}q_{\alpha, \beta}(t, x, y)
\qquad
\forall
(t, x, y)\in (0, 1]\times E\times E.
\end{equation}
Let $(\widehat{T}_t)_{t\ge 0}$ be the dual semigroup of $(T_t)_{t\ge 0}$.
By \eqref{e:CKS}, one can easily show that for any $f\in\mathcal{B}_{b}(E)$, $T_{t}f$ and $\widehat{T}_{t}f$ are bounded continuous functions on $E$, that $T_{t}$ and $\widehat{T}_{t}$ are bounded operators from $L^{2}(E,m)$ into $L^{\infty}(E,m)$, and that $(T_{t})_{t\ge 0}$ and $(\widehat{T}_{t})_{t\ge 0}$ are strongly continuous semigroups on $L^{2}(E,m)$.
Let
${\bf L}$ and ${\widehat{\bf L}}$ be the generators of $(T_t)_{t\ge 0}$
and $(\widehat{T}_t)_{t\ge 0}$ respectively.
Let $\sigma(\mathbf{L})$ and $\sigma(\widehat{\mathbf{L}})$ denote the spectrum of $\mathbf{L}$ and $\widehat{\mathbf{L}}$ respectively.
It follows from \eqref{e:CKS} and Jentzsch's theorem (\cite[Theorem V.6.6, p. 337]{Sch})
that the common value $
-\lambda_{1}
:=\sup{\rm Re} (\sigma({\bf L}))=
\sup{\rm Re} (\sigma(\widehat{{\bf L}}))$ is an eigenvalue of multiplicity 1 for both
${\bf L}$ and ${\widehat{\bf L}}$, and that an eigenfunction $h$ of ${\bf L}$
associated with
$-\lambda_{1}$
is bounded continuous and
can be chosen strictly positive on $E$ and satisfies $\|h\|_{L^2(E, m)}=1$,
and that an eigenfunction $\widehat{h}$ of ${\widehat{\bf L}}$
associated with
$-\lambda_{1}$
is bounded continuous and
can be chosen strictly positive on $E$ and satisfies $(h, \widehat{h})=1$.
Thus Assumption 2 and 3.(ii) are satisfied.
It follows from \eqref{e:CKS} and the equations $e^{-\lambda_{1}}h=T_{1}h$, $e^{-\lambda_{1}}\widehat{h}=\widehat{T}_{1}h$
that there exists a constant
$C_{6}>1$
such that
$$
C^{-1}_{6}
\delta_E(x)^{\beta}\le h(x)\le
C_{6}
\delta_E(x)^{\beta},\quad
C^{-1}_{6}
\delta_E(x)^{\beta}\le \widehat{h}(x)\le
C_{6}
\delta_E(x)^{\beta}\qquad\forall x\in E.$$
It follows from this, \eqref{e:CKS} and the semigroup properpty that
the semigroups $(T_t)_{t\ge 0}$ and $(\widehat{T}_t)_{t\ge 0}$ are intrinsically ultracontractive.
For the definition of intrinsic ultracontractivity, see \cite{KS}.
Let $\widetilde{\mathfrak{S}}_{t}f(x):=e^{\lambda_{1}t}h(x)^{-1}T_{t}(fh)(x)$ for $f\in\mathcal{B}^{+}(E)$, $t\ge 0$ and $x\in E$. Then $\widetilde{\mathfrak{S}}_{t}$ admits a density $\widetilde{p}(t,x,y)$ with respect to the probability measure $h(y)\widehat{h}(y){\rm d}y$ which is related to $q(t,x,y)$ by
$$\widetilde{p}(t,x,y)=\frac{e^{\lambda_{1}t}q(t,x,y)}{h(x)\widehat{h}(y)}\qquad\forall (t,x,y)\in (0,+\infty)\times E\times E.$$
Now it follows from
\cite[Theorem 2.7]{KS} that Assumption 4 is satisfied.  As in the previous example,
Assumption
3.(iii) will be satisfied if the function $\pi(x, y)$ satisfies
$$
\int_E\pi(x, y)h(y){\rm d}y
\le ch(x) \qquad
\forall x\in \{z\in E:\gamma(z)>0\}
$$
for some constant $c>0$, where $\gamma(z)$ is as given in Subsection \ref{sec1}.

One concrete example of $\xi$ is the killed symmetric $\alpha$-stable process
in $E$. In this case, \eqref{e:CKSassump} is satisfied with $\beta=\alpha/2$, a fact
which was first proved in \cite{CKSjems}.

Another concrete example of $\xi$
is the censored symmetric $\alpha$-stable process
in $E$ introduced in \cite{BBC} when $\alpha\in (1, 2)$. In this case, \eqref{e:CKSassump} is satisfied with $\beta=\alpha-1$, a fact
which was first proved in \cite{CKSptrf}.

In fact, by using \cite{CKS}, one could also include the case when $E$ is a $d$-set,
$\alpha\in (0, 2)$ and $\xi$ is an $\alpha$-stable-like process in $E$ introduced in
\cite{CK}. We omit the details.
}
\end{example}

\begin{example}
{\rm
Suppose that $E=\mathbb{R}^{d}$, $m$ is the Lebesgue measure on $\mathbb{R}^{d}$, $\alpha\in (0,2)$ and that $\xi=(\xi_t,\Pi_{x})$ is a Markov process corresponding to the Feynman-Kac transform
of a $d$-dimensional isotropic $\alpha$-stable process with killing potential
$\eta(x)=|x|^\beta$ ($\beta>0$).
Let $J(x)=J(|x|)$ be the L\'evy density of the isotropic $\alpha$-stable process, i.e., $J(x)=c(d,\alpha)|x|^{-d-\alpha}$ for some positive constant $c(d,\alpha)$ depending only on $d$ and $\alpha$.
It is known that $\xi$ has a L\'{e}vy system $(N,t)$ where $N(x,\mathrm{d}y)=2J(y-x)\mathrm{d}y$.
Let $(\mathcal{E},\mathcal{F})$ be the Dirichlet form of $\xi$. Then $\mathcal{E}$ has the following form
$$\mathcal{E}(u,v)=\int_{\R^d}\int_{\R^d}(u(x)-u(y))(v(x)-v(y))J(y-x)\mathrm{d}x\mathrm{d}y+\int_{\R^d}u(x)v(x)|x|^\beta \mathrm{d}x\quad \forall u,v\in\mathcal{F}.$$
Suppose that the branching mechanisms $\phi^{L}$ and $\phi^{NL}$ are as given in Subsection 2.1. For simplicity we assume $a(x)\equiv 0$ and Assumption 0 holds.
Let $\pi(x)=\pi(|x|)$ be a
probability density
on $\mathbb{R}^{d}$ such that the function $\pi(x)/J(x)$ is bounded from above.
We assume that the probability kernel $\pi(x,\mathrm{d}y)$ has a density
$\pi(x,y)=\pi(y-x)$ with respect to the Lebesgue measure and that the function $\gamma(x)\equiv \gamma$ is a constant.
Define $\gamma(x, y):=\gamma\pi(y-x)$. Then $\int_{\mathbb{R}^{d}}\gamma(x,y)\mathrm{d}x=\gamma\int_{\mathbb{R}^{d}}\pi(y-x)\mathrm{d}x=\gamma$, and Assumption 1 is trivially satisfied. Define
$$
F(x, y):=\frac{\gamma(x,y)}{2J(y-x)}=\frac{\gamma\pi(y-x)}{2J(y-x)}\quad\forall x,y\in\R^{d}.
$$
Then $F$ is a bounded function on $\R^{d}\times\R^{d}$ vanishing on the diagonal, and thus Assumption 3.(i2) is satisfied.
Define
\begin{align*}
\widetilde{A}_t&:=\sum_{s\le t}\log(1+F(\xi_{s-}, \xi_s))-2\int^t_0\int_{\R^d}F(\xi_s, y)J(y-\xi_{s})\mathrm{d}y\mathrm{d}s\\
&=\sum_{s\le t}\log(1+F(\xi_{s-}, \xi_s))-\gamma\int^t_0\int_{\R^d}\pi(y-\xi_{s})\mathrm{d}y\mathrm{d}s\\
&=\sum_{s\le t}\log(1+F(\xi_{s-}, \xi_s))-\gamma t.
\end{align*}
It follows from \cite[Theorem 4.8]{CS} (see also \cite[p. 275]{CS2}) that
the bilinear form corresponding to the symmetric semigroup
$$
\widetilde{T}_tf(x)=\Pi_{x}\left[\mathrm{e}^{\widetilde{A}_t}f(\xi_t)\right]\quad\forall t\ge 0,\ x\in\R^{d},\ f\in\mathcal{B}_{b}(\R^{d})
$$
is
\begin{align*}
\widetilde{\mathcal{Q}}(u, v)&=\int_{\R^d}\int_{\R^d}(u(x)-u(y))(v(x)-v(y))(1+
F(x, y))J(y-x)\mathrm{d}xdy\\
&+\int_{\R^d}u(x)v(x)|x|^\beta \mathrm{d}x\quad\forall u,v\in\mathcal{F}.
\end{align*}
This is the bilinear form of the stable-like L\'evy process killed with the potential
$\eta(x)=|x|^{\beta}$. Now we apply
\cite[Examples 4.5 and 4.8]{KL} to get that the symmetric semigroup $(\widetilde{T}_t)_{t\ge 0}$ is intrinsically ultracontractive.
Define
$$
A_t:=\sum_{s\le t}\log(1+F(\xi_{s-}, \xi_s)).
$$
Again by \cite[Theorem 4.8]{CS} (see also \cite[p. 275]{CS2}) that
the bilinear form corresponding to the symmetric semigroup
\begin{equation*}
T_tf(x)=\Pi_x\left[e^{A_t}f(\xi_t)\right]\quad\forall t\ge 0,\ x\in\R^{d},\ f\in\mathcal{B}_{b}(\R^{d})
\end{equation*}
is
\begin{align*}
\mathcal{Q}(u, v)&=\mathcal{E}(u,v)-2\int_{\R^d}\int_{\R^d}u(y)v(x)F(x,y)J(x-y)\mathrm{d}x\mathrm{d}y\\
&=\mathcal{E}(u,v)-\int_{\R^{d}}\int_{\R^{d}}u(y)v(x)\gamma(x,y)\mathrm{d}x\mathrm{d}y\quad\forall u,v\in\mathcal{F}.
\end{align*}
We observe that $T_{t}f=\mathrm{e}^{\gamma t}\widetilde{T}_{t}f$. So the semigroup $(T_{t})_{t\ge 0}$ is also intrinsically ultracontractive.
Applying similar argument as in Example \ref{examp:CKS}, one can show that Assumptions 2 and 4 are satisfied.
Finally, Assumption 3(iii) (and thus Assumption 3(ii))
will be satisfied if the function $\pi(x)$ satisfies that
$$\int_{\R^{d}}\pi(y-x)h(y)\mathrm{d}y\le c h(x)\quad \forall x\in\R^{d}$$
for some positive constant $c$. Here $h$ is the strictly positive eigenfunction associated with the principal eigenvalue of the generator of the semigroup $(T_{t})_{t\ge 0}$.
}
\end{example}

\section*{Appendix}
\noindent\textit{Proof of Proposition \ref{prop3}: }
We will prove \eqref{l1} first. We claim that
\begin{align}
\widetilde{\Pi}_{x}\left[f(\tau_{1},\widetilde{\xi}_{\tau_{1}-},\widetilde{\xi}_{\tau_{1}})1_{\{\tau_{1}\le t\}}\right]
&=\widetilde{\Pi}_{x}\left[\int_{0}^{t\wedge \tau_{1}}q(\widetilde{\xi}_{s}){\rm d}s\int_{E}f(s,\widetilde{\xi}_{s},y)\pi^{h}(\widetilde{\xi}_{s},{\rm d}y)\right].\label{prop3.1}
\end{align}
It is easy to see from the construction of $\widetilde{\xi}$ that
\begin{equation}
\mbox{LHS of }\eqref{prop3.1}=\Pi^{h}_{x}\left[\int_{0}^{t}q(\xi^{h}_{s})e_{q}(s){\rm d}s\int_{E}f(s,\xi^{h}_{s},y)\pi^{h}(\xi^{h}_{s},{\rm d}y)\right].\label{prop3.7}
\end{equation}
On the other hand, by Fubini's theorem, we have
\begin{align}
\mbox{RHS of }\eqref{prop3.1}&=\int_{0}^{t}{\rm d}s\,\widetilde{\Pi}_{x}\left[\int_{E}q(\widetilde{\xi}_{s})f(s,\widetilde{\xi}_{s},y)
\pi^{h}(\widetilde{\xi}_{s},{\rm d}y)1_{\{s<\tau_{1}\}}\right]\nonumber\\
&=\int_{0}^{t}{\rm d}s\,\Pi^{h}_{x}\left[\int_{E}q(\widehat{\xi}_{s})f(s,\widehat{\xi}_{s},y)
\pi^{h}(\widehat{\xi}_{s},{\rm d}y)1_{\{s<\widehat{\zeta}\}}\right]\nonumber\\
&=\int_{0}^{t}{\rm d}s\,\Pi^{h}_{x}\left[e_{q}(s)\int_{E}q(\xi^{h}_{s})f(s,\xi^{h}_{s},y)\pi^{h}(\xi^{h}_{s},{\rm d}y)\right]\nonumber\\
&=\Pi^{h}_{x}\left[\int_{0}^{t}q(\xi^{h}_{s})e_{q}(s){\rm d}s\int_{E}f(s,\xi^{h}_{s},y)\pi^{h}(\xi^{h}_{s},{\rm d}y)\right].\label{prop3.8}
\end{align}
Combining \eqref{prop3.7} and \eqref{prop3.8} we arrive at the claim \eqref{prop3.1}.
Note that applying the shift operator $\widetilde{\theta}_{\tau_{n}}$ to $f(\tau_{1},\widetilde{\xi}_{\tau_{1}-},\widetilde{\xi}_{\tau_{1}})1_{\{\tau_{1}\le t\}}$ gives $f(\tau_{n+1},\widetilde{\xi}_{\tau_{n+1}-},\widetilde{\xi}_{\tau_{n+1}})1_{\{\tau_{n+1}\le t\}}$. Using the strong Markov property of $\widetilde{\xi}$ and Fubini's theorem, we can prove by induction that for all
$n\ge 2$,
\begin{equation}
\widetilde{\Pi}_{x}\left[f(\tau_{n},\widetilde{\xi}_{\tau_{n}-},\widetilde{\xi}_{\tau_{n}})1_{\{\tau_{n}\le t\}}\right]
=\widetilde{\Pi}_{x}\left[\int_{t\wedge \tau_{n-1}}^{t\wedge \tau_{n}}q(\widetilde{\xi}_{s}){\rm d}s\int_{E}f(s,\widetilde{\xi}_{s},y)\pi^{h}(\widetilde{\xi}_{s},{\rm d}y)\right].\nonumber
\end{equation}
Thus by the above equality, Fubini's theorem and the fact that $\widetilde{\Pi}_{x}(\lim_{n\to+\infty}\tau_{n}=+\infty)=1$, we have
\begin{align}
\widetilde{\Pi}_{x}\left[\sum_{\tau_{i}\le t}f(\tau_{i},\widetilde{\xi}_{\tau_{i}-},\widetilde{\xi}_{\tau_{i}})\right]
&=\widetilde{\Pi}_{x}\left[\sum_{i=1}^{+\infty}f(\tau_{i},\widetilde{\xi}_{\tau_{i}-},\widetilde{\xi}_{\tau_{i}})1_{\{\tau_{i}\le t\}}\right]\nonumber\\
&=\widetilde{\Pi}_{x}\left[\lim_{n\to+\infty}\int_{0}^{t\wedge \tau_{n}}q(\widetilde{\xi}_{s}){\rm d}s\int_{E}f(s,\widetilde{\xi}_{s},y)\pi^{h}(\widetilde{\xi}_{s},{\rm d}y)\right]\nonumber\\
&=\widetilde{\Pi}_{x}\left[\int_{0}^{t}q(\widetilde{\xi}_{s}){\rm d}s\int_{E}f(s,\widetilde{\xi}_{s},y)\pi^{h}(\widetilde{\xi}_{s},{\rm d}y)\right].\nonumber
\end{align}
Hence we have proved  \eqref{l1}.
We next show \eqref{l2}. It is easy to see that
\begin{multline}
\widetilde{\Pi}_{x}\left[\left(\sum_{\tau_{i}\le t}f(\tau_{i},\widetilde{\xi}_{\tau_{i}-},\widetilde{\xi}_{\tau_{i}})\right)\left(\sum_{\tau_{j}\le t}g(\tau_{j},\widetilde{\xi}_{\tau_{j}-},\widetilde{\xi}_{\tau_{j}})\right)\right]\\
\shoveleft{=\widetilde{\Pi}_{x}\left[\sum_{\tau_{i}\le t}f g(\tau_{i},\widetilde{\xi}_{\tau_{i}-},\widetilde{\xi}_{\tau_{i}})\right]}\\
+\sum_{i=1}^{+\infty}\sum_{j=i+1}^{+\infty}
\left\{\widetilde{\Pi}_{x}\left[f(\tau_{i},\widetilde{\xi}_{\tau_{i}-},\widetilde{\xi}_{\tau_{i}})
g(\tau_{j},\widetilde{\xi}_{\tau_{j}-},\widetilde{\xi}_{\tau_{j}})1_{\{\tau_{j}\le t\}}\right]\right.\\
+\left.\widetilde{\Pi}_{x}\left[g(\tau_{i},\widetilde{\xi}_{\tau_{i}-},\widetilde{\xi}_{\tau_{i}})
f(\tau_{j},\widetilde{\xi}_{\tau_{j}-},\widetilde{\xi}_{\tau_{j}})1_{\{\tau_{j}\le t\}}\right]\right\}.\label{prop3.9}
\end{multline}
By the strong Markov property and \eqref{prop3.7}, we have for $j\ge 2$,
\begin{align}
&\widetilde{\Pi}_{x}\left[f(\tau_{1},\widetilde{\xi}_{\tau_{1}-},\widetilde{\xi}_{\tau_{1}})
g(\tau_{j},\widetilde{\xi}_{\tau_{j}-},\widetilde{\xi}_{\tau_{j}})1_{\{\tau_{j}\le t\}}\right]\nonumber\\
&=\widetilde{\Pi}_{x}\left[f(\tau_{1},\widetilde{\xi}_{\tau_{1}-},\widetilde{\xi}_{\tau_{1}})1_{\{\tau_{1}\le t\}}\left.\widetilde{\Pi}_{\widetilde{\xi}_{\tau_{1}}}\left(g(\tau_{j-1}+s,\widetilde{\xi}_{\tau_{j-1}-},\widetilde{\xi}_{\tau_{j-1}})1_{\{\tau_{j-1}+s\le t\}}\right)\right|_{s=\tau_{1}}\right]\nonumber\\
&=\Pi^{h}_{x}\big[\int_{0}^{t}{\rm d}s\int_{E}\pi^{h}(\xi^{h}_{s},{\rm d}y)q(\xi^{h}_{s})e_{q}(s)f(s,\xi^{h}_{s},y)
\widetilde{\Pi}_{y}
\left(g(\tau_{j-1}+s,\widetilde{\xi}_{\tau_{j-1}-},\widetilde{\xi}_{\tau_{j-1}})1_{\{\tau_{j-1}\le t-s\}}\right)
\big],
\label{prop3.3}
\end{align}
and for $j>i\ge 2$,
\begin{align}
&\widetilde{\Pi}_{x}\left[f(\tau_{i},\widetilde{\xi}_{\tau_{i}-},\widetilde{\xi}_{\tau_{i}})
g(\tau_{j},\widetilde{\xi}_{\tau_{j}-},\widetilde{\xi}_{\tau_{j}})1_{\{\tau_{j}\le t\}}\right]\nonumber\\
&=\widetilde{\Pi}_{x}\left[1_{\{\tau_{1}\le t\}}\left.\widetilde{\Pi}_{\widetilde{\xi}_{\tau_{1}}}\left(f(\tau_{i-1}+s,\widetilde{\xi}_{\tau_{i-1}-},\widetilde{\xi}_{\tau_{i-1}})g(\tau_{j-1}+s,\widetilde{\xi}_{\tau_{j-1}-},\widetilde{\xi}_{\tau_{j-1}})1_{\{\tau_{j-1}+s\le t\}}\right)\right|_{s=\tau_{1}}\right]\nonumber\\
&=\Pi^{h}_{x}\big[\int_{0}^{t}{\rm d}s\int_{E}\pi^{h}(\xi^{h}_{s},{\rm d}y)  q(\xi^{h}_{s})e_{q}(s)\nonumber\\
&\quad\quad \cdot\widetilde{\Pi}_{y}
\left(f(\tau_{i-1}+s,\widetilde{\xi}_{\tau_{i-1}-},\widetilde{\xi}_{\tau_{i-1}})g(\tau_{j-1}+s,\widetilde{\xi}_{\tau_{j-1}-},\widetilde{\xi}_{\tau_{j-1}})1_{\{\tau_{j-1}\le t-s\}}\right)\big].\label{prop3.4}
\end{align}
By \eqref{prop3.3}, Fubini's theorem, the strong Markov property of $\widetilde{\xi}$, \eqref{l1} and \eqref{prop3.7},
\begin{align}
&\sum_{j=2}^{+\infty}\widetilde{\Pi}_{x}\left[f(\tau_{1},\widetilde{\xi}_{\tau_{1}-},\widetilde{\xi}_{\tau_{1}})
g(\tau_{j},\widetilde{\xi}_{\tau_{j}-},\widetilde{\xi}_{\tau_{j}})1_{\{\tau_{j}\le t\}}\right]\nonumber\\
&=\Pi^{h}_{x}\big[\int_{0}^{t}{\rm d}s\int_{E}\pi^{h}(\xi^{h}_{s},{\rm d}y)q(\xi^{h}_{s})e_{q}(s)f(s,\xi^{h}_{s},y)\notag\\
&\quad\quad\quad\quad\cdot\widetilde{\Pi}_{y}\left(\sum_{j=2}^{+\infty}g(\tau_{j-1}+s,\widetilde{\xi}_{\tau_{j-1}-},\widetilde{\xi}_{\tau_{j-1}})1_{\{\tau_{j-1}\le t-s\}}\right)\big],\nonumber\\
&=\Pi^{h}_{x}\left[\int_{0}^{t}{\rm d}s\int_{E}\pi^{h}(\xi^{h}_{s},{\rm d}y)q(\xi^{h}_{s})e_{q}(s)f(s,\xi^{h}_{s},y)\widetilde{\Pi}_{y}
\left(\sum_{\tau_{k}\le t-s}g(\tau_{k}+s,\widetilde{\xi}_{\tau_{k}-},\widetilde{\xi}_{\tau_{k}})\right)\right]\nonumber\\
&=\Pi^{h}_{x}\big[\int_{0}^{t}{\rm d}s\int_{E}\pi^{h}(\xi^{h}_{s},{\rm d}y)q(\xi^{h}_{s})e_{q}(s)f(s,\xi^{h}_{s},y)\notag\\
&\quad\quad\quad\quad\cdot\widetilde{\Pi}_{y}
\left(\int_{0}^{t-s}{\rm d}r\int_{E}\pi^{h}(\widetilde{\xi}_{r},{\rm d}z)q(\widetilde{\xi}_{r})g(r+s,\widetilde{\xi}_{r},z)\right)\big]\nonumber\\
&=\widetilde{\Pi}_{x}\left[f(\tau_{1},\widetilde{\xi}_{\tau_{1}-},\widetilde{\xi}_{\tau_{1}})\left.\widetilde{\Pi}_{\widetilde{\xi}_{\tau_{1}}}
\left(\int_{0}^{t-s}{\rm d}r\int_{E}\pi^{h}(\widetilde{\xi}_{r},{\rm d}z)q(\widetilde{\xi}_{r})g(r+s,\widetilde{\xi}_{r},z)\right)\right|_{s=\tau_{1}}
1_{\{\tau_{1}\le t\}}\right]\nonumber\\
&=\widetilde{\Pi}_{x}\left[f(\tau_{1},\widetilde{\xi}_{\tau_{1}-},\widetilde{\xi}_{\tau_{1}})
\int_{\tau_{1}}^{t}{\rm d}r\int_{E}\pi^{h}(\widetilde{\xi}_{r},{\rm d}z)q(\widetilde{\xi}_{r})g(r,\widetilde{\xi}_{r},z)
1_{\{\tau_{1}\le t\}}\right].\nonumber
\end{align}
Similarly, by \eqref{prop3.4}, Fubini's theorem, the strong Markov property of $\widetilde{\xi}$, \eqref{l1} and \eqref{prop3.7}, we can prove by induction that for $i\ge 1$,
\begin{align}
&\sum_{j=i+1}^{+\infty}\widetilde{\Pi}_{x}\left[f(\tau_{i},\widetilde{\xi}_{\tau_{i}-},\widetilde{\xi}_{\tau_{i}})
g(\tau_{j},\widetilde{\xi}_{\tau_{j}-},\widetilde{\xi}_{\tau_{j}})1_{\{\tau_{j}\le t\}}\right]\nonumber\\
&=\widetilde{\Pi}_{x}\left[f(\tau_{i},\widetilde{\xi}_{\tau_{i}-},\widetilde{\xi}_{\tau_{i}})
\int_{\tau_{i}}^{t}{\rm d}r\int_{E}\pi^{h}(\widetilde{\xi}_{r},{\rm d}z)q(\widetilde{\xi}_{r})g(r,\widetilde{\xi}_{r},z)
1_{\{\tau_{i}\le t\}}\right].\nonumber
\end{align}
By this, Fubini's theorem, the strong Markov property of $\widetilde{\xi}$ and \eqref{l1}, we get
\begin{align}
&\sum_{i=1}^{+\infty}\sum_{j=i+1}^{+\infty}\widetilde{\Pi}_{x}\left[f(\tau_{i},\widetilde{\xi}_{\tau_{i}-},\widetilde{\xi}_{\tau_{i}})
g(\tau_{j},\widetilde{\xi}_{\tau_{j}-},\widetilde{\xi}_{\tau_{j}})1_{\{\tau_{j}\le t\}}\right]\nonumber\\
&=\widetilde{\Pi}_{x}\left[\sum_{i=1}^{+\infty}f(\tau_{i},\widetilde{\xi}_{\tau_{i}-},\widetilde{\xi}_{\tau_{i}})
\int_{\tau_{i}}^{t}{\rm d}r\int_{E}\pi^{h}(\widetilde{\xi}_{r},{\rm d}z)q(\widetilde{\xi}_{r})g(r,\widetilde{\xi}_{r},z)
1_{\{\tau_{i}\le t\}}\right]\nonumber\\
&=\widetilde{\Pi}_{x}\left[\sum_{\tau_{i}\le t}f(\tau_{i},\widetilde{\xi}_{\tau_{i}-},\widetilde{\xi}_{\tau_{i}})
\int_{\tau_{i}}^{t}{\rm d}r\int_{E}\pi^{h}(\widetilde{\xi}_{r},{\rm d}z)q(\widetilde{\xi}_{r})g(r,\widetilde{\xi}_{r},z)\right]\nonumber\\
&=\widetilde{\Pi}_{x}\left[\sum_{\tau_{i}\le t}f(\tau_{i},\widetilde{\xi}_{\tau_{i}-},\widetilde{\xi}_{\tau_{i}})
\left.\widetilde{\Pi}_{\widetilde{\xi}_{\tau_{i}}}\left(\int_{0}^{t-s}q(\widetilde{\xi}_{r}){\rm d}r\int_{E}g(s+r,\widetilde{\xi}_{r},z)
\pi^{h}(\widetilde{\xi}_{r},{\rm d}z)\right)\right|_{s=\tau_{i}}\right]\nonumber\\
&=\widetilde{\Pi}_{x}\big[\int_{0}^{t}{\rm d}s\int_{E}\pi^{h}(\widetilde{\xi}_{s},{\rm d}y)q(\widetilde{\xi}_{s})f(s,\widetilde{\xi}_{s},y)
\widetilde{\Pi}_{y}\left(\int_{0}^{t-s}{\rm d}r\int_{E}\pi^{h}(\widetilde{\xi}_{r},{\rm d}z)q(\widetilde{\xi}_{r})g(s+r,\widetilde{\xi}_{r},z)\right)
\big].\label{prop3.6}
\end{align}
Combining \eqref{prop3.9} and \eqref{prop3.6}, we arrive at
\eqref{l2}.
\qed
\medskip

\section*{Acknowledgements}
 The research of Y.-X. Ren is supported by NSFC
(Grant Nos.  11671017 and 11731009), and  LMEQF.
The research of R. Song is supported in part by a grant from the Simons
Foundation (\#429343).
The research of T. Yang is supported by NSFC
(Grant Nos. 11501029 and 11731009).
The authors thank Professor Zhen-Qing Chen for his helpful comments on the construction of
the concatenation process.

\end{document}